\documentclass[lang = american]{ems-icm} 

\usepackage{enumerate}
\usepackage{undertilde}
\usepackage[normalem]{ulem}
\usepackage{subfigure}
\usepackage{enumitem}
\usepackage[11pt]{moresize}

\usepackage{xurl}

\numberwithin{equation}{section}

\begin{document}

\volume{?} 

\title{What is a random surface?}
\titlemark{Random surface}

\emsauthor{1}{Scott Sheffield}{S.~Sheffield}

\emsaffil{1}{77 Mass.\ Ave.\ Cambridge, MA, USA \email{sheffield@math.mit.edu}}

\begin{abstract}
Given $2n$ unit equilateral triangles, there are finitely many ways to glue each edge to a partner.  We obtain a random {\em sphere-homeomorphic} surface by sampling uniformly from the gluings that produce a topological sphere. As $n \to \infty$ these random surfaces (appropriately scaled) converge in law. The limit is a ``canonical'' sphere-homeomorphic random surface, much the way Brownian motion is a canonical random path.
\vspace{.1in}

Depending on how the surface space and convergence topology are specified, the limit is the {\em Brownian sphere}, the {\em peanosphere}, the {\em pure Liouville quantum gravity sphere}, or a certain {\em conformal field theory}. All of these objects have concise definitions, and are all in some sense equivalent, but the equivalence is highly non-trivial, building on hundreds of math and physics papers over the past half century.

\vspace{.1in}
More generally, the ``continuum random surface embedded in $d$-dimensional Euclidean space'' makes a kind of sense for $d \in (-\infty, 25)$ even when $d$ is not a positive integer; and this can be extended to higher genus surfaces, surfaces with boundary, and surfaces with marked points or other decoration.

\vspace{.1in}

These constructions have deep roots in both mathematics and physics, drawing from classical graph theory, complex analysis, probability and representation theory, as well as string theory, planar statistical physics, random matrix theory and a simple model for two-dimensional quantum gravity.
\vspace{.1in}

We present here an informal, colloquium-level overview of the subject, which we hope will be accessible to both newcomers and experts. We aim to answer, as cleanly as possible, the fundamental question. What is a random surface?

\end{abstract}

\maketitle

As a random path, Brownian motion is {\em canonical} in the sense that it is uniquely characterized by certain symmetries, and {\em universal} in the sense that it is a limit of many kinds of discrete random walks. This paper will describe a similarly canonical and universal {\em random surface}. This random surface has several different formulations. To tell the story in a fanciful way, imagine a dialog, first about random paths, then about random surfaces.

\vspace{.03in}
\noindent {\textbf INSTRUCTOR:} Consider the simple random walk on $\mathbb Z$. At each time step a coin toss decides whether position goes up or down. If you shrink the graph horizontally by a factor of $C$ and vertically by a factor of $\sqrt{C}$, then the $C \to \infty$ limit is {\em Brownian motion}.
\vspace{-.3in}

\begin{figure}[ht!]
\begin{center}
\includegraphics[scale=0.2]{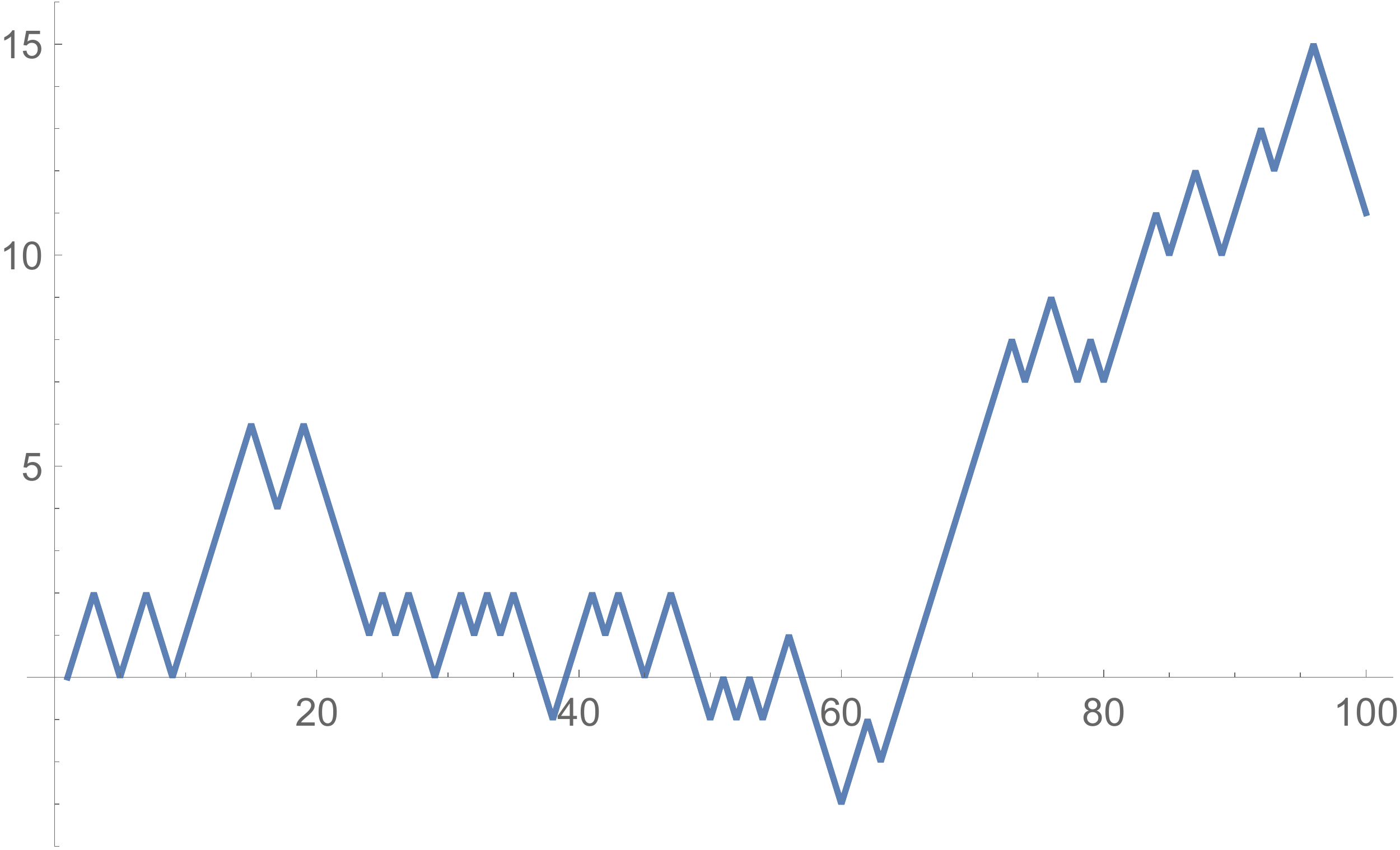}
\hspace{.4in} \includegraphics[scale=0.2]{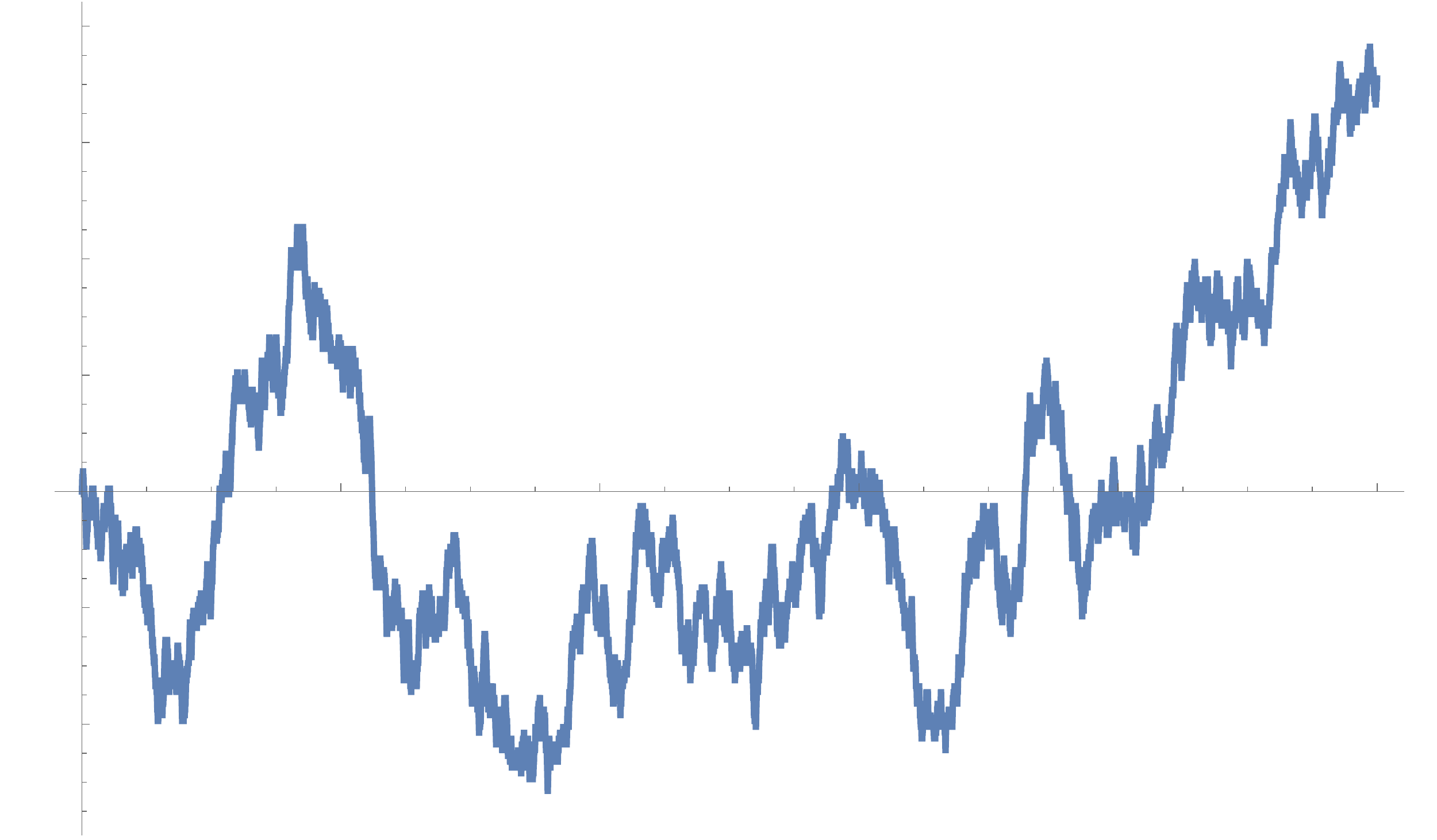}
\end{center}
\end{figure}

\vspace{-.6in}

\noindent {\textbf STUDENT:} Great! But can you define Brownian motion directly in the continuum?
\vspace{.03in}

\noindent {\textbf INSTRUCTOR:} Sure! Fix $0=t_0\! < \!t_1\! <\! \ldots\! <\! t_n$. Specify the joint law of $B(t_1), \ldots, B(t_n)$ by making increments $B(t_k)\! -\! B(t_{k-1})$ independent normals with mean $0$, variance $t_k \!-\! t_{k-1}$. Extend to countable dense set (Kolmogorov extension), then all $t$ (Kolomogorov-\v{C}entsov).
\vspace{.03in}

\noindent {\textbf STUDENT:} Are there other natural ways to characterize Brownian motion?
\vspace{.03in}

\noindent {\textbf INSTRUCTOR:} Brownian motion is {\em canonical} in that it is the only random path with certain symmetries (like stationarity/independence of increments). It is {\em universal} in that (per central limit theorem) it is a limit of many discrete walks. It comes up everywhere.
\vspace{.03in}

\noindent {\textbf STUDENT:} What if I want a random path embedded in $\mathbb R^d$?
\vspace{.03in}

\noindent
{\textbf INSTRUCTOR:} Use a vector $\Bigl(B_1(t), B_2(t), \ldots, B_d(t) \Bigr)$ of independent Brownian motions.
\vspace{.08in}

The student is happy. Now imagine a similar dialog for random surfaces.
\vspace{.08in}

\noindent {\textbf INSTRUCTOR:} Take a uniformly random triangulation of the sphere with $n$ triangles. The picture below is an example with 30,000 triangles by Budzinski, using some software to embed the surface in three dimensions and give us a view. The $n \to \infty$ limit of this object is a random fractal surface called the {\em Brownian sphere}. It is also a {\em peanosphere} and a {\em pure Liouville quantum gravity sphere} and a {\em conformal field theory}.

\vspace{-.2in}
\begin{center}
\includegraphics[scale=0.04]{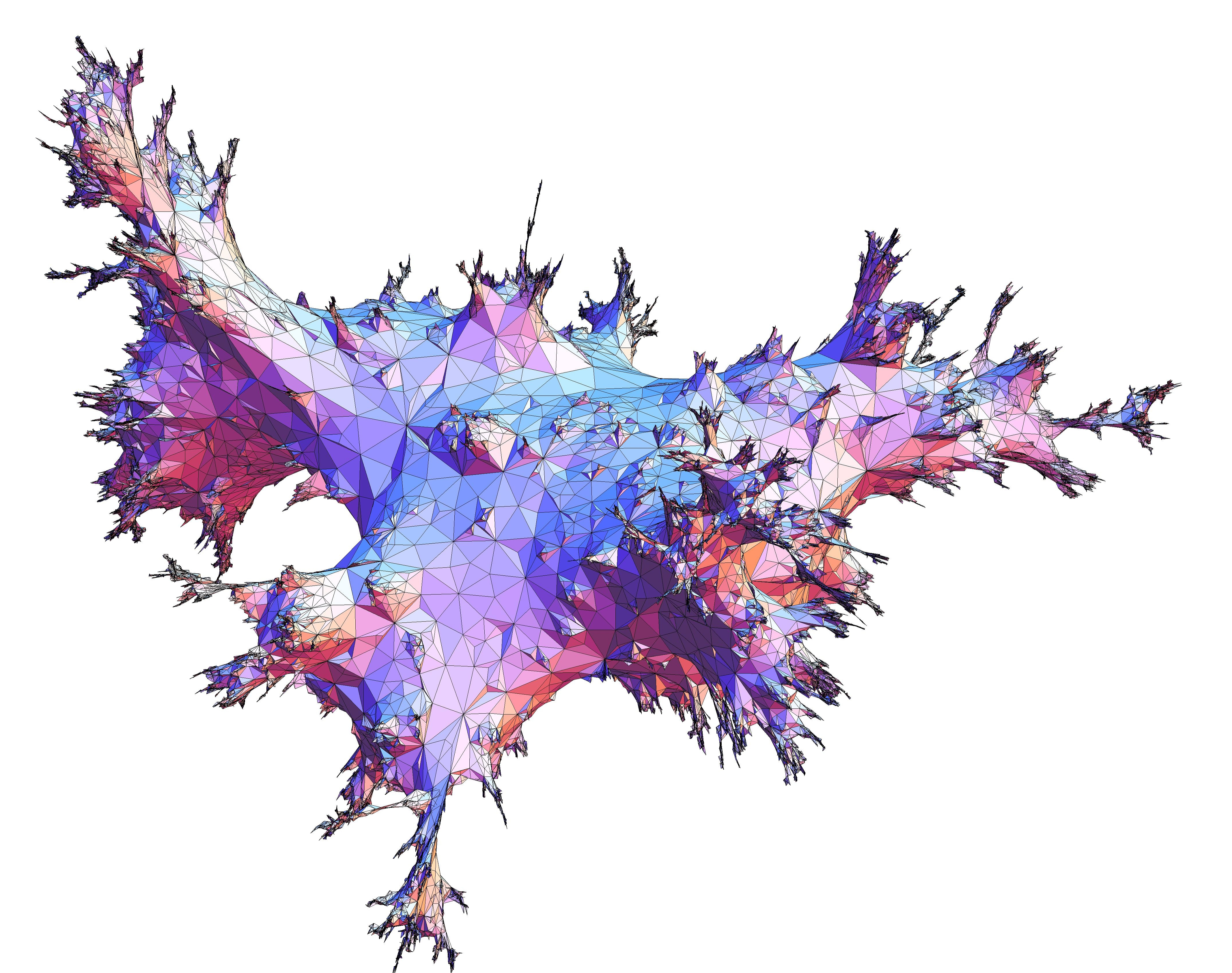}
\end{center}

\noindent {\textbf STUDENT:} You just listed four things! Which is the $n \to \infty$ limit of Budzinski's picture?
\vspace{.03in}

\noindent {\textbf INSTRUCTOR:} They all are! The difference comes down to the topology of convergence and the features of the limit presumed to be measurable. Think of them as different aspects of the same universal object. Four blind mathematicians feel the surface of an elephant and describe four different things:
\vspace{-.1in}
\begin{center}
\includegraphics[scale=0.6]{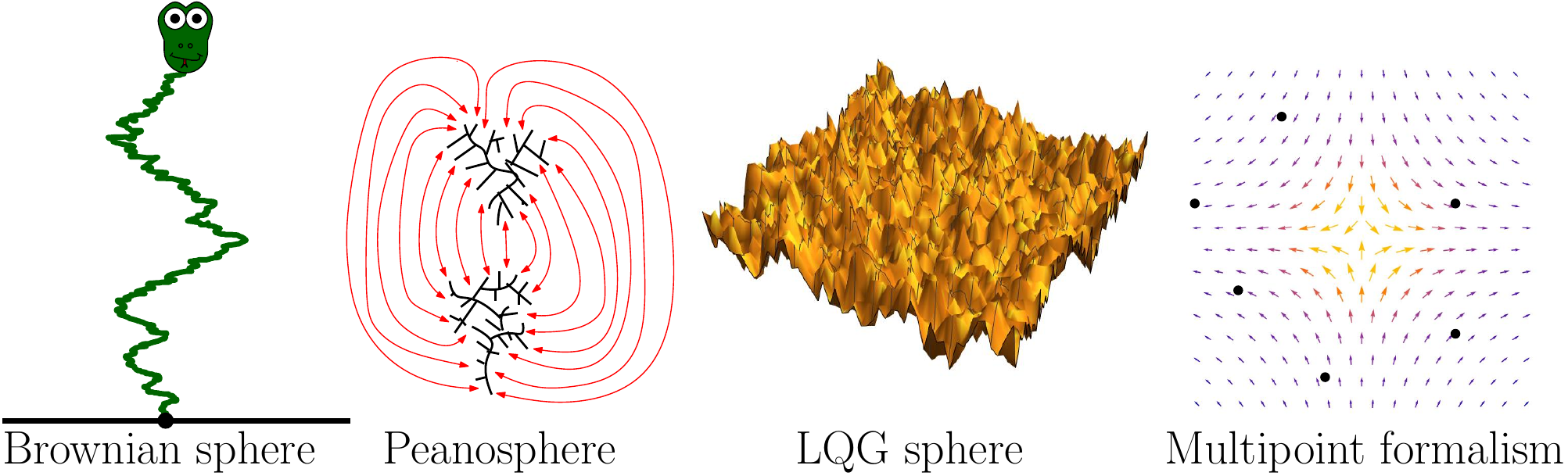}
\end{center}

\vspace{-.4in}

\begin{itemize}[leftmargin=.08in]
\item \textbf{Brownian sphere:} a {\em random metric measure space} constructed from the so-called {\em Brownian snake} as explained in surveys by Le Gall, Miermont and Baez \cite{gall2014random, le2019browniangeometry, miermont2009random, miermont2014aspects,baez2021brownian}.

\item \textbf{Peanosphere:} a {\em mating of continuum random trees} that encodes both a surface and an {\em extra tree and/or collection of loops} drawn on top of it, as explained in surveys by Gwynne, Holden and Sun and by Biane \cite{gwynne2019mating, biane2021mating}.

\item \textbf{Liouville quantum gravity sphere:} a {\em random fractal Riemannian surface} in which areas, lengths and other measures are given by exponentials of a {\em Gaussian free field} $\phi$, as explained in surveys by Berestycki, Ding, Dub\'edat, Duplantier, Garban, Gwynne, Miller, Powell and Werner \cite{berestycki2015introduction, garban2013quantum, gwynne2020notices, miller2018liouville, berestycki2021gaussian, ding2021introduction, powell2022lecture, duplantier2014ICM}.

\item \textbf{Conformal field theory:} a {\em collection of multipoint functions} representing  (regularized) integrals of products of the form $\mathbb \prod e^{\alpha_i \phi(x_i)}$ w.r.t.\ a certain infinite measure, as explained mathematically by David, Guillarmou, Kupiainen, Rhodes and Vargas \cite{david2016liouville, kupiainen2020integrability, guillarmou2021segals, guillarmou2020conformal, vargas2017lecture} who also survey the physics literature. The infinite measure is the {\em Polyakov measure} which is the product of an unrestricted-area measure on LQG spheres (with defining field $\phi$) and Haar measure on the M\"{o}bius group PSL$(2, \mathbb C)$ (to select an embedding in $\mathbb C$).  
\item[] \textbf{\em Technical point:} The {\em unit area} Brownian/Peano/LQG-sphere is a sample from a probability measure $dS$ on a space of unit area surfaces. There is a natural {\em infinite} measure on {\em unrestricted-area} surfaces (with $k \geq 0$ {\em marked points}) given by $A^{-7/2 + k}dAdS$ where $(A,S)$ is the rescaling of $S$ with area $A$. This is natural because the number of triangulations (quadrangulations, etc.) with $n$ faces and $k$ marked points scales like $C \beta^n n^{-7/2+k}$ for model-dependent constants $C$ and $\beta$. Weighting the counting measure by $\beta^{-n}$ we obtain a discrete measure that (appropriately rescaled) converges to the measure above as the area-per-triangle $\epsilon$ goes to zero. If we replace $\beta$ by the ``off-critical'' $\beta(1+ \epsilon \mu)$ then the limit is $A^{-7/2 + k}e^{-\mu A} dAdS$, which is finite if $\mu >0$ and $k \geq 3$. Polyakov included the $e^{-\mu A}$ factor and motivated it in a different way, using the {\em Liouville equation} as we explain later.

\end{itemize}

\vspace{.03in}
\noindent {\textbf STUDENT:} Do I really have to learn all four viewpoints?
\vspace{.03in}

\noindent {\textbf INSTRUCTOR:} A lot of good work has been done by people fluent in only one of the four. But all four have important applications. For example, in the Brownian sphere construction {\em distances} are easy to define, and you can show that the surface has fractal dimension $4$. So to cover the surface with metric balls of radius $1/n$, you need $\approx n^4$ balls.
In the peanosphere construction, one sees the random trees and loops that naturally live on top of a random surface, and most easily makes contact with loop-decorated discrete models. The LQG approach allows one to define {\em Brownian motion on a random surface} as well as other conformally defined objects (like SLE curves). Conformal field theory is a staple of physics. It ties into quantum field theory (if one analytically continues from Euclidean to Minkowski space and replaces probability measures with quantum wave functions).
\vspace{.03in}

\noindent {\textbf STUDENT:} Are these surfaces characterized by simple axioms like Brownian motion is?
\vspace{.03in}

\noindent {\textbf INSTRUCTOR:} Yes. But each viewpoint comes with its own axioms, its own history, its own motivation, its own surveys. The proofs that they agree are long and involved.
\vspace{.03in}

\noindent {\textbf STUDENT:} Can any of the viewpoints describe a random surface embedded in $\mathbb R^d$?
\vspace{.03in}

\noindent {\textbf INSTRUCTOR:} Sure. You basically weight the law of the surface by the ``number of ways'' to embed it in $\mathbb R^d$. This amounts to weighting by the $d$th power of a certain ``partition function'' and this changes the law of the surface itself (even ignoring the embedding). You can apply such a weighting even when $d$ is not a positive integer, so there is actually a one parameter family of random surfaces parameterized by $d$, with the scaling limit of Budzinski's picture corresponding to $d=0$. These surfaces are ``rougher'' when $d$ is large and ``smoother'' when $d$ is small, converging to the Euclidean sphere as $d \to -\infty$. They are defined as random metric spaces for any $d < 25$, but are only finite-diameter and finite-volume if $d \leq 1$. The family can equivalently be parameterized by related quantities that come up in LQG theory like $Q>0$ (where $d = 25-6Q^2$) or $\gamma$ (where $Q=2/\gamma + \gamma/2$ --- note that $\gamma$ is only real if $d \leq 1$ so that $Q \geq 2$) or by the peanosphere correlation coefficient.
\vspace{.03in}

\noindent {\textbf STUDENT:} I'm getting a bit lost. Can you give me the four definitions you promised?
\vspace{.03in}

The remainder of this paper will aim to do just that. We present a narrative-style introduction to each of the four viewpoints above, emphasizing the distinctive intellectual heritage behind each approach, as well as the relationships between them and the relevant recent developments.

Along the way, we will introduce several other natural objects: random trees (like the continuum random tree), random distributions (like the Gaussian free field) and random non-self-crossing curves (like the Schramm-Loewner evolution). We aim to make the exposition accessible to outsiders and newcomers, as well as to researchers who have expertise in one or more of the viewpoints but who might appreciate a high-level overview of the others. Readers interested in a still less technical account might take a look at the recent Quanta Magazine articles on this subject \cite{quantamagazine, quantamagazinestable, quantamagazineCFT}.

It is ambitious to try to tell four stories in one article, but we'll do our best to convey at least the main ideas. The style is informal---no proofs---but we provide accurate definitions and the surveys mentioned above contain more detail. This paper is targeted primarily at a mathematical audience and is meant to be broadly accessible. This subject has deep roots in physics (especially string theory, quantum field theory, and planar statistical physics) but we will do our best to avoid any terminology that would be difficult for non-physicists to understand. Let us also stress that our reference list is long but very far from complete, biased by both the limits of the author's knowledge and the narrative focus of the paper. As such, it is more of a sampling than an exhaustive survey. Many highlights of the subject are not mentioned here.

For concreteness and simplicity we will focus mostly on {\em sphere-homeomorphic} random surfaces. But all of the constructions in this paper have analogs that look the same locally but have different global topology. For example, in addition to the Brownian sphere one has a Brownian disk, a Brownian plane, a Brownian torus, Brownian surfaces of arbitrary genus and/or arbitrarily many boundary components, etc. \cite{curien2014brownian, bettinelli2017compact, budzinski2021local, le2019brownian, riera2021brownian, miermont2021compact}. Variants like these can be defined within all four viewpoints. One would expect that these definitions are consistent from one viewpoint to another (so that e.g.\ a Brownian three-holed torus is somehow equivalent to a corresponding LQG three-holed torus or Peano-three-holed-torus). Significant work along these lines has been done, but the set of questions one can ask is very large and the program has not been completed yet.

\section{Brownian sphere: a random metric measure space} \label{sec::browniansphere}
Our first construction is a random surface called the {\em Brownian sphere} (a.k.a.\ {\em Brownian map}). The Brownian sphere is a random sphere-homeomorphic metric measure space. It has been described in longer survey articles by Le Gall and Miermont \cite{gall2014random, le2019browniangeometry, miermont2009random, miermont2014aspects, miermont2021compact} and in a shorter overview by Baez \cite{baez2021brownian}.

Both the peanophere and the Brownian sphere can be constructed by ``gluing together'' a pair of ``continuum trees'' along their outer boundaries, which produces a sphere decorated by a space-filling curve (that somehow snakes in between the two trees)---but the law of the pair of trees is different in the two settings. The idea of gluing together trees to obtain a sphere may seem counterintuitive, but we will see that it is well motivated by the discrete models.

Readers familiar with complex dynamics may also recall that the ``mating'' of two dendritic Julia sets (these are tree-like sets with empty interior) is a topological sphere, see e.g\ \cite{milnor2004pasting, yampolsky2001mating, timorin2010topological, aspenberg2009mating} or look up the online video animations by Arnaud Ch\'eritat. The Brownian sphere and peanosphere constructions are random versions of this phenemonon. Let us begin by defining {\em continuum random trees} which are sometimes also called {\em Brownian trees}.

\subsection{Continuum random trees (a.k.a.\ Brownian trees)}

Brownian motion is a random function $B(t)$ defined for $t \in [0,\infty)$. But Brownian motion (like the Brownian sphere) has many {\em variants} with similar local behavior. For example, the {\em Brownian bridge} is a Brownian motion on $t \in [0,1]$ somehow {\em conditioned} on the (zero probability) event that $B(1)=0$. The {\em Brownian excursion} is a Brownian bridge with $B(0)=B(1)=0$ that is further conditioned to stay positive on $[0,1]$.

If one starts with the graph of a Brownian excursion $B$, and identifies two points whenever they are connected by a horizontal chord under the graph (see below), one obtains a random metric space $\mathcal T$ called the {\em continuum random tree}, which was first constructed by Aldous in 1991 \cite{aldouscontinuum1, aldouscontinuum2, aldouscontinuum3}, and which plays a role in the construction of both the Brownian sphere and the peanosphere.
\begin{center}
\includegraphics[scale=.9]{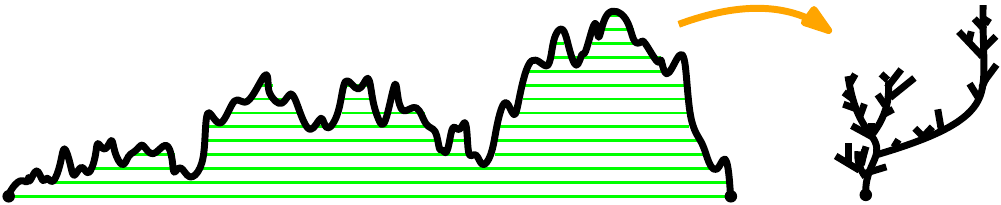}
\end{center}
\vspace{-.1in}
Note that the function taking $t$ to the graph point $\bigl(t,B(t)\bigr)$ induces a function $f$ from $[0,1]$ to $\mathcal T$ that ``traces the boundary of the tree clockwise.'' The pushforward of Lebesgue measure on $[0,1]$ endows $\mathcal T$ with a natural measure.  Moreover, one can define the {\em distance} between $f(s)$ and $f(t)$ to be $$\Bigl( B(s) - \inf_{r \in [s,t]} B(r) \Bigr) + \Bigl( B(t) - \inf_{r \in [s,t]} B(r)\Bigr).$$
This measures how far ``down and up'' one has to travel within $\mathcal T$ to get from $f(s)$ and $f(t)$. Since $\mathcal T$ comes with a natural measure and a natural metric, we call it a  {\em random metric measure space}.

\subsection{Planar map bijections that motivate Brownian sphere and peanosphere}
The Brownian sphere has its historical roots in the study of planar maps.  A {\em planar map} is a finite graph together with an embedding in the complex sphere $\mathbb C \cup \{\infty\}$, where two embeddings are considered equivalent if there is an orientation-preserving homeomorphism of $\mathbb C \cup \{\infty\}$ taking one to the other. For example, the two figures below are isomorphic as graphs but represent different planar maps, as there is no orientation-preserving homeomorphism of $\mathbb C \cup \{\infty\}$ mapping one to the other.

\begin{center}
\includegraphics[scale=.9]{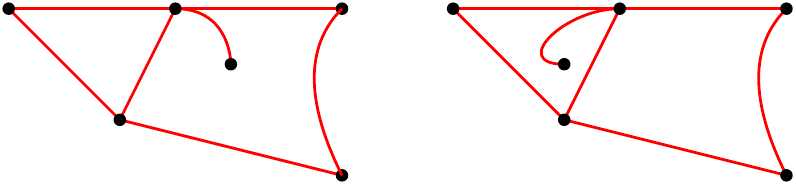}
\end{center}
\vspace{-.1in}
Combinatorially, a planar map is fully determined by the graph together with the ``clockwise cyclic ordering'' of the edges surrounding each vertex. In particular, the number of distinct planar maps with $n$ edges is finite. To eliminate the ambiguity that comes from having non-trivial automorphisms, it is sometimes helpful to specify a ``root'' by fixing an oriented edge. If an orientation-preserving homeomorphism of $\mathbb C \cup \{\infty\}$ takes a planar map to itself (mapping edges to edges and vertices to vertices) and fixes the oriented edge, one can show that it must induce the identity map on the whole set of vertices and edges. To see why, imagine that your car starts driving along the oriented edge. From there you can describe how to get to any other vertex or edge with a set of directions like ``Drive until you reach a vertex, then take the furthest road to your left, then drive until you reach another vertex, then take the third road from the left,'' and so forth. This method of specifying vertices and edges is preserved by homeomorphisms of $\mathbb C \cup \{\infty\}$. In other words, fixing a root gives us a way to uniquely specify all other vertices (just as a Cartesian coordinate system gives us a way to uniquely specify points on a planar lattice). This implies that fixing the root eliminates all non-trivial planar map automorphisms: for example, a triangle has three non-trivial orientation-preserving automorphisms---corresponding to the three rotations---but if the triangle is assigned a root, then only the identity automorphism would fix the root.

Perhaps the most famous planar map problem of all is the ``4 color conjecture'' (a.k.a.\ Guthrie's problem) which captivated mathematicians worldwide from its formulation in 1852 to its computer-assisted proof by Appel and Haken in 1976. See e.g.\ \cite{heawood1890map, coxeter1959four, fritsch1998four,  rogers2008four} for historical accounts. The problem is to show, given any finite planar map, that it is possible to color each vertex one of four colors in such a way that no two neighbors have the same color. Note that if one has a map of South America, say, one can create a {\em planar map} by putting a vertex in the center of each country and drawing an edge connecting two countries wherever they share a non-zero-length border. The four color theorem allows one to color the countries in such a way that no two countries sharing a non-zero-length border have the same color. One of the key contributors to research on this problem, and to graph theory in general, was William Tutte.

By Tutte's own account \cite{tutte1998graph} (see also \cite{noy2021expected}) he was motivated by the four color theorem when he began thinking about enumerating {\em all} planar maps. He wrote his famous ``census of'' series in 1962 and 1963 \cite{tutte1962slice, tutte1962census, tutte1962hamiltonian, tutte1963census} which among other things included a remarkable formula for the number of rooted planar maps with $n$ edges:  namely $\frac{2}{n+2} \cdot \frac{3^n }{n+1} {2 n \choose n}$.  Using e.g.\ Stirling's formula (or the local central limit theorem) one can show that
 ${2n \choose n} = \frac{4^n}{\sqrt{\pi} n^{1/2}} \bigl(1 + O(\frac{1}{n}) \bigr)$
 which implies that Tutte's formula grows asymptotically like $C \beta^{n} n^{-7/2+k}$ (which is the formula from the {\em technical point} in our introductory dialog) with $C = \frac{2}{\sqrt{\pi}}$ and $\beta = 12$ and $k=1$. This is consistent with the formula from the introduction because the root effectively plays the role of a single marked point and we recall that $k$ represents the number of marked points.

These papers were followed by another series by Mullin \cite{mullin1964enumeration, mullin1965counting, mullin1966enumeration, mullin1967enumeration} which included a formula for the number of rooted planar maps with $n$ edges decorated by a distinguished spanning tree: namely, $\frac{(2n)!(2n+2)!}{n![(n+1)!]^2(n+2)!}$. This grows asymptotically like $C \beta^{n} n^{-4+k}$ with $k=1$ and $\beta = 16$. The addition of the ``spanning-tree decoration'' somehow changes the growth rate in a fundamental way, changing not only $\beta$ and $C$ but also the power of $n$ (effectively replacing the $7/2$ with $4$).

Both of these formulas have interesting bijective proofs: the former, in its first form due to Cori and Vauquelin in 1981 \cite{cori1981planar}, was later advanced and popularized by Schaeffer \cite{schaeffer1998conjugaison}. The latter was essentially due to Mullin in 1967 \cite{mullin1967enumeration}, see also the explanation by Bernardi \cite{bernardi2006bijective}, and it reduces the problem of counting spanning-tree-decorated rooted planar maps to a problem about walks in $\mathbb Z_+^2$ (which can be separately addressed e.g.\ with reflection arguments). These bijections motivate the definition of the Brownian map and the peanosphere, respectively. For now let us describe the bijections side by side, beginning with the Mullin bijection.

Suppose that $(X_n, Y_n)$ is a simple walk in $\mathbb Z^2_+$ starting and ending at the origin. Then if we fix a large enough value of $C$, the graphs of $X_n$ and $C-Y_n$ (linearly interpolated to $\mathbb R$) will not intersect; see the figure below. Draw a vertical red line at each time increment. Then declare two points in the graph of $X_n$ (resp.\ $C -Y_n$) to be equivalent if one can draw a horizontal chord connecting them that does not go above the $X_n$ graph (resp.\ below the $C-Y_n$ graph).  In other words we ``identify'' each pair of points connected by a horizontal (blue, black or green) line segment. We also glue together the leftmost and rightmost red lines. After this gluing is done, the region below $X_n$ collapses to become a tree (shown with black edges and green vertices) as does the region above $C-Y_n$ (shown with black edges and blue vertices) and we are left with a planar map with black/red edges and blue/green vertices.  (The left and right red curves are glued together, which makes the topological disk they enclose into a sphere.)  This planar map is a {\em triangulation} in which each triangle contains two red edges and one black edge.

\begin{center}
\includegraphics[scale=.9]{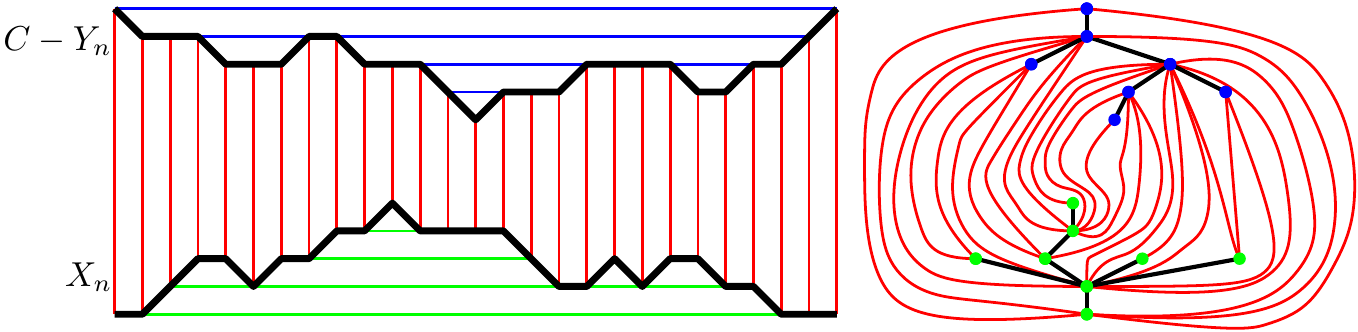}
\end{center}
\vspace{-.1in}
If the black edges are erased, one is left with a red {\em quadrangulation} $\mathcal Q$ which is bipartite: the blue and green vertices are the two partite classes.  Each quadrilateral in $\mathcal Q$ has one blue-to-blue diagonal.  These diagonals together form a planar map $M$. Conversely, given the planar map $M$ it is not hard to show that one can recover $\mathcal Q$ (by adding a vertex in the center of each face of $M$ and connecting it to all of the boundary vertices of that face). The blue tree $T$ is a spanning tree of $M$.  The green-to-green diagonals of $\mathcal Q$ form the dual graph $M^*$, and the green tree $T^*$ is the dual spanning tree. The story above gives a bijection between:
\begin{itemize}
\item Simple walks $(X_n, Y_n)$ in $\mathbb Z_+^2$ of length $2N$ that start and end at the origin.
\item Pairs $(M,T)$ where $M$ is a rooted planar map and $T$ is a spanning tree of $M$. Here {\em rooted} means that a distinguished vertex of $M$ and a distinguished incident vertex of $M^*$ are fixed (to be the roots of $T$ and $T^*$ --- i.e.\ the two vertices attached to the leftmost red line, or the equivalent rightmost red line). 
\end{itemize}
Choosing a uniformly {\em random} walk $(X_n, Y_n)$ of the type described above produces a uniformly {\em random} $(M,T)$ pair. The probability of a given map $M$ is proportional to the number of spanning trees $M$ has.

But what if we want to simply count the total number of rooted maps $M$ (or the total number of quadrangulations $\mathcal Q$) instead of the number of $(M,T)$ pairs? This is what the Cori-Vauquelin-Schaeffer bijection does. It can be seen as similar to the Mullin bijection but with a few key differences. First, instead of requiring $(X_n,Y_n)$ to traverse lattice edges we at each step allow $Y_n$ to change by $\pm 1$ and $X_n$ by either $0$ or $\pm 1$. Second, instead of perfectly horizontal green chords, we draw chords that are one unit higher on the right than on the left. We draw one such chord leftward starting at each vertex on the graph of $X_n$, which means that we have to add an extra vertex of minimal height as shown.

\begin{center}
\includegraphics[scale=1]{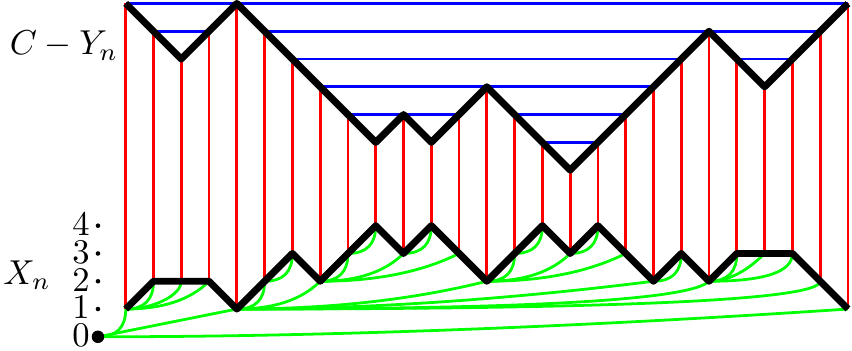}
\end{center}
\vspace{-.1in}
Third, we consider only $(X_n,Y_n)$ pairs for which the above picture has a special property: namely, whenever two red vertical lines are incident to the same blue chord, their lower endpoints have the {\em same height}.
\begin{center}
\includegraphics[scale=.8]{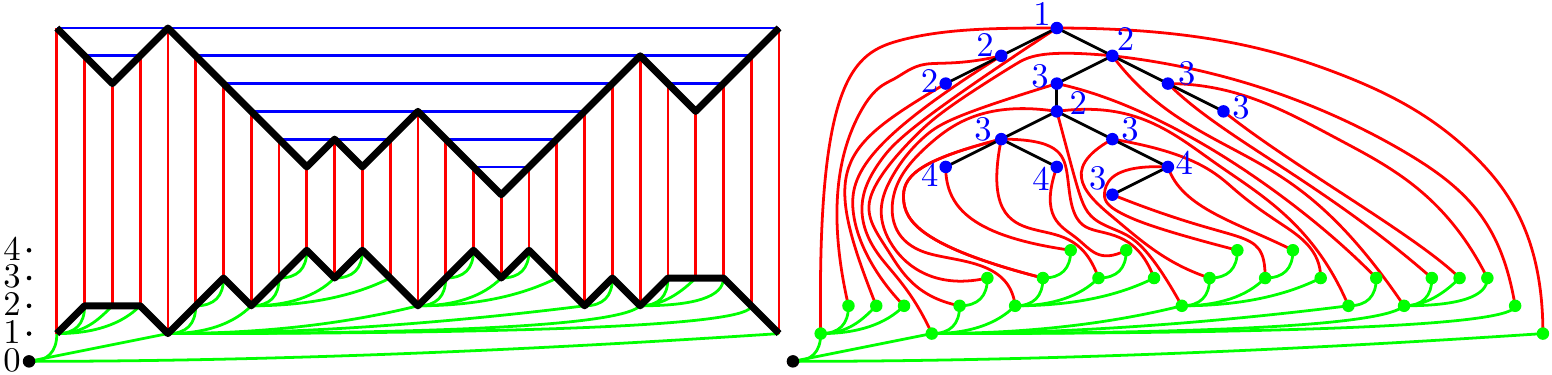}
\end{center}
\vspace{-.1in}
If we glue together the edges on the $C-Y_n$ graph connected by blue chords (the same way as in the Mullin bijection) we obtain a tree, and the third condition above is equivalent to the condition that if two red lines start at the same blue vertex on that tree then they terminate at green vertices of the same height. Thus we can label each vertex in the tree by the height of the green vertex (or vertices) it connects to. We then obtain a planar map by {\em gluing} the black and green trees on the right to one another: precisely, we glue each green vertex in the tree on the right to the blue vertex it is connected to by a red edge, effectively shrinking each red edge to a point. (Also, the outermost two green edges are understood to be glued/identified with each other, so that the topological disk they surround becomes a sphere.)
\begin{center}
\includegraphics[scale=.7]{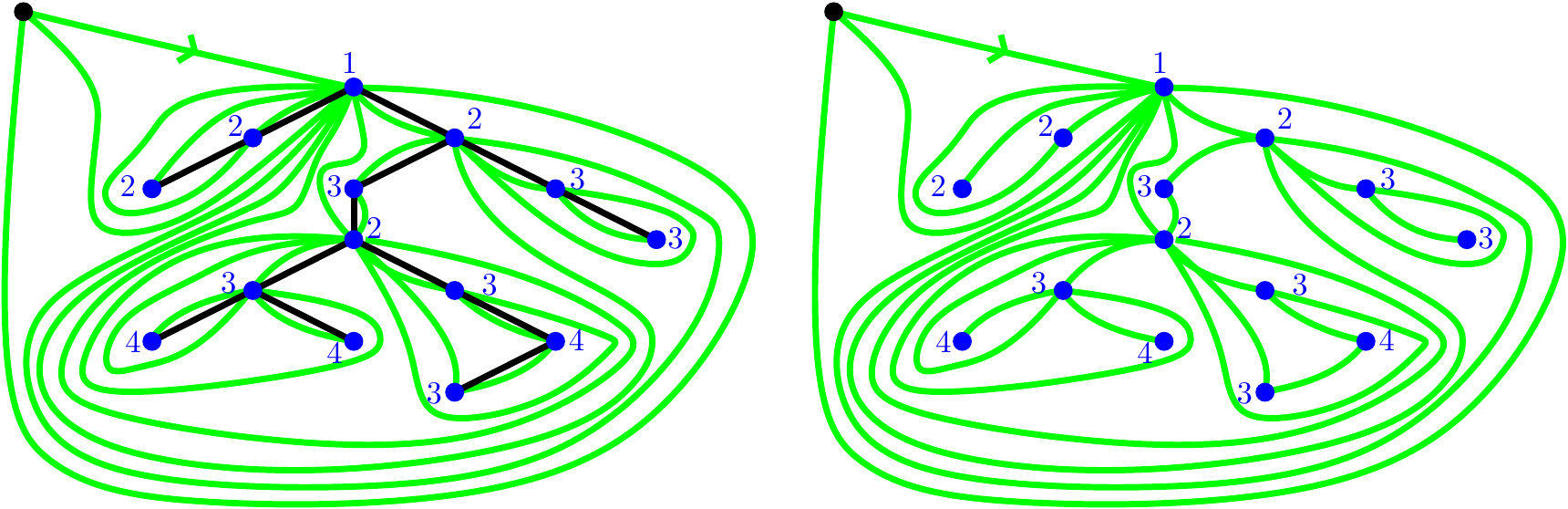}
\end{center}
\vspace{-.1in}
We can think of the graph on the left above as obtained by starting with the region below the $X_n$ graph, then gluing together two black edges of the $X_n$ graph whenever the edges of the $C-Y_n$ graph just above them correspond to the same edge of the upper tree. Whenever two of these black edges are glued together, they are either both horizontal (each part of a triangle with two green edges) or one increasing and one decreasing (respectively part of a two-gon with one green edge and a quadrilateral with three green edges). In each case, once the two black edges are glued together (and erased) we are left with a green quadrilateral.

The construction above yields a bijection between
\begin{itemize}
\item Well-labeled rooted planar trees $(T,\ell)$ (here $\ell$ maps vertices of $T$ to positive integers, where the root has label one and adjacent vertices differ by $0$ or $\pm 1$).
\item Rooted quadrangulations $\mathcal Q$.
\end{itemize}

To prove that the above is a bijection, one has to show that one can reconstruct the green and black trees given just the green-edge quadrangulation $Q$. Roughly speaking, one reconstructs the original green tree as follows. First we shorten the green edges slightly so that they are no longer connected to each other at the vertices, but then we reconnect some of them: precisely, we connect each edge (directed ``downward'' towards its lower-label endpoint) to the successor downward-pointing edge that involves ``turning maximally to the left.'' Once this is done, the green tree is a tree of ``leftmost geodesics'' built from the directed green edges (each directed in the direction of decreasing distance to the root). Given the green tree, one can construct the $X_n$ in the figure above and it is not too hard to check that one can reconstruct the upper tree as well (since one knows which green vertices have to be equivalent---and hence connected by red-blue-red paths).

We remark that once the black and green trees are glued together, they have vertices (but no edges) in common (this is different from the Mullin bijection) but in a sense the two trees still do not ``cross'' each other.  To try to see why, it is interesting to imagine a path that crosses through the red edges in order from left to right (but does not intersect the green or black trees). Once we take the quotient (w.r.t.\ to the equivalence class where all points on the same red edge are considered equivalent) this becomes a path that winds between the green and black trees, intersecting the trees at vertices but never crossing a branch of either tree.

This bijection (like the Mullin bijection) constructs a quadrangulation $\mathcal Q$ by gluing together a pair of trees. The microscopic details of the gluing are a bit different: in the Mullin bijection the red edges form $\mathcal Q$ and in the Cori-Vauquelin-Schaeffer bijection the green edges form $\mathcal Q$. There are many variants and analogs of the CVS bijection, see for instance \cite{BDGmobiles, bouttier2004planar, bernardi2012unified, albenque2015generic, ambjornbudd, schaeffer1998conjugaison, bettinelli2014scaling} and the references therein. At the global level, the main difference between the CVS bijections (and its variants) and the Mullin bijection (and its variants) is this: in the CVS bijection, one of the trees (the green one) is the so-called {\em leftmost geodesic tree of edges} and is determined by $\mathcal Q$ itself. The second tree (the black one) is in some sense ``dual'' to the geodesic tree (and is also determined by $\mathcal Q$ itself). By contrast, in the Mullin bijection the spanning tree is not determined by $\mathcal Q$ and $\mathcal Q$ is not uniform among all quadrangulations.

\subsection{Unconstrained variant: when $\mathcal Q$ is both pointed and rooted}
There is a variant of the Mullin bijection in which we relax the restriction that $X_n$ is non-negative, see below.

\begin{center}
\includegraphics[scale=.9]{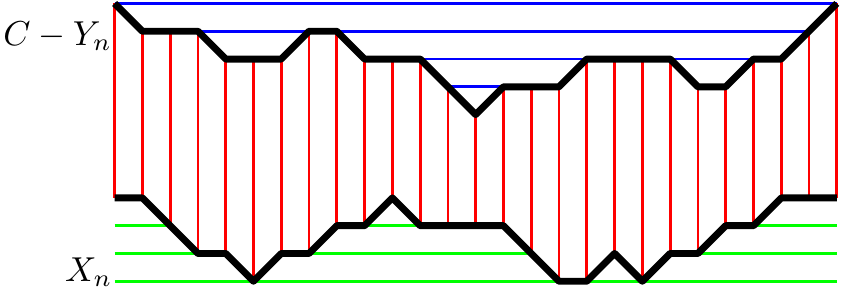}
\end{center}
\vspace{-.1in}
Here we imagine that the left and right sides of the above rectangle are glued to one another (so that both $X_n$ and $Y_n$ then become indexed by a circle).  We can then identify points on the same blue chord or the same green (possibly wrapping around) chord, just as before, and we obtain a pair of trees with red edges between them.  A key difference from the original construction is that in this construction the root of the lower tree (corresponding to the minimum of $X_n$) and the root of the upper tree (corresponding to the minimum of $Y_n$) are no longer required to be adjacent.  The minima of $X_n$ and $Y_n$ may occur at different places, and we can think of the vertex corresponding to the $X_n$ minimum as an extra marked point. (In principle, one could allow both $X_n$ and $Y_n$ to take negative values, in which case the two tree roots---and the point described by the leftmost/rightmost red line---would effectively describe three marked points.)

It turns out that the CVS construction also has a similar variant, which corresponds to relaxing the constraint that all of the blue labels are positive. If we relax this constraint, then the corresponding $X_n$ process can be positive and negative, and in particular no longer has to have a minimum at the same place that $Y_n$ does.  This means that the root of the upper tree and the root of the geodesic tree are no longer required to be adjacent. With this constraint relaxed, the set of possible labeled upper trees is easier to count: one has a Catalan number $\frac{1}{n+1} {2n \choose n}$ of rooted planar trees, and given the tree, there are exactly $3^n$ ways to choose the labels if we fix the root label to be zero. (We are free to do this since adding a constant to the labels does not affect the map construction, and hence only the labels modulo additive constant are relevant.) The formula $\frac{3^n}{n+1} {2n \choose n}$ can also be obtained by starting with Tutte's formula and multiplying by $n+2$ (there are $n+2$ possible to choose an extra marked vertex) and dividing by a factor of $2$ (which is related to the choice of root orientation---this is explained in many places, see e.g.\ \cite{schaeffer1998conjugaison} or the short explanation in Section 3 of \cite{caraceni2020polynomial}).

There is another way to think about the trajectory $(X_n, Y_n)$ obtained from a labeled tree.  For every point on a labeled tree, there is a labeled path from that point back to the root.  As one traces the tree clockwise, one obtains a corresponding sequence of labeled paths; the labeled path seen at time $n$ can be drawn as a ``snake'' (a vertical-to-horizontal function defined on $[0,Y_n]$) with the horizontal coordinate indicating the label.

The figures below (read from left to right) are steps in a Markov process on a space of ``snakes.''  Each figure is a sequence of edges, each of which goes up one unit and $-1$, $0$ or $1$ units to the right.  To move from one figure to the next, one first tosses a coin to decide whether to delete the top edge (with probability $1/2$) or to add an up-left edge (probability $1/6$), up-right edge (probability $1/6$) or up edge (probability $1/6$).

\begin{center}
\includegraphics[scale=.7]{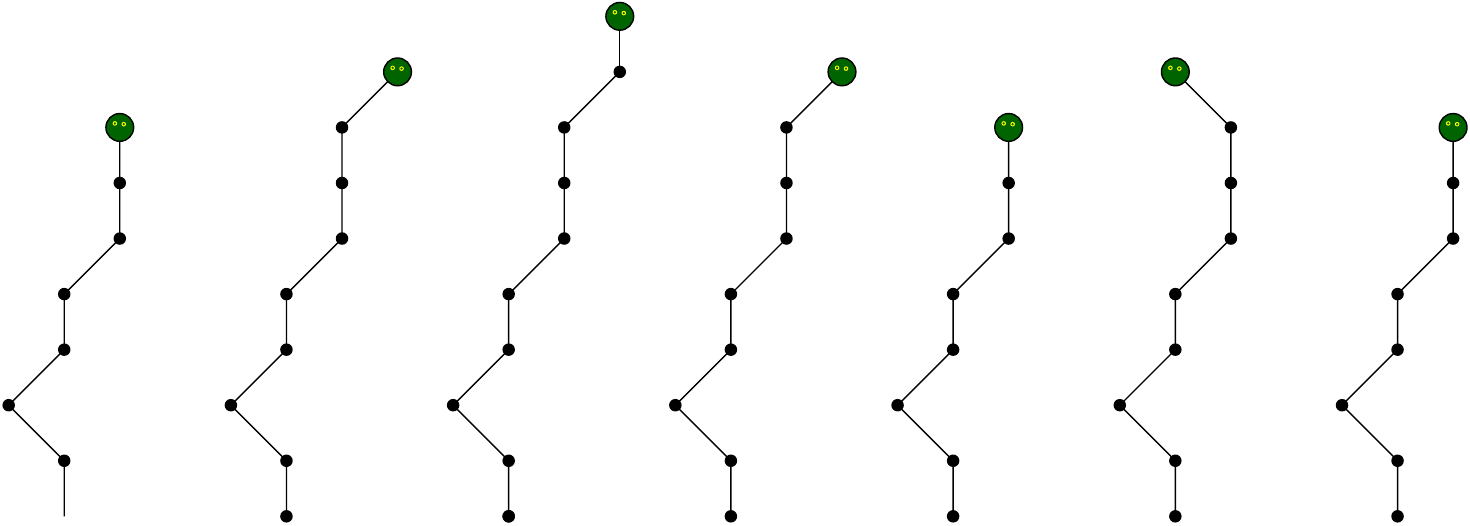}
\end{center}
\vspace{-.1in}
If the head of the snake starts at height zero---and we condition on the head height staying non-negative and returning to zero at time $2n$---then we can let the head height represent the process $Y_n$ and the head horizontal location represent the process $X_n$, and it is not hard to see that this corresponds precisely to the $(X_n, Y_n)$ construction from the labeled tree.   If one rescales this process vertically by a factor of $C$ and horizontally by $\sqrt{C}$ one obtains a continuum {\em Brownian snake} model in the limit.  The Brownian snake was introduced by Le Gall in the 1990s (see e.g.\ \cite{le1995brownian}) not long after Aldous introduced the continuum random tree (a.k.a.\ Brownian tree).  The terms {\em Brownian snake} and {\em Markov snake} were coined by Dynkin and Kuznetsov (who credit Le Gall for the construction) in 1995 \cite{dynkin1995markov}.
\begin{center}
\includegraphics[scale=.7]{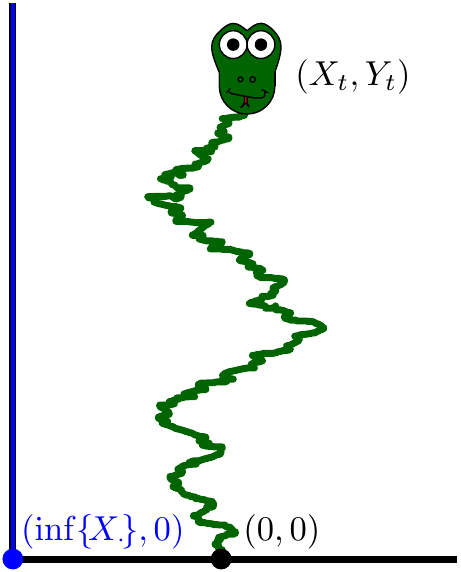}
\end{center}
\vspace{-.1in}
Now the law of the continuum pair $(X_t, Y_t)$ can be described as follows: one first samples $Y_t$ as a Brownian excursion on $[0,1]$ and then takes a quotient of the graph of $Y_t$ to produce a Brownian tree $\mathcal T$, together with the natural map $f$ taking $t$ to $\mathcal T$ (which traces around the boundary of $\mathcal T$ clockwise as $t$ varies from $0$ to $1$).  One then generates a Brownian motion $\zeta$ indexed by $\mathcal T$ (taken to be zero at the root of $\mathcal T$).  This in turn determines a process $X_t$ by $X_t = \zeta(f(t))$. The minimum of this process corresponds to the location of the root of the geodesic tree.  The picture of the whole snake at time $t$ (viewed as a vertical-to-horizontal function as above) is described by the restriction of $\zeta$ to the path in $\mathcal T$ from $f(0)$ and $f(t)$.

\subsection{Passing to the continuum}

One would expect both the Mullin and the Cori-Vauquelin-Schaeffer bijections to have limits given by the gluing of two fractal trees defined by continuum processes $X_t$ and $Y_t$, something like the figure below:

\begin{center}
\includegraphics[scale=.9]{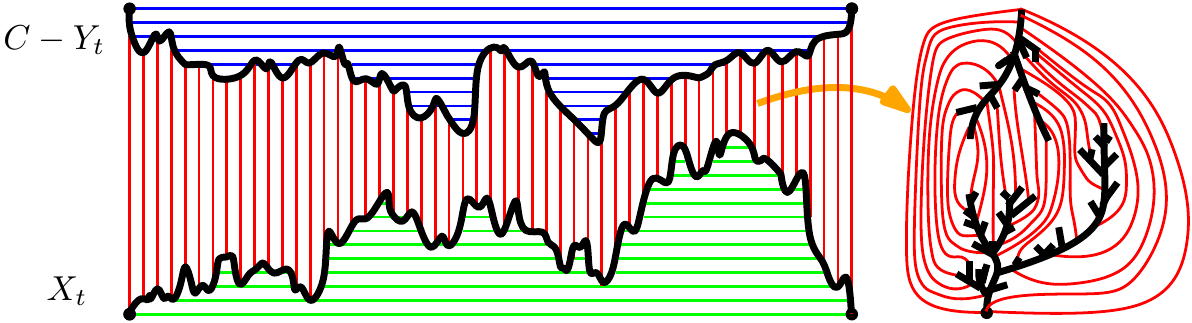}
\end{center}
\vspace{-.1in}
 However, we would expect the law of $(X_t,Y_t)$ to be rather different in the two scenarios. In the Mullin limit, $X_t$ and $Y_t$ are independent Brownian excursions, and the picture is a topological quotient of two independent continuum random trees (w.r.t.\ the equivalence that identifies two points if they are connected by a red curve). The topological sphere constructed this way is a special case of the {\em peanosphere} that we will discuss in the next section.  In the CVS limit, $X_t$ and $Y_t$ correspond to the horizontal and vertical components of the head location in the Brownian snake process.  The {\em Brownian sphere} (a.k.a.\ {\em Brownian map}) is defined to be the metric space obtained by starting with the tree generated by $X_t$ (viewed as a metric space) and then taking the metric quotient w.r.t.\ the equivalence that identifies points if they are connected (by the red curves) to the same point on the upper tree defined by $Y$.

The Brownian map was originally defined a bit differently. Marckert and Mokkadem \cite{marckert2006limit}, building on Schaeffer's work \cite{schaeffer1998conjugaison} (also citing e.g.\ the early work of Krikun and of Angel and Schramm \cite{angel2003uniform, krikun2005uniform}), argued that the limit of the discrete models could be defined in some sense (since the $(X_n,Y_n$) process had a limit) and coined the term {\em Brownian map} to describe the limiting object; however they did not show convergence in the space of random metric spaces. Le Gall showed that subsequential limits of random quadrangulations exist w.r.t.\ to the topology obtained from the {\em Gromov-Hausdorff metric} (a natural ``metric on the space of all compact metric spaces'') and are a.s.\ homeomorphic to a metric quotient of the Brownian tree \cite{le2007topological}. Le Gall and Paulin subsequently showed that all such subsequential limits are a.s.\ homeomorphic to the $2$-sphere \cite{le2008scaling}, see also the sphericity proof by Miermont \cite{miermont2008sphericity}.  The term {\em Brownian map} was sometimes used (see e.g.\ \cite{le2010scaling}) to describe any one of the subsequential scaling limits of the discrete models, and in this language, the question ``Do random quadrangulations have a (non-subsequential) Gromov-Hausdorff scaling limit?'' was formulated as the question: ``Is the Brownian map unique?'' This problem was solved independently by Le Gall and by Miermont \cite{le2013uniqueness, miermont2013brownian}, who also showed that this limit exists in law (without passing to a subsequence) and that the limiting random metric space is indeed the Brownian-snake-based metric quotient defined above.

Since this foundational work, Brownian surfaces have been the subject of a sizable literature. Many different kinds of random planar maps (quadrangulations, triangulations, various variants, etc.) have been shown to have the Brownian sphere as a scaling limit, strengthening the idea that the Brownian map (like Brownian motion) is something ``universal''  \cite{addario2017scaling, beltran2013quadrangulations, marzouk2018scaling, abraham2016rescaled, addario2021convergence, carrance2021convergence, bettinelli2014scaling, abraham2015random, albenque2020scaling, addario2019convergence}. There has also been significant work done on the continuum model. For example, the behavior of the geodesics in the Brownian sphere is quite interesting (geodesics ``merge into'' one another in ways that are very unlike what one sees for Euclidean spheres) and much has been done to understand their behavior \cite{gall2010geodesics, bettinelli2016geodesics, angel2017stability, miller2020geodesics, bouttier2009confluence, bouttier2008three, bouttier2003geodesic}.

\subsection{An axiomatic approach}
Like Brownian motion, the Brownian sphere is uniquely characterized by certain axioms, as shown in a long paper by the author and Miller \cite{axiomaticbrownianmap} which draws on related ideas by Bertoin, Budd, Curien, Kortchemski, Krikun, Le Gall, Miermont and others \cite{curien2016hull, cactusI, le2015brownian,krikun2005uniform,  angel2003uniform, bertoin2018martingales, bertoin2018random, le2018subordination, curien2019peeling}. We won't give a fully precise description here, but let us summarize the rough idea. For Brownian motion, there is a Markov property, which states that {\em given} $B(t)$ the conditional laws of $B$ restricted to $[0,t]$ and $B$ restricted to $[t, \infty)$ are independent.  Now consider an unrestricted-area Brownian sphere with two marked points $x$ and $y$. The natural Markov property in this setting states that {\em given} the boundary length of the filled metric ball of radius $r$ centered at $x$ (here the ``filled metric ball'' is obtained by starting with the ordinary metric ball and adding the components of its complement that don't contain $y$) the conditional laws of the ball and its complement are independent and depend on the boundary length in a scale invariant way. Furthermore, if we cut a filled metric ball into slices---using geodesics from evenly spaced points around the boundary to the center---then the slices are conditionally independent of each other given the length of their intersections with the metric ball boundary. The main result of \cite{axiomaticbrownianmap} (very roughly speaking) is that if a measure on sphere-homeomorphic metric measure spaces satisfies these properties, then it must be the Brownian sphere.

\section{Peanosphere: a random mating of trees}
\subsection{Basic definition}
In Section~\ref{sec::browniansphere}
 we saw that the Mullin bijection gives a bijection between lattice walks in $\mathbb Z_+^2$ (starting and ending at zero) and pairs $(M,T)$ where $M$ is a rooted planar map and $T$ is a spanning tree on $M$. The $X$ and $Y$ coordinates of the walk encode, respectively, a tree and a dual tree, which are somehow stitched together to create $M$. This is actually a special case of a much more general idea. The idea of considering planar maps {\em together with} an extra structure on the map (a spanning tree, a collection of loops, a distinguished edge subset, a bipolar orientation, etc.) has long been a staple of this subject, both on the mathematics side and the physics side. The extra structure is sometimes called a {\em decoration} or (on the physics side) a {\em statistical physics model} or a {\em matter field}.

It turns out that {\em many} kinds of planar map decorations can be used to produce a spanning-tree/dual-spanning-tree pair in some way. And the resulting trees are often encoded by {\em some kind of} lattice walk conditioned to stay in $\mathbb Z_+^2$, so that in the fine mesh scaling limit, one obtains {\em some kind of} Brownian motion conditioned to stay in $\mathbb R_+^2$ (starting and ending at the origin).  However, in general, the kind of Brownian motion involved (before imposing the quadrant constraint) may be one in which $X_t$ and $Y_t$ are {\em correlated}.  That is, the ``diffusion matrix'' may be such that $\mathrm{Var} (X_t) = \mathrm{Var}(Y_t) = t$ and $\mathrm{Cov}(X_t, Y_t) = \rho t$ for some possibly non-zero correlation coefficient $\rho$. (There is no loss in assuming $\mathrm{Var} (X_t) = \mathrm{Var}(Y_t) = t$, since multiplying $X$ or $Y$ by a constant does not change the topological surface construction.)

The limit of the Mullin construction corresponds to $\rho =0$, but in general one might consider any $\rho$ strictly between $-1$ (where the Brownian trees are perfectly negatively correlated) and $1$ (where the Brownian trees are perfectly correlated).  For each $\rho$ there is a natural way to make sense of this Brownian motion {\em conditioned} to stay in $\mathbb R_+^2$ (starting at the origin at time $0$ and returning there at time $1$).  One can then use the corresponding $(X_t,Y_t)$ pair to generate a pair of trees that can be ``mated together'' in the manner described in the previous section.

Formally, then, the peanosphere is just a random pair of measure-endowed metric planar trees that, when glued together along their boundaries, make a topological sphere, as described in the author's work with Duplantier and Miller \cite{matingoftrees, finitetree} and in a more recent survey by Gwynne, Holden and Sun  \cite{gwynne2019mating}, along with a survey by Biane about the discrete side \cite{biane2021mating}.  The sphere comes equipped with a non-self-crossing, space-filling path (a.k.a.\ {\em Peano curve}) which in some sense traces the interface between the two continuum trees. (This is the motivation behind the term {\em peanosphere} which was originally proposed by Richard Kenyon in private communication.) It is also equipped with a measure, which is the pushforward of Lebesgue measure on the parameterizing interval $[0,1]$. Note that no matter the value of $\rho$, the law of $(X_t,Y_t)$ is unchanged if we swap the roles of $X_t$ and $Y_t$, which (by contrast) is very much not the case for the $(X_t, Y_t)$ pair used to define the Brownian map.

Moore in 1925 gave a very general criterion for determining when a topological quotient space is topologically a sphere \cite{moore1925concerning} and as explained in \cite{matingoftrees, finitetree} one can verify directly that these criteria are satisfied in our setting, so that the object obtained by ``gluing the two trees together along their outer boundaries'' is indeed a topological sphere.  The {\em peanosphere} does not {\em a priori} come with a simple metric space structure like the Brownian sphere does, because neither the tree nor the dual tree is a tree of geodesics. On the other hand, both the tree and the dual tree (as determined by $X_t$ and $Y_t$) can be viewed as metric measure spaces themselves, so lengths of arcs {\em within} these trees are well defined.

\subsection{Percolation and other simple kinds of decoration}

A major reason to study the {\em peanosphere} is that it helps one understand random $(M,T)$ pairs, where $M$ is a random planar map and $T$ is some kind {\em extra decoration} on $M$.  Indeed, most of the scaling limit results known for decorated random planar maps make use of the peanosphere construction in some way. But before we discuss that topic, let us say a few words about randomly decorated {\em deterministic} planar maps---i.e., let us suppose that $M$ is fixed to be a grid or lattice, and we are choosing some way to ``decorate'' $M$, e.g.\ with a coloring or spanning tree.

In statistical physics, it is often interesting to try to find the {\em simplest possible models} that exhibit behaviors (phase transitions, correlation decays, fractal patterns, etc.) that one might see in complicated real world systems. By understanding these models thoroughly and mathematically, one can hope to get a glimpse of {\em why} similar behaviors appear in more complex systems. We will describe a few of these very simple models here.

First, let us discuss {\em lattice percolation} which was introduced to the mathematics literature by Broadbent and Hammersley in \cite{broadbent1957percolation}.  In this model, one starts with a set of edges---or vertices or faces---and tosses an independent coin for each one to decide whether it is ``open'' or ``closed.'' This is easiest to visualize by computer. For example, the following Mathemetica code generates the picture below it:

\vspace{-.1in}
{\small
\begin{verbatim}n=40;Graphics[{Table[If[RandomInteger[1]==1,Line[{{i,j},{i,j+1}}]],{i,0,n},{j,0,n-1}], 
  Table[If[RandomInteger[1]==1,Line[{{i,j},{i+1,j}}]],{i,0,n-1},{j,0,n}]}]
\end{verbatim}
}
\vspace{-.2in}
\begin{center}
\includegraphics[scale=.4]{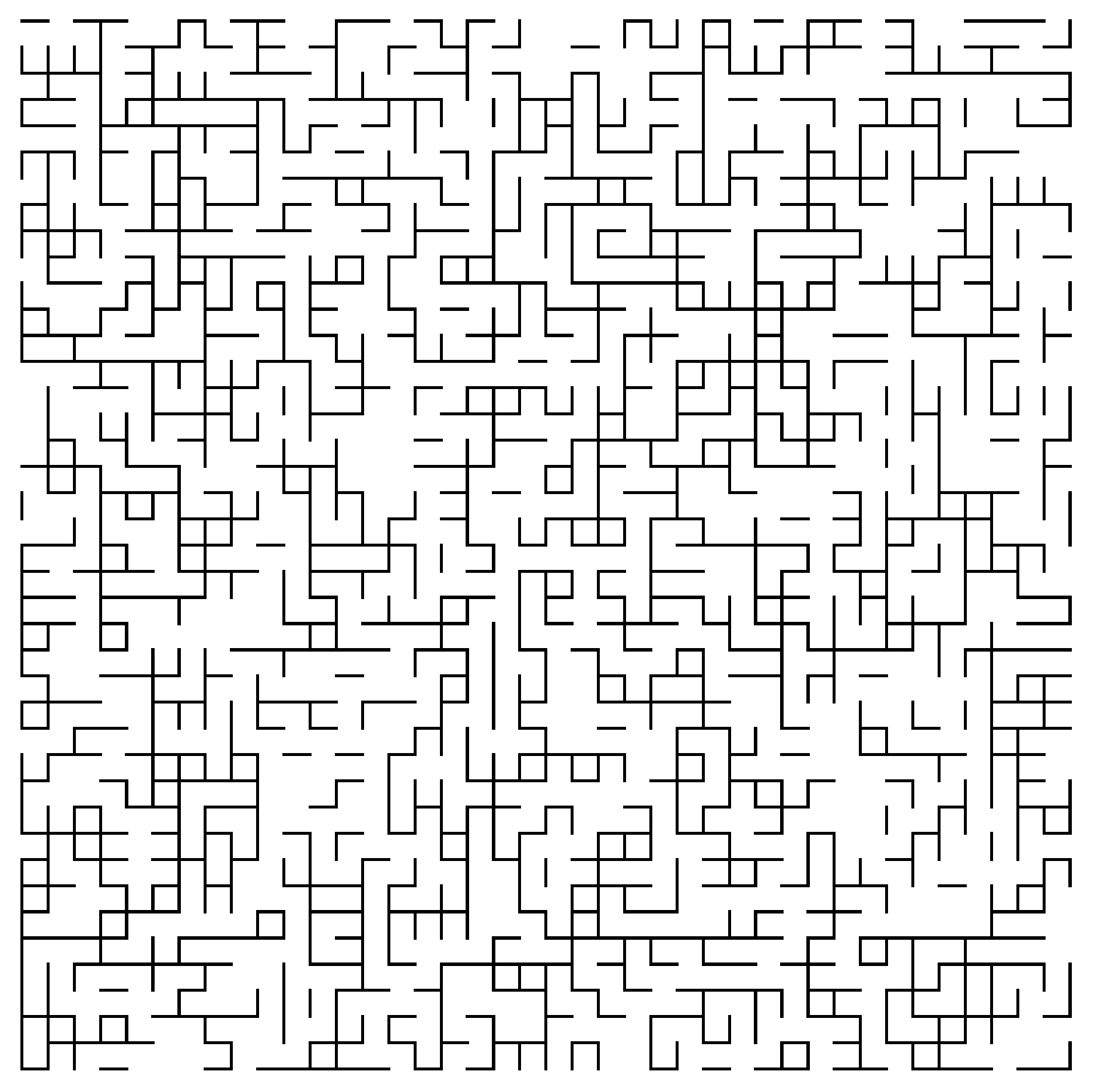}
\end{center}

The above is a $40$ by $40$ grid where one tosses a fair coin independently for each edge to decide whether to display it or not.  The displayed edges are the open edges; the edges not displayed are closed. A connected component of open edges is called a {\em cluster}. One can think of this as a simplistic model for a porous medium---with open edges representing pathways that current or fluid can flow through, and closed edges representing obstructions. The figure above is called {\em bond percolation} (with parameter $p=1/2$) because it is the edges (``bonds'') that are assigned to be open or closed.  An alternative is {\em site percolation} which independently assigns each vertex---or face if one uses the ``dual perspective''---to be open with probability $p$. By way of illustration, the following code generates the figure below it:

\vspace{-.1in}

{\small
\begin{verbatim}n=40; Graphics[Table[{If[(i-n)(j-n)==0,Blue, If[i j==0,Yellow,If[RandomInteger[1]==1,
Yellow,Blue]]],RegularPolygon[i{-Sqrt[3],-1}+j{-Sqrt[3],1},{1,0},6]},{i,0,n},{j,0,n}]]
\end{verbatim}
}
\vspace{-.2in}
\begin{center}
\includegraphics[scale=.6]{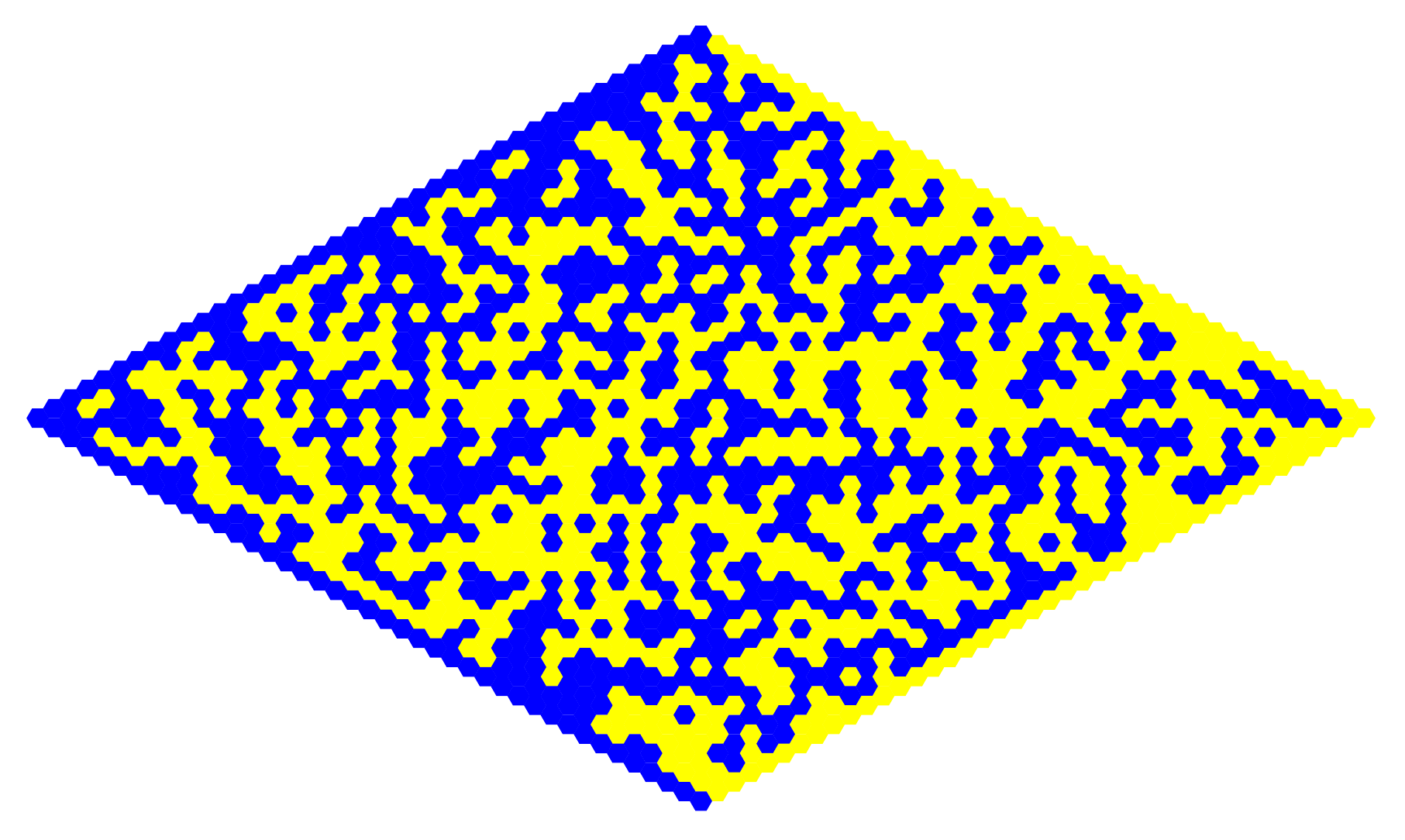}
\end{center}
\vspace{-.1in}
This figure represents percolation on the faces of a hexagonal grid.  For each interior face, an independent coin toss decides whether it is open (blue) or closed (yellow).  Here we have also imposed some deterministic boundary conditions, fixing the colors to be blue along one boundary arc and yellow on the complementary arc.  The ``blue-yellow'' edges (i.e.\ the edges on the boundary between blue and yellow) form many finite length cycles and one long path, which starts at the lower corner and ends at the upper corner. The long path is the boundary between the outermost yellow cluster (which contains half of the boundary) and the outermost blue cluster (which contains the other half of the boundary).  In the fine mesh limit, it converges in law to a random fractal curve called the chordal {\em Schramm-Loewner evolution with parameter $\kappa = 6$} (SLE$_6$) as proved in the breakthrough work of Smirnov \cite{smirnov2001critical} and Camia and Newman \cite{camia2007critical}, see also the expository paper by Sun \cite{sun2011conformally}. We will say more about SLE curves in Section~\ref{sec::SLE}.

Both bond percolation and site percolation have natural {\em variants} in which the edge or vertex subsets are sampled from different (non-uniform) probability measures.  One of the simplest such models is the {\em Fortuin-Kasteleyn (FK) random cluster model}, which was introduced by Fortuin and Kasteleyn in 1972  \cite{fortuin1972random, fortuin1972random2, fortuin1972random3}. It has a ``partition function'' (defined below) that is equivalent to the famous {\em Tutte polynomial} introduced by Tutte in 1954 \cite{tutte1954contribution}. See the historical accounts in \cite{crapo1969tutte, farr2007tutte} and the overviews by Bernardi and Welsh \cite{bernardi2008tutte, welsh2000potts} of the Tutte polynomial and its (surprisingly numerous) applications.

In the FK random cluster model, one begins with a finite graph $G=(V,E)$ and chooses a random $T \subset E$ (which need not necessarily be a spanning tree). The probability of a given $T$ is proportional to a {\em weight} given by $W_G(T) := q^{k(T)} w^{|T|}$ where $q$ and $w$ are positive constants, and $k(T)$ represents the number of {\em clusters} (i.e., connected components) of the graph $(V,T)$. Note that $k(t)$ is maximized and equal to $|V|$ if $T$ is empty. It is minimized and equal to $1$ if $(V,T)$ is connected.

If $q=1$, then the FK random cluster model is just ordinary bond percolation with $p = w/(1+w)$.  Generally, a large $w$ gives a bias toward $T$ having more edges, and a small $w$ gives a bias toward $T$ having fewer edges. If $q$ is very large, then the probability measure is biased in favor of $(V,T)$ being highly {\em disconnected} (i.e., having lots of components). If $w$ is close to zero, then the bias is in favor of $T$ being more {\em connected} (i.e., having few components).  The uniform spanning tree can be obtained as a limiting case of the FK random cluster model.  (Taking $w \to 0$ quickly forces $T$ to be a.s.\ connected; taking $p \to 0$ slowly makes $|T|$ as small as possible, and the connected $T$ that minimize $|T|$ are spanning trees.) The {\em partition function} of the FK random cluster model is defined to be the sum of the weights, taken over all $T \subset E$.
 $$Z_G(q,w) = \sum_{T \subset E} W_G(T),$$
 which again is equivalent to the Tutte polynomial (up to a certain coordinate change).  An important point to stress is that if $q=1$, the partition function $W_G(T)$  depends only on the number of edges in $G$.
 
For certain parameters, the FK random cluster model is also closely related to the hugely influential {\em Ising model} which was presented by Ising in 1925 but actually introduced earlier by Lenz in 1920, see e.g.\ the historical account at \cite{brush1967history}. In the Ising model, not all of the possible face colorings are equally likely; rather one {\em weights} the probability of each coloring by a constant to the number of adjacent pairs on which the colors disagree.

The decorations we have discussed so far (bond percolation, site percolation, uniform spanning tree, FK random cluster model, Ising model) are some of the very simplest statistical physics models. But there are others.  For example, one can decorate a graph $G$ with a {\em bipolar orientation}, which is a way of {\em orienting} the edges in $G$ so that only one vertex $v$ (the sink) has all incident edges oriented toward $v$, and only one vertex $w$ (the source) has all incident edges oriented away from $w$. The number of bipolar orientations is also encoded in the Tutte polynomial, see e.g.\ Bernardi's explanation in \cite{bernardi2008tutte}. One can also decorate $G$ with a so-called {\em Schnyder wood} or an instance of the Gaussian free field or a Brownian loop soup (all of which will be mentioned below). The list goes on.  We will not properly survey the literature on decorated planar maps here, but a few examples of work in this area include \cite{borot2011recursive, bouttier2004planar, BDGmobiles, bernardi2011counting, bernardi2009intervals, chen2020perimeter, bousquet2003degree, bernardi2007bijective, bouttier2005combinatorics}.

\subsection{Decorated random planar maps} \label{sec::decmap}


In this section we give some examples of how one can generate a  pair of spanning trees if one starts with another type of decoration. The figures below were used in \cite{burgerbijection} to explain how to generate a pair of trees from an instance of the FK random cluster model. By way of setup, one starts with a planar map $M$ in black, adds a vertex in the center of each face and draws a green edge connecting it to each vertex on the boundary of the face, thereby producing a green quadrangulation $Q$.  The planar map $M^*$ (which is dual to $M$) is obtained by connecting red-to-red within each green quadrilateral of $Q$, and is shown with dotted lines.  (Here $M$ itself is recovered from $Q$ by connecting blue to blue within each quadrilateral.)

\begin {center}
\includegraphics [width=3in]{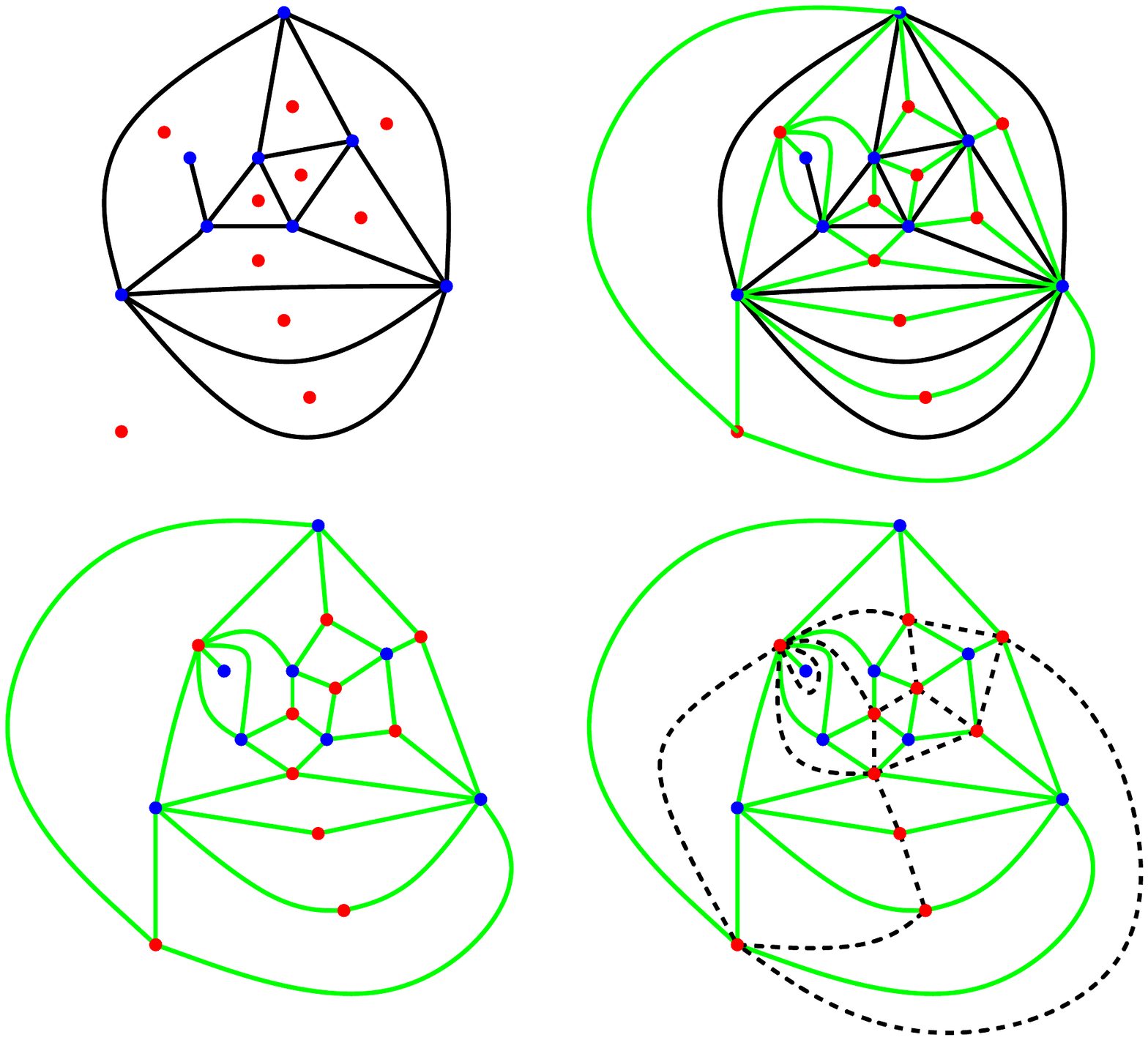}
\end {center}

Next, the figures below illustrate an edge subset $T$ in solid black (which happens to form a spanning tree) with the dual $T^*$ shown as dotted black lines.  (The dual consists of a red-to-red dotted line contained within each quadrilateral that {\em does not} have a blue-to-blue solid line.) The interface between the two trees is represented by a red curve that passes through all of the green edges one at a time, and starts at a location between the root and dual root (shown as large blue and red vertices).
\begin {center}
\includegraphics [width=3in]{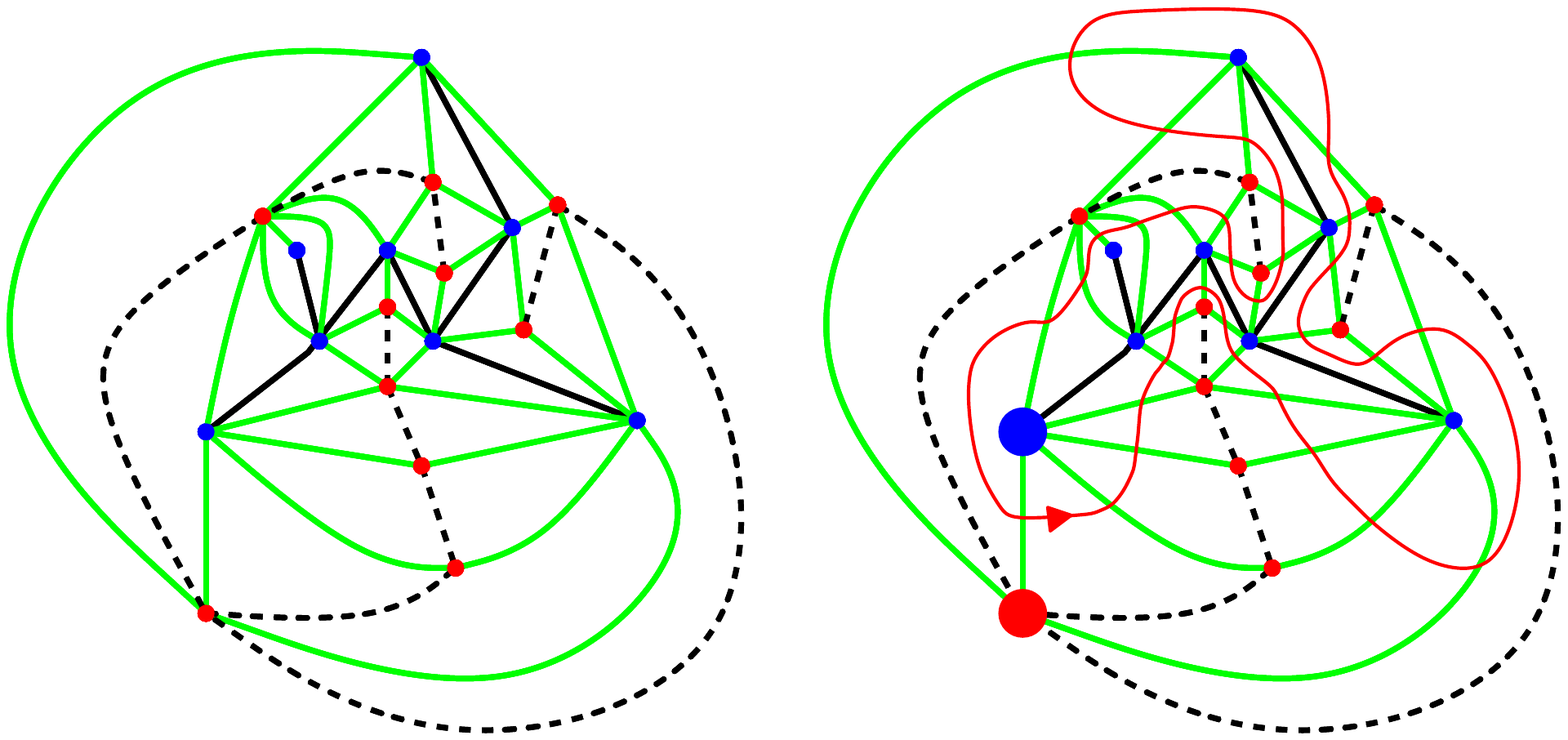}
\end {center}

Finally, the figure below on the left illustrates an example of $T$ (solid black lines) that is not a spanning tree (it contains cycles and there is one blue vertex that is not connected to any others).  When one draws the red curve (which passes through green lines but never through black solid/dotted lines) one finds that it is not connected---in fact it has five different components.  In the figure on the right, four ``swaps'' have been made, where the edge in a quadrilateral was changed to a dual edge, or vice versa.  The new edges/dual-edges produced by the swaps are colored in yellow. There is a simple algorithm (due to Bernardi in \cite{bernardi2008tutte}, see the author's explanation in \cite{burgerbijection}) for choosing which edges to swap. (Very roughly speaking: one tries to form a single red loop that crosses all the green edges one at a time, as above---but whenever the path is about to ``seal off'' a region without traversing it, one replaces a black blue-to-blue edge with a dotted dual red-to-red edge, or vice versa, to prevent this from happening.) Each swap reduces the number of red loops by one---and at the end, one is left with a single red loop, and a tree/spanning-tree pair.  This is called the {\em hamburger-cheeseburger} construction because of its relationship to a simple inventory accumulation model in which two kinds of products (``burgers'') are created and consumed  \cite{burgerbijection}.

\begin {center}
\includegraphics [width=5in]{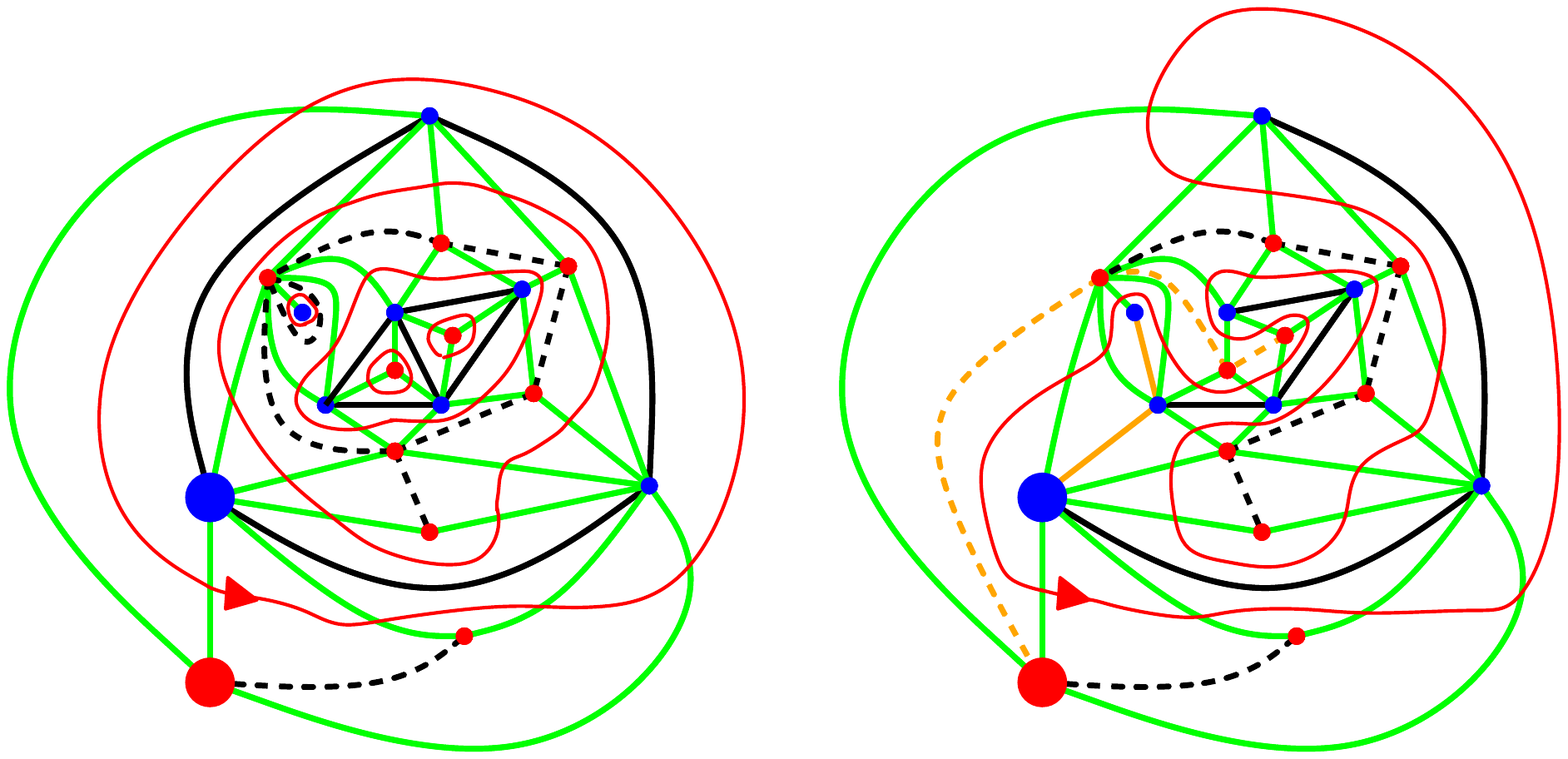}
\end {center}

At the end of the day, \cite{burgerbijection} uses some algebra (which is too detailed to explain here) to show that the lattice path encoding this pair of trees converges to a two-dimensional Brownian motion where the correlation coefficient $\rho$ depends on the FK-cluster parameter in a simple way. While the work in  \cite{burgerbijection} was done for infinite graphs, Gwynne, Mao and Sun have understood the scaling limits for the FK model on a finite random planar map, and have given much stronger forms of convergence, which also encode the structure of the loops themselves in a direct way \cite{gwynne2019scaling, gwynne2017scaling, gwynne2015scaling}. See also \cite{berestycki2017critical, chen2017basic} for more results about critical FK random maps. Since then, work by Bernardi, Holden, and Sun has explained that a simpler, tailor-made bijection applies for site-percolation-decorated random triangulations,  in which case $\rho = .5$.

With site or bond percolation, the partition functions depend only on the number of vertices or edges, and once one conditions on this the marginal law of $M$ is the same as in the undecorated model \cite{bernardi2018percolation}. This suggests that the peanosphere corresponding to $\rho=.5$ should be somehow equivalent to the Brownian map (together with extra randomness encoding the percolation structure on the Brownian map). And this is indeed the case---see Section~\ref{sec::relationships}.

Other bijections have been found for bipolar orientations and Schnyder woods, and so-called active spanning trees \cite{bipolar, li2017schnyder, biane2021mating, gwynne2018active, bernardi2009intervals}. All of these models have different $\rho$ values. For an extreme case, note that if we naively set $\rho=1$, then we are effectively gluing two identical Brownian trees to each other, and we simply obtain the same Brownian tree back.  On the other hand, there is a way of taking the $\rho \to 1$ limit that produces a non-trivial object, as shown in a remarkable paper by Aru, Holden, Powell and Sun \cite{aru2021mating}, which shows that there is a non-trivial construction that is far nicer than anyone could have expected, see also \cite{aru2019critical}.  (The $\rho=-1$ limit is expected to correspond to an ordinary Euclidean sphere.)

As a short preview for the experts, let us now note that there is a general relationship between $\rho$ and the $\kappa$ parameter from SLE theory, a parameter we will say more about in Section~\ref{sec::SLE}. Using this relationship, one finds for example that the bipolar orientation model corresponds to SLE$_{12}$ and the Schnyder wood model to SLE$_{16}$. This is because (see Section~\ref{sec::relationships}) there is a natural way to ``conformally map'' the peanosphere onto the Euclidean sphere, and when this is done the space-filling Peano curve becomes a space-filling form of SLE$_\kappa$  (see Sections~\ref{sec::SLE} and \ref{sec::CLE}) with $\kappa>4$, where $\rho = -\cos (4\pi/ \kappa)$ \cite{matingoftrees}. Here $\kappa$ ranges from $4$ to $\infty$ as $\rho$ ranges from $1$ to $-1$.  Let us note that the peanosphere surface itself turns out to be equivalent to an LQG sphere (see Section~\ref{sec::lqgsphere}) with relevant parameters satisfying $\gamma^2 = 16/\kappa$ and $Q = 2/\gamma + \gamma/2$ and $d = 25-6Q^2$. The special value $\rho = .5$  corresponds to $\kappa=6$ (and $\gamma = \sqrt{8/3}$ and $d = 0$). We will say more about these quantities in subsequent sections.

\subsection{Computing the scaling exponent}

As mentioned above, it turns out that if we set $\rho=.5$ then the peanosphere is a random surface that is equivalent to the Brownian sphere (decorated by a percolation model that does not change the law of the surface itself) although this is far from obvious at first glance. In this case, the two continuum trees are somehow midway between perfectly uncorrelated (as in the Mullin limit) and perfectly correlated. Even before constructing that bijection, it is not hard to argue heuristically that $.5$ is the only correlation coefficient that is consistent with the $C \beta^{n} n^{-7/2+k}$ formula, just as $0$ is the only correlation coefficient consistent with  $C \beta^{n} n^{-4+k}$. In general, one can derive a relationship between $\rho$ and the $b$ in $C \beta^{n} n^{-b+k}$.

We won't give the full details here, but let us sketch the rough idea, just to show that the relationship between $b$ and $\rho$ is not something mysterious. What is the probability that a simple random walk on $\mathbb Z^2$, started at the origin, remains in the wedge of angle $\theta$ until time $n$, and is back at the origin at time $n$? By way of illustration, part of a wedge of angle $\pi/4$ is shown below on the left (below the generating Mathematica code) along with its conformal image under the map $z \to z^4$.

\vspace{-.1in}
{\small
\begin{verbatim}f[z_]=z^4;F[{x_,y_}]={Re[f[x+I y]],Im[f[x+I y]]};n=40;
{ParametricPlot[Table[{{t(1-j/n)+j/n,j/n},{j/n,t j/n}},{j,0,n}],{t,0,1}],
ParametricPlot[Table[{F[{t(1-j/n)+j/n,j/n}],F[{j/n,t j/n}]},{j,0,n}],{t,0,1}]}
\end{verbatim}
}
\vspace{-.2in}
\begin{center}
\includegraphics[scale=.2]{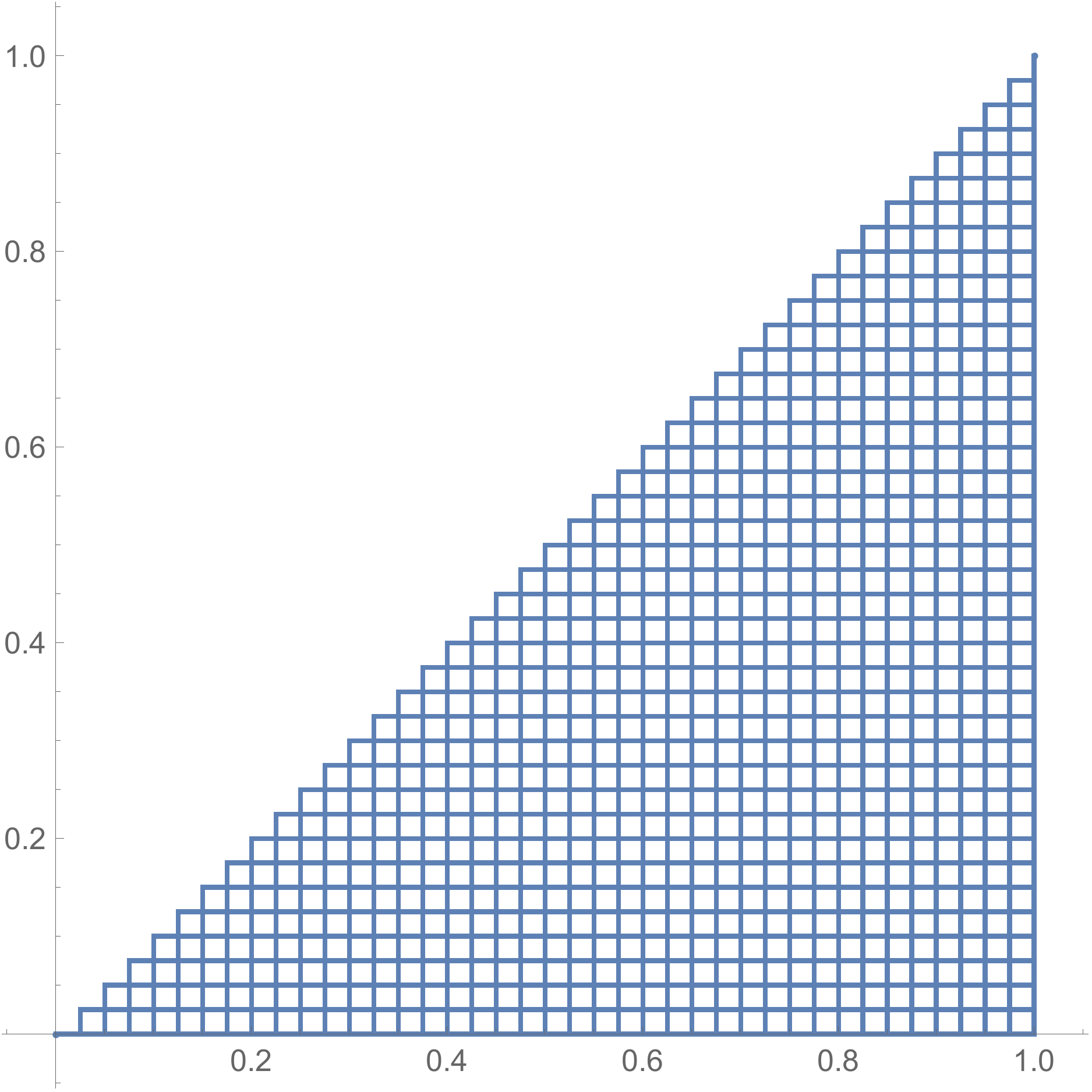} \includegraphics[scale=.3]{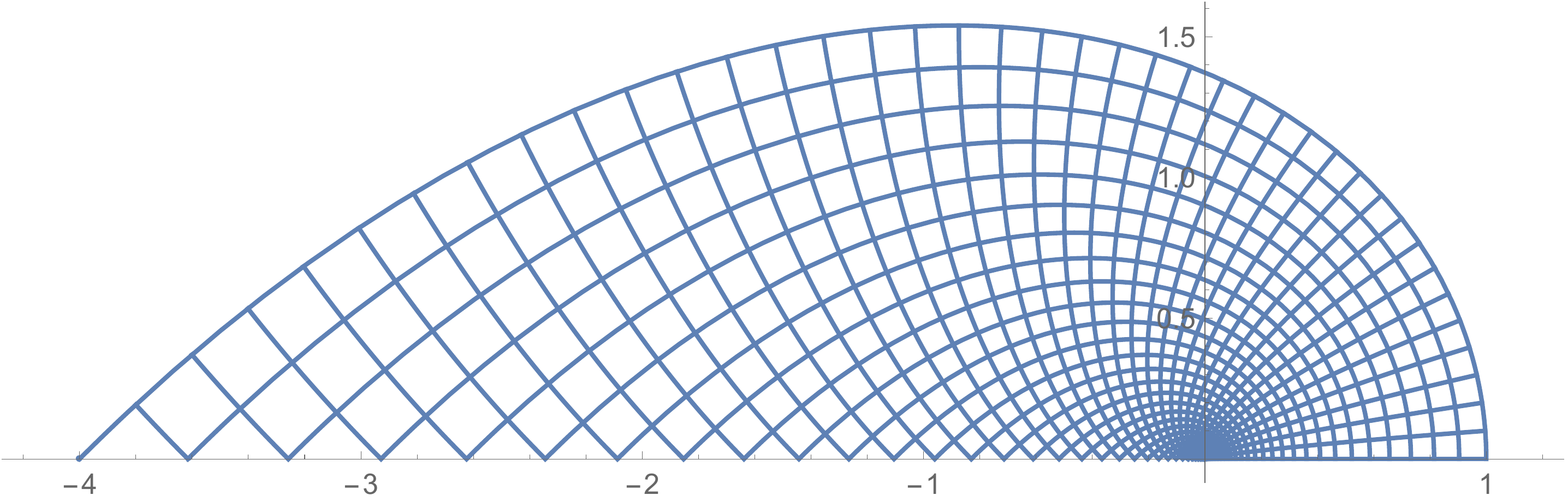}
\end{center}
Here is the heuristic. (References to rigorous arguments along these lines are given in \cite{bipolar}.) Suppose one does a simple random walk in a wedge of angle $\theta$, with lattice size $1/\sqrt{n}$, and we want the probability of the walk starting at the origin and returning in $n$ steps (without leaving the wedge).  Using the fact that the random walk approximates Brownian motion, and Brownian motion is conformally invariant under the map shown above, it appears that the odds of escaping to a macroscopic distance within $n/3$ steps should be of order $(1/\sqrt{n} )^{\pi/\theta} = n^{-\pi/(2\theta)}$.  The same should hold if one selects the walk increments in reverse (starting at time $n$) for order $n/3$ steps; and then there is an order $1/n$ chance that the middle third lines up correctly, so overall one expects a probability of order $n^{-\pi/\theta - 1}$.

The same should also hold if the simple random walk is replaced by another walk that has Brownian motion as a scaling limit. If $\rho \not = 0$ then it is necessary to {\em stretch} or {\em squash} the grid by some amount (in the $(1,-1)$ direction) in order to produce a walk with Brownian motion as a scaling limit, and one can easily work out the angle $\theta$ as a function of $\rho$.

The formula  $n^{-\pi/\theta - 1}$ corresponds to $n^{-b+k}$ where $k=1$ and $b = \pi/\theta + 2$ so that $\theta = \pi/(b-2)$.  If we have $b=4$ then $\theta = \pi/2$, which makes sense since the Mullin bijection corresponds to a simple random walk in a quadrant. More generally, if we have $X_t = \tilde B_t + a B_t$ and $Y_t=\tilde B_t - a B_t$ (where $B_t$ and $\tilde B_t$ are independent standard one-dimensional Brownian motions) we have a correlation coefficient $\rho = (1-a^2)/(1+a^2)$ where $a$ is the factor by which the space is scaled in the $(1,-1)$ direction.  Then $\theta$ is the angle obtained when we squash the standard positive quadrant by a factor of $a$ in that direction. Precisely, $1/a = \tan(\theta/2)$ so $a = \cot(\theta/2)$ and $\rho = (1-\cot(\theta/2)^2)/(1+\cot(\theta/2))^2 = -\cos(\theta) = -\cos\bigl( \frac{\pi}{b-2}\bigr) $. Plotting this relationship, we see that the $b$ value increases from $3$ to $\infty$ as $\rho$ decreases from $1$ to $-1$.
\begin{center}
\includegraphics[scale=.3]{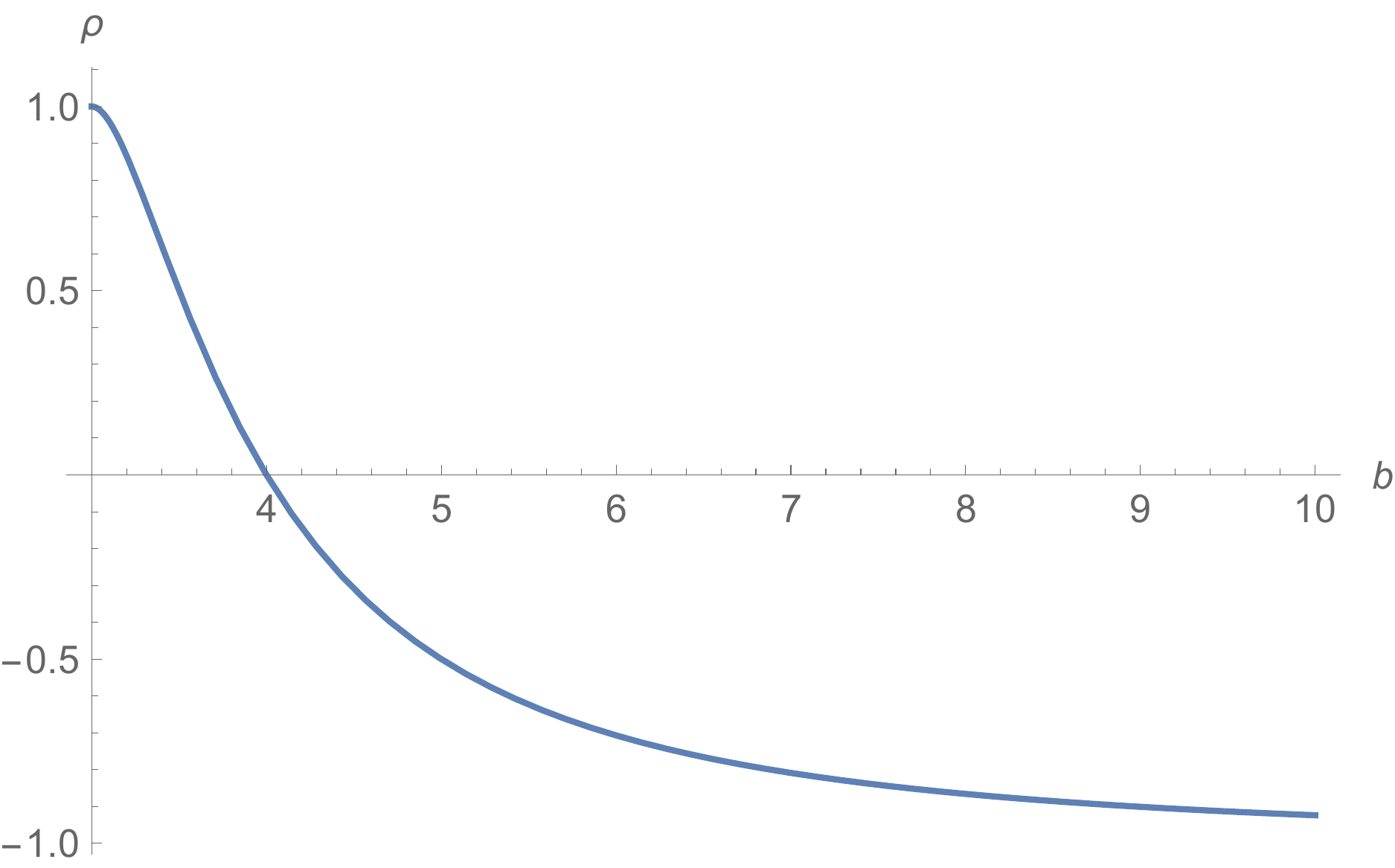}
\end{center}

\subsection{Schramm-Loewner evolution} \label{sec::SLE}

Schramm-Loewner evolution is by itself a large subject. It has been the topic of several other sets of ICM lecture notes, see e.g.\ the ICM notes by Duplantier, Lawler, Schramm, Smirnov, and Werner \cite{duplantier2014ICM, smirnov2007towards, smirnov2010discrete, werner2006conformal, lawler2018conformally, schramm2011conformally}. Longer expository introductions to the subject include e.g.\ \cite{berestycki2014lectures, lawler2004introduction, werner2004wendelin, kager2004guide}.

Here we will briefly explain what SLE curves are. First let us present the axioms that motivated the definition. In 1999 \cite{schrammsle}, Schramm set out to construct---for any simply connected domain $D$ with boundary points $a$ and $b$---a random non-self-crossing chordal curve $\eta$ connecting $a$ and $b$. Schramm insisted that the definition of an SLE curve have two properties. First, the definition had to be {\em conformally invariant} meaning that if $\psi$ is a conformal (i.e.\ analytic and one-to-one) map taking $D$ to a domain $\psi(D)$ then the image of $\eta$ under $\psi$ should have the law of an SLE in $\psi(D)$ from $\psi(a)$ to $\psi(b)$ (up to a time change).  Second, the path should be {\em Markovian} in the sense that {\em given} $\eta$ up to a stopping time $\tau$, the conditional law of {\em the rest} of $\eta$ is (up to a time change) that of an SLE in $D \setminus \eta\bigl([0, \tau] \bigr)$ from $\eta(\tau)$ to $b$.  Schramm showed that there was only a one-parameter family of ways to define SLE if one insists on these properties. Schramm indexed this family by a parameter $\kappa \in [0, \infty)$.

For completeness, let us now give Schramm's more explicit definition of SLE (though we won't say too much more about it here). By conformal invariance, it is enough to define the law of $\eta$ for one domain and one pair of boundary points.  It turns out to be convenient to work with the upper half plane $\mathbb H \subset \mathbb C$ with $a=0$ and $b =\infty$. For any time $t$, we define the function $g_t$ to be the unique conformal map from the unbounded component of $\mathbb H \setminus \eta([0,t])$ to $\mathbb H$ that satisfies $$\lim_{|z| \to \infty} g_t(z) -z = 0.$$  Schramm defined SLE in a rather indirect way: namely, he constructed the analytic functions $g_t$, and then used these functions to deduce what $\eta$ must be. The $g_t$ are defined by setting $g_0(z)=z$ and then requiring that for any fixed $z \in \mathbb H$, the ODE
$$\partial_t g_t(z) = \frac{2}{g_t- W_t}$$
is satisfied up until the smallest time $t$ at which $z$ is hit (or cut off from infinity by) the curve $\eta([0,t])$, where $W_t := B_{\kappa t}$ is a standard Brownian motion sped up by a factor of $\kappa$. This requirement determines the functions $g_t$ which in turn determine $\eta$.

In some sense, the larger $\kappa$ is (and hence the faster $W_t$ is moving up and down) the ``windier'' the curve becomes. In fact, Rohde and Schramm showed in \cite{rohde2005basic} that $\eta$ is a.s.\ a simple curve when $\kappa \in [0,4]$, that $\eta$ a.s.\ hits (but does not cross) itself when $\eta \in (4,8)$ and that $\eta$ is a.s.\ a {\em space-filling} curve when $\kappa \geq 8$. Beffara showed that the Hausdorff dimension of the range of $\eta$ is linear in $\kappa$, a.s\ given by $\min\{ 1+\kappa/8, 2 \}$ \cite{beffara2008dimension}. (Beffara builds on a related dimension calculation in \cite{rohde2005basic}.) In particular, since the scaling limit of the percolation picture above corresponds to $\kappa = 6$ the Hausdorff dimension of the curve is almost surely $1+6/8 = 7/4$.

At least on the square lattice, scaling limits of critical Ising model interfaces and the FK-Ising interfaces (which correspond to the FK cluster model for a particular choice of parameters) have scaling limits given by SLE$_3$ and SLE$_{16/3}$ respectively, as shown in a remarkable series of papers by Smirnov and co-authors including Chelkak, Duminil-Copin, Hongler, and Kemppainen \cite{chelkak2014convergence, smirnov2010conformal, chelkak2012universality}.  Loop-erased random walks and uniform spanning tree boundaries have scaling limits given by forms of SLE$_{2}$ and SLE$_8$ respectively, as shown in a remarkable paper by Lawler, Schramm and Werner \cite{lawler2011conformal}. Level lines of the so-called Gaussian free field are given by SLE$_4$ as shown by Schramm and the author \cite{dgffcontour, continuumgff}.

\subsection{Conformal loop ensembles} \label{sec::CLE}

What happens if we consider a percolation model with all blue boundary conditions, and then consider the collection of {\em all} of the loops that form boundaries between blue and yellow regions?
\vspace{-.2in}
\begin{center}
\includegraphics[scale=.6]{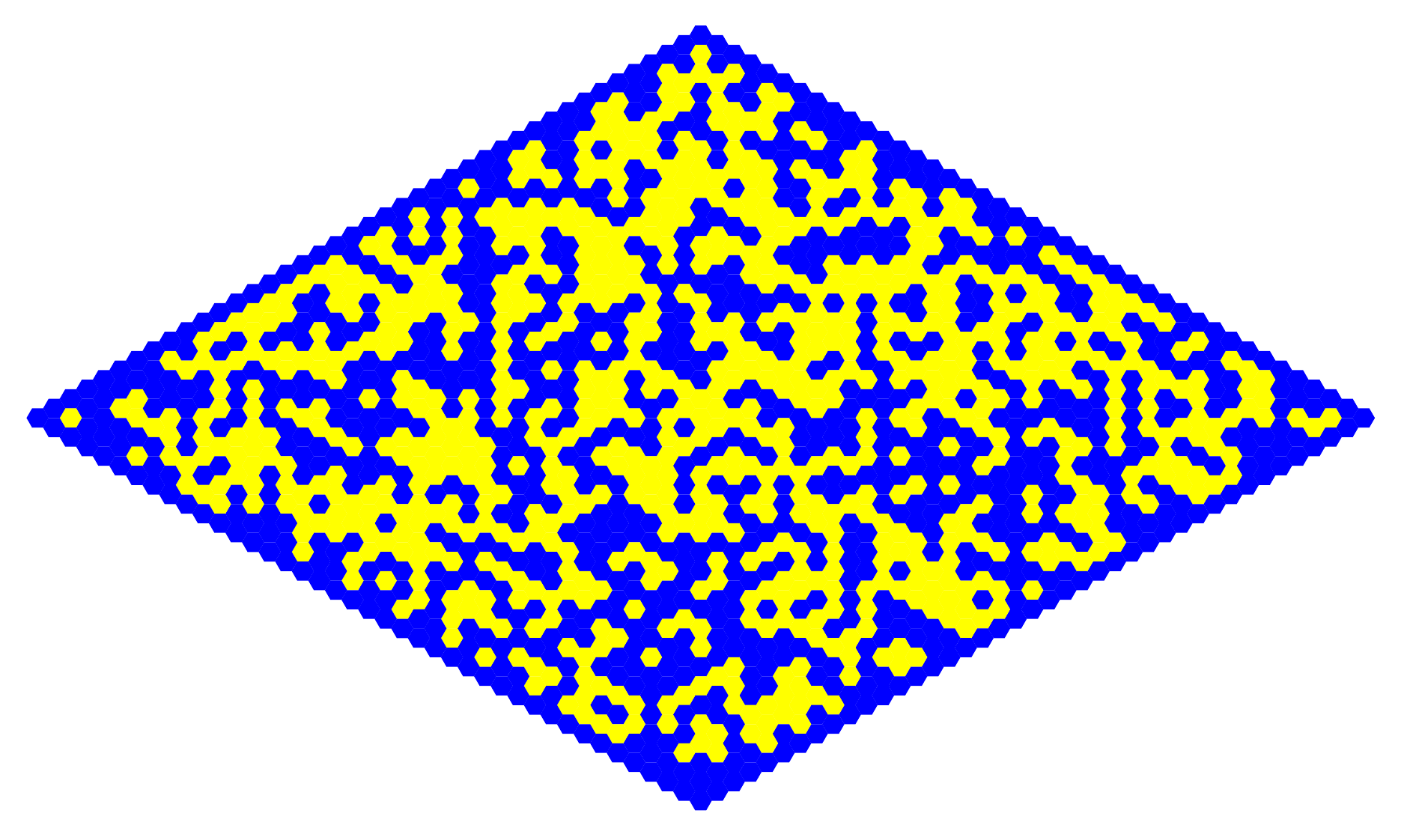}
\end{center}
\vspace{-.1in}

Camia and Newman showed in \cite{camia2006two} that in the fine-mesh scaling limit, these loops converge to a random collection of continuum loops, called a {\em conformal loop ensemble} (CLE) with parameter $\kappa=6$. Conformal loop ensembles are defined for any $\kappa \in (8/3,8]$, see the CLE construction by the author in \cite{explorationcle} and/or the axiomatic characterization of simple CLE loops by the author and Werner in \cite{sheffield2012conformal}.

We will not give a formal definition of CLE here.  But we wish to stress one point. As explained in  \cite{sheffield2012conformal}, it turns out that there is a natural way to use an instance of CLE$_\kappa$ for $\kappa \in (4,8)$ to generate a {\em space-filling} version of SLE$_\kappa$ (even though ordinary SLE$_\kappa$ is not space-filling for $\kappa$ in this range) which is somehow the continuum version of the procedure used in Section~\ref{sec::decmap} to combine multiple red loops into one loop.  One can also reverse the procedure and recover the loops from the space-filling curve. Very roughly speaking, one tries to follow the CLE interfaces, but any time the curve separates a region from the target point, one violates the rules and fills up that region before continuing.  The space-filling curve divides space into a continuum tree-and-dual-tree pair. Although we are not giving details here, this implies that, at least in the continuum, if one wants to understand {\em loop-decorated} random surfaces one can equivalently try to understand {\em tree-decorated} random surfaces.  We remark that these continuum trees also have an interpretation as coalescing rays within a so-called {\em imaginary geometry} \cite{IG1, IG2, IG3, IG4, dubedat2009sle} where the dual tree corresponds to the coalescing tree of rays drawn in the opposite direction.

\section{Liouville quantum gravity sphere: a random Riemannian geometry} \label{sec::lqgsphere}

The {\em Liouville quantum gravity sphere} is a random sphere-homeomorphic space whose law depends on the single parameter $d$ (or $Q$ or $\gamma$) as mentioned earlier.  When $d=0$ the LQG sphere is called the {\em pure LQG sphere}.  The pure LQG sphere differs from the Brownian sphere in that the former {\em a priori} comes endowed with a conformal structure but no metric space structure, while the latter is {\em a priori} endowed with a metric space structure but no conformal structure. 

A recent summary of LQG surfaces and the Gaussian free field can be found in the lecture notes by Berestycki and Powell \cite{berestycki2015introduction, berestycki2021gaussian}. See also the Ast\'erisque summary of LQG surfaces and the KPZ formula by Garban \cite{garban2013quantum}, the {\em Notices of the AMS} overview article by Gwynne \cite{gwynne2020notices}, and the ICM proceedings article by Miller. A more recent set of ICM lecture notes  by Ding, Dub\'edat and Gwynne  gives an overview of recent works that have established a metric space structure for general values of $d < 25$ \cite{ding2021introduction}.

To summarize the latter point, recently various researchers (such as Ang, Basu, Bhatia, Ding, Dub\'edat, Dunlap, Falconet, Ganguly, Gwynne, Holden, Miller, Pfeffer, Remy, Sep{\'u}lveda and Sun)  have contributed to a spectacular international program to show that  LQG spheres with $d \not = 0$ can also be given a canonical metric measure space structure, and to prove some basic properties about the resulting random metric spaces \cite{basu2021environment, gwynne2021conformal, gwynne2019kpz, pfeffer2021weak, gwynne2020liouville, gwynne2020dimension, dubedat2021metric,dubedat2020weak, ding2020tightness,ang2020volume, ding2021uniqueness,ding2021distance}. The reader might start by looking at the metric constructions by Gwynne and Miller or by Ding and Gwynne \cite{gwynne2021existence, ding2021uniqueness}. Viewed as random metric spaces, general LQG spheres can thus also be viewed as generalizations of the Brownian map. On the other hand, we stress that when $d \not = 0$ there is no reason to believe that the geodesic-tree/dual-tree pair is described by anything as simple as the Brownian snake. That is, as far as we know, the simple Brownian map construction given in Section~\ref{sec::browniansphere} has no simple analog corresponding to $d \not =0$.

\subsection{The Gaussian free field}
Brownian motion is a natural random function from $\mathbb R$ to $\mathbb R$.  A generalization of Brownian motion called the {\em Gaussian free field} (GFF) is a random (generalized) function from $\mathbb R^d$ to $\mathbb R$ for any $d$. Our use of the term GFF in this paper will be limited to $d=2$.

There are many ways to define the Gaussian free field, see the author's survey \cite{GFFsurvey}. It is a scaling limit of discrete random functions from $\mathbb Z^2$ to $\mathbb R$ much as Brownian motion is the scaling limit of random functions from $\mathbb Z$ to $\mathbb R$. One particularly concise definition is as follows.

Fix a bounded planar domain $D$. If $f$ and $g$ are functions on $D$ whose gradients lie in $L^2$ then we can write $(f,g)_\nabla = (2\pi)^{-1} \int_D \nabla f(z) \cdot \nabla g(z) dz$ for the {\em Dirichlet inner product} of $f$ and $g$. Let $H(D)$ be the Hilbert space closure of the space of compactly supported smooth functions w.r.t.\ this inner product.  Then the {\em Gaussian free field} on $D$ with zero boundary conditions is the sum $\sum \alpha_i f_i$ where the $f_i$ are an orthonormal basis for $H(D)$ and the $\alpha_i$ are independent standard normal random variables (mean zero, variance one). The sum a.s.\ does not converge pointwise or in $H(D)$ but it a.s.\ does converge in the space of generalized functions (a.k.a.\ {\em distributions}) \cite{GFFsurvey}.

The GFF can also be defined as the Gaussian random distribution with covariance given by Green's function $G(x,y)$. Here $G(x, \cdot)$ is given by $-\log |x-\cdot|$ minus the harmonic extension of $-\log|x-\cdot|$ from $\partial D$ to $D$.  The symbols $h$ and $\phi$ are both commonly used to describe an instance of the Gaussian free field, depending on the context. We will use $\phi$ in this paper.

The GFF can also be defined on the whole plane, where one simply has $G(x,y) = -\log|x-y|$. In this setting the field $\phi$ is only defined up to additive constant.  But one may nonetheless write $\textrm{Cov} \Bigl( (\phi, f), (\phi,g) \Bigr) = \int \int-\log|x-y| f(x) g(y)dxdy$ as long as both $f$ and $g$ have mean zero. The integration-by-parts identity $$(f,g)_\nabla = \frac{1}{2\pi} \int \nabla f(z) \cdot \nabla g(z) dz = \frac{-1}{2\pi} f(z) \Delta g(z) dz = \frac{-1}{2\pi}(f, -\Delta g)$$ is frequently used. 

\subsection{Conformal parameterizations}
The LQG-sphere has a long history.  On the physics side, LQG surfaces come up in certain formulations of string theory and 2D quantum field theories based on the Einstein equations (which in two dimensions reduce to the very simple Liouville equation). This literature is rich and complex, with foundational contributions by Belavin, Br\'ezin, David, Di Francesco, Distler, Dorn, Duplantier, Eynard, Fateev, Itzykson, Kawai, Kazakov, Knizhnik, Kostov, Migdal, Otto, Parisi, Polchinski, Polyakov, Segal, Seiberg, Teschner, Witten, the Zamolodchikov brothers, Zinn-Justin, Zuber and many others. (This list is far from exhaustive.) We will not attempt to properly survey the physics literature in this paper, but we point the reader to the long list of references in \cite{lqgkpz} (or the articles cited in Section~\ref{sec::cft}) as a place to start.

On the mathematics side, one might begin with Gauss \cite{gauss2005general} who explained in 1827 how curvature could be understood as an intrinsic property of a two-dimensional surface, independently of how the surface was ``embedded'' in a higher dimensional space.  The Riemann mapping theorem (formulated by Riemann in 1851, proved by Osgood in 1900) and the more general Riemann uniformization theorem (conjectured by Klein in 1893, Poincar\'{e} in 1892, proved by Poincar\'{e} in 1907, Koebe in 1907) also play a central role \cite{gray1994history, walsh1973history}.

It is standard in differentiable geometry to define a surface (or two-dimensional manifold) by covering the surface with a ``chart'' of open sets that can each be diffeomorphically mapped to a planar domain.  Within one of these open sets, parameterized by pairs $(x,y)$, the ``metric'' can be written $A(x,y)dx^2 + B(x,y) dxdy + C(x,y)dy^2$. The parameterization is said to be {\em conformal} if $A=C$ and $B=0$, and the Riemann uniformization theorem implies that one can always find a conformal parameterization.  If we treat the parameterizing domain $U$ as a subset of $\mathbb C$ and write $z=x+iy$ then we can write the metric as $e^{\rho(z)}(dx^2+dy^2)$ and the associated area measure as $e^{\rho(z)}dz$ where $dz$ is Lebesgue measure on $U$.  (In some conventions the definition of $\rho$ may differ by a factor of two; the issue is whether $e^{\rho(z)}$ is interpreted as the length multiplier or the area multiplier.)

The Gaussian curvature is $-e^{-\rho(z)}\Delta \rho(z)$ so the integral of the Gaussian curvature over the set parameterized by a region $R$ is equal to the integral of $\Delta \rho$ over $R$.  In particular, the Gaussian curvature is zero if and only if $\Delta \rho = 0$ so that $\rho$ is harmonic. A function $\rho$ on a bounded domain $D$ is harmonic if and only if it minimizes $\int \Bigl(\nabla \rho(z) \cdot \nabla \rho(z)\Bigr) dz = \int \Bigl( -\rho(z) \Delta \rho(z)\Bigr)dz$ given its boundary conditions. More generally, the Gaussian curvature is equal to the constant $K$ if and only if  $\Delta \rho = - K e^{\rho}$. The latter equation is called {\em Liouville's equation} and was formulated by Liouville in 1838 \cite{liouville1838note}. A function with constant curvature minimizes $\int\Bigl(\nabla \rho(z) \cdot \nabla \rho(z)\Bigr) + K e^{\rho(z)}dz$, which is a linear combination of the Dirichlet energy $(\rho,\rho)_\nabla$ and the overall surface area. (Polyakov used the latter quantity to define the so-called {\em Liouville action}, see Section~\ref{sec::embeddedpolyakov}.)

If we want to choose a random perturbation of a flat metric, we need to find a random function $\rho$ that is in some sense a random perturbation of a harmonic function. Roughly speaking, one ``randomly perturbs'' a harmonic function by replacing the assertion that ``the Dirichlet energy of $\rho$ is minimal given its boundary values'' with ``the probability of a given $\rho$ is proportional to the exponential of minus the Dirichlet energy of $\rho$.'' Formally this means taking $\rho$ to be a constant $\gamma\in(0,2)$ times the Gaussian free field.  The induced area measure then takes the form $\nu_\phi := e^{\gamma \phi(z)} dz$.  This cannot make literal sense (since $\phi$ is a random distribution) but there are various ways to make sense of this through regularization.

The first rigorous construction of a random measure with the law of $\nu_\phi$ is due to H{\o}egh-Krohn \cite{hoegh1971general} who constructed the object for $\gamma \in [0,\sqrt{2}]$.  H{\o}egh-Krohn was motivated by earlier works in the quantum field theory literature that made sense of $V(\phi)$ where $\phi$ was the Gaussian free field and $V$ was a polynomial. The exponential of the free field (viewed as a quantum field theory) was studied in several papers over the subsequent decade, cited by prominent quantum field theorists such as Glimm, Jaffe, and Simon.

 In the 1980s Kahane derived a similar construction for all $\gamma \in [0,2)$ and used the term {\em Gaussian multiplicative chaos} to describe the associated random measures \cite{kahane1985chaos}. Kahane was motivated by the multiplicative cascades popularized by Mandelbrot. (See e.g.\ \cite{duplantier2014critical} for $\gamma = 2$.) Neither H{\o}egh-Krohn nor Kahane {\em interpreted} these random measures as the pullback to a planar domain of the area measure on a random surface parameterized by that domain.  The author's work with Duplantier \cite{lqgkpz} constructed the measure $\nu_\phi$ as a weak limit as $\epsilon \to 0$ of the random measures $\epsilon^{\gamma^2/2} e^{\gamma \phi_\epsilon(z)} dz$, where $\phi_\epsilon(z)$ is the mean value of $\phi$ on the ball boundary $\partial B_\epsilon(z)$, see also \cite{sheffield2016field}.  This construction showed that the measure $\nu_\phi$ is a function of the Gaussian free field $\phi$.

Computer visualization can help us generate some intuition about what the measures $\nu_\phi$ look like. The following Mathematica code generates the two figures shown below.

\vspace{-.15in}
{\small
\begin{verbatim}K = 8; fieldmultiplier = 1.5; squarefraction = .001;
phi=Re[Fourier[Table[(InverseErf[2 Random[]-1]+I InverseErf[2 Random[]-1])*If[j+k == 2,0,
     1/Sqrt[(Sin[(j-1)*Pi/2^K]^2+Sin[(k-1)*Pi/2^K]^2)]],{j,2^K},{k,2^K}]]]; 
CO = squarefraction Sum[M[[i,j]], {i,1,2^K}, {j,1,2^K}]; M=Exp[fieldmultiplier phi];
{ListPlot3D[phi],Graphics[Table[Table[If[Sum[M[[2^k m+i,2^k n+j]],{i,1,2^k},{j,1,2^k}]<CO,
    {Hue[k/8],EdgeForm[Thin],Rectangle[{2^k m, 2^k n},{2^k m+2^k,2^k n+2^k}]}],
  {m,0,2^(K-k)-1},{n,0,2^(K-k)-1}], {k,0,K-1}]]}  \end{verbatim}
}

\begin{center}
\includegraphics[scale=.32]{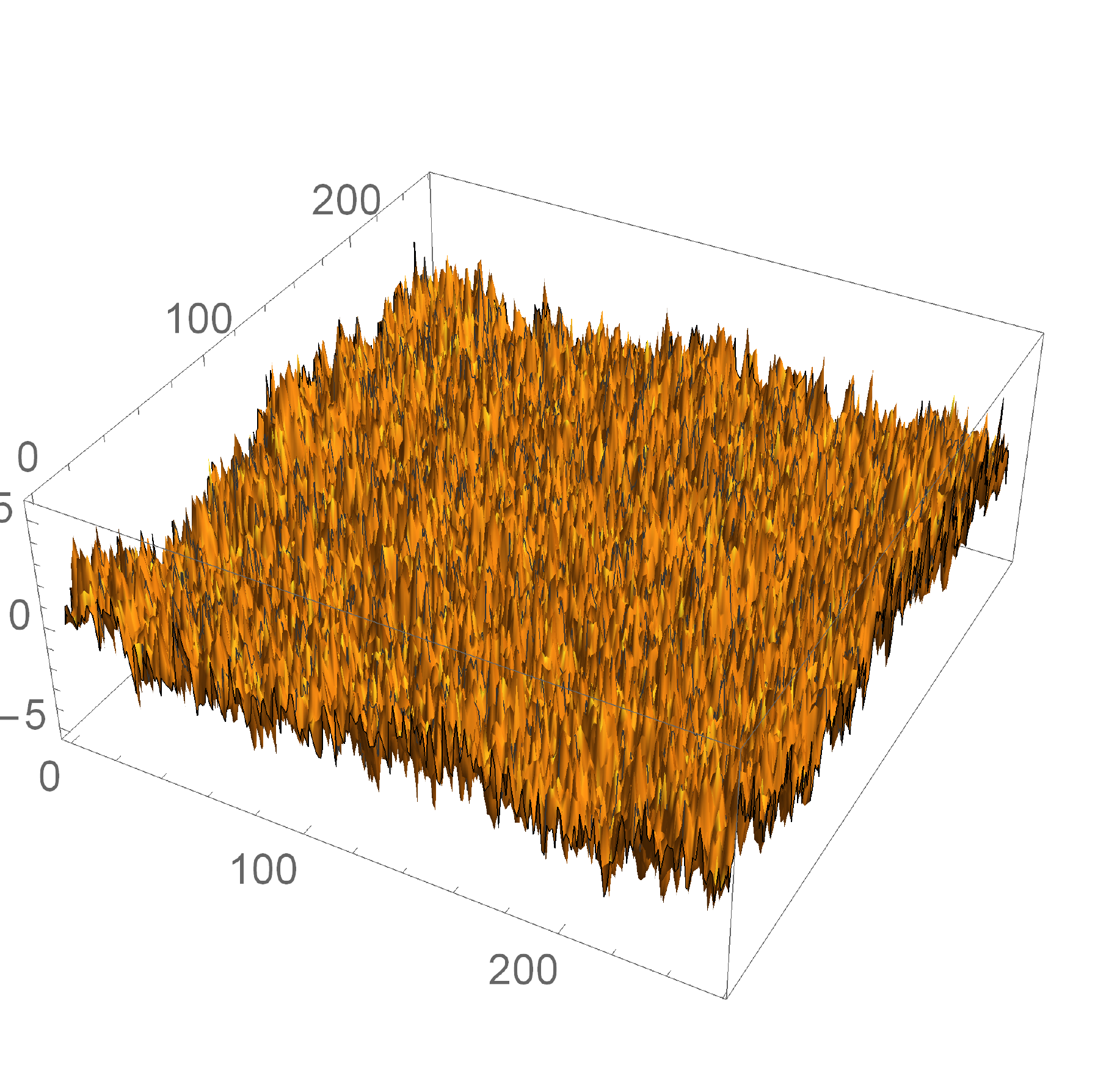} \includegraphics[scale=.28]{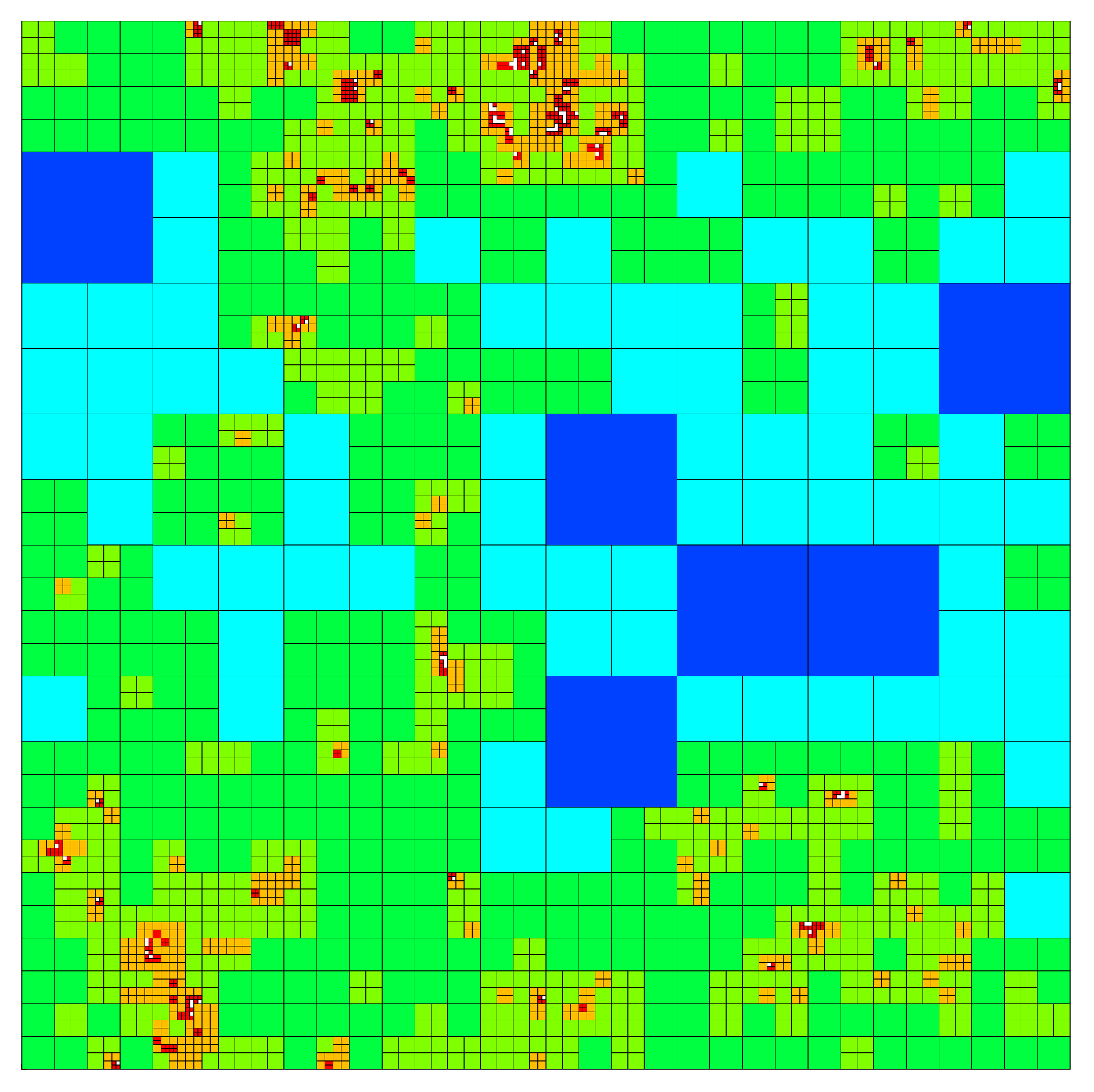}
\end{center}
The left figure is an instance $\phi$ of a discrete version of the Gaussian free field on the $256 \times 256$ torus (lines 2 and 3 encode the discrete GFF; see the explanation in \cite{GFFsurvey}). The exponential of {\em fieldmultiplier} times $\phi$ (called M in the code) describes a random measure on the torus. To obtain the picture on the right one starts with the whole square, then divides it into four equal squares, then divides each of those into four equal squares and so on, except that one stops dividing whenever one reaches a square where the area in the square is less than some constant cutoff (taken in the code to be {\em squarefraction} times the total area). The squares are colored according to their Euclidean size.  Each square $S$ shown in the picture above has area less than the cutoff---but its dyadic parent $S'$ must have area greater than the cutoff (as otherwise $S'$ would have been drawn after $S$, and would have covered $S$). In this sense one expects that all of the squares shown are ``about the same size'' in the random geometry.  The blue squares correspond to regions where $\phi$ is smaller on average (and the measure is less dense) while the orange and red squares correspond to high density areas. One can easily paste the code above into Mathematica and experiment with different variants: for example, shown below are the figures obtained by taking a fieldmultiplier of $.5$ (left) and $1$ (right) instead of $1.5$ above. The measure represented by the left figure is much closer to Euclidean measure.

\begin{center}
\includegraphics[scale=.35]{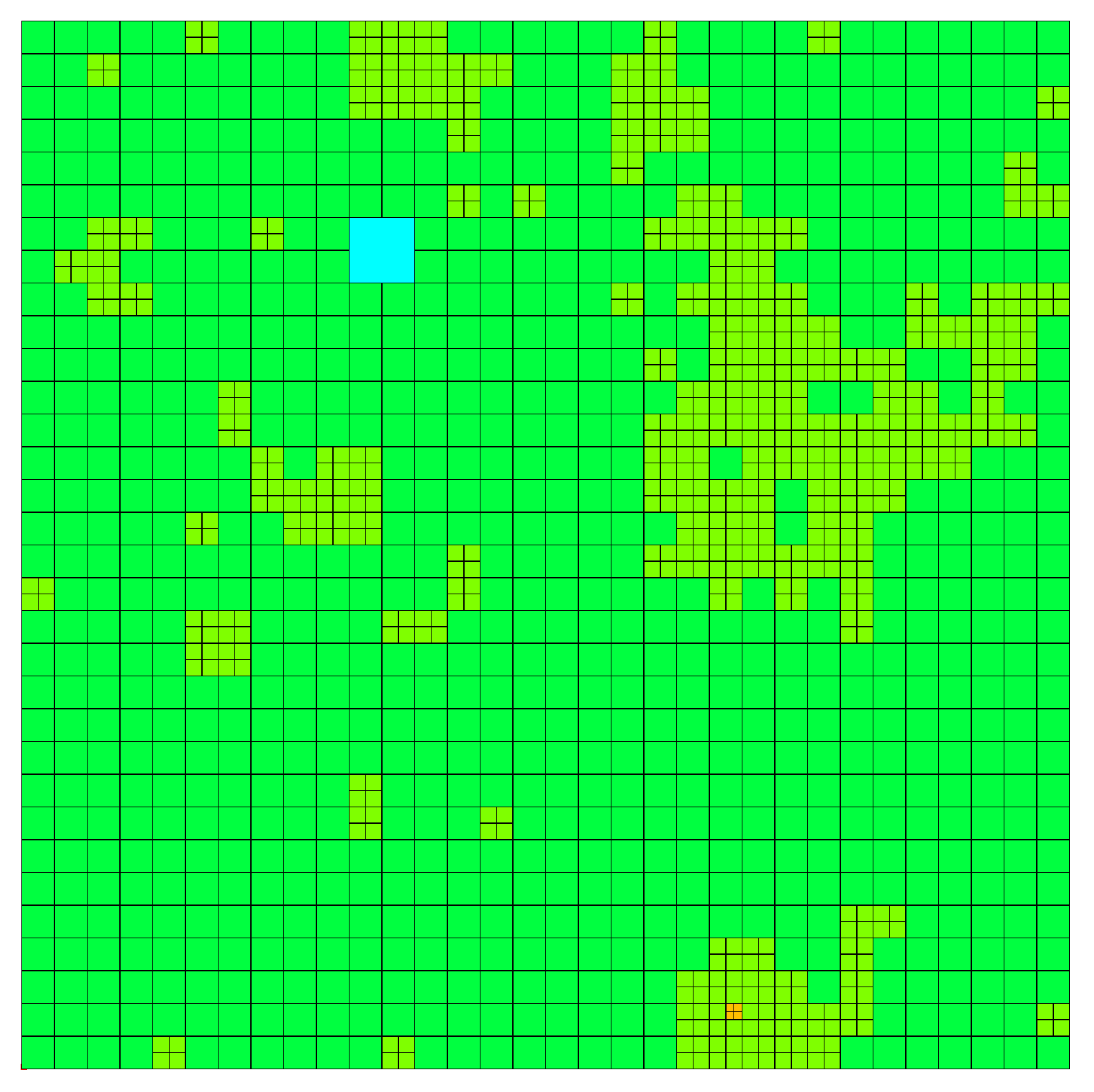} \includegraphics[scale=.35]{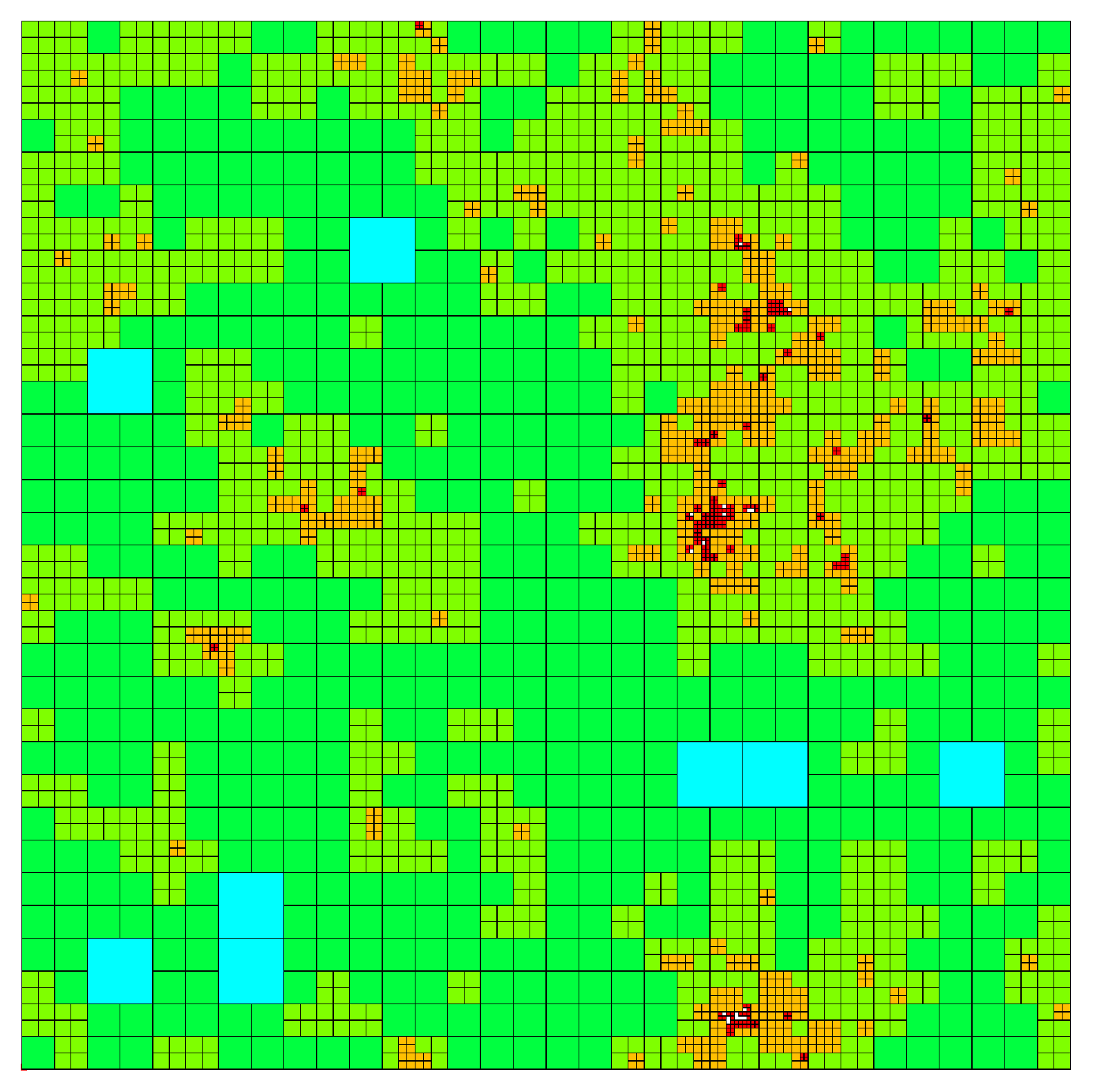}
\end{center}

\subsection{LQG surfaces}
We use the term {\em LQG surface} broadly to describe any surface whose area measure is $e^{\gamma \phi(z)}dz$ where $\phi$ has a law that is locally absolutely continuous w.r.t.\ the Gaussian free field. Formally an LQG surface is an equivalence class of pairs $(D,\phi)$, where $(D, \phi)$ and $(\tilde D, \tilde \phi)$ are equivalent if they are related as in the following diagram with $Q = \gamma/2 + 2/\gamma$.

\begin{center}
\includegraphics[scale=1]{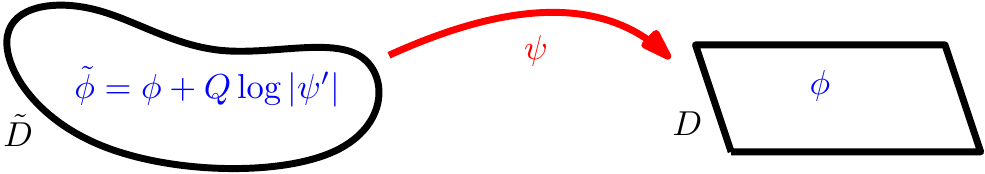}
\end{center}
\vspace{-.1in}
Using the definition $\nu_\phi := \lim_{\epsilon \to 0} \epsilon^{\gamma^2/2} e^{\gamma \phi_\epsilon(z)} dz$
one can show that the $\nu_{\tilde \phi}$ area measure on $\tilde D$ above must a.s.\ agree with the pullback of the $\nu_{\phi}$ area measure on $D$. Roughly speaking, this is because $\exp (\gamma Q \log |\psi'|) = |\psi'|^2 |\psi'|^{\gamma^2/2}$. The $|\psi'|^2$ is part of the usual change of variables formula, and the $|\psi'|^{\gamma^2/2} $ accounts for the fact that $\psi$ maps circles of diameter $\epsilon$ to (approximate) circles of size diameter $|\psi'|\epsilon$, which affects the $\epsilon^{\gamma^2/2}$ factor in the regularization.

In addition to the area measure $\nu_\phi$ it is possible to define a boundary length measure (in the case that $\phi$ is a GFF with free boundary \cite{lqgkpz}) or a fractal length associated to a line or an SLE curve or other fractal set, see \cite{lqggff}.  All of these measures, as well as the LQG {\em distance function} discussed above, are preserved by the coordinate change---which makes sense since $(D,\phi)$ and $(\tilde D, \tilde \phi)$ literally represent the same surface. The image under $\psi$ of the so-called {\em Liouville Brownian motion} in $D$ as defined in \cite{berestycki2015diffusion, garban2016liouville} is a Liouville Brownian motion in $\tilde D$.

Furthermore, if $x_1, \ldots, x_k$ are points in $D$, and we write $\tilde x_i = \psi(x_i)$, then we say that $(D, \phi, x_1, \ldots, x_k)$ and $(\tilde D, \tilde \phi, \tilde x_1, \ldots, \tilde x_k)$ represent the same ``LQG surface with $k$ marked points.'' Finally, let us stress that the above definition of LQG surface makes sense for any $Q\geq 0$, not only for $Q \geq 2$. (The values $Q \geq 2$ are those that have the form $Q= \gamma/2 + 2/\gamma$ for some $\gamma>0$, or equivalently those for which $d = 25-6Q^2<1$.) One just has to accept that the {\em area measure} is not well defined when $Q < 2$. However, as we mentioned earlier, the metric space structure (along with the fractal measure of some other sets) continues to be well defined when $Q<2$ (although the diameter of the surface becomes infinite when $Q<2$, which corresponds to $d>1$).

\subsection{Constructing the LQG sphere}

One simple way to define a unit-area LQG sphere is to consider a GFF $\phi$ on a simply connected bounded domain $D$ with boundary conditions given by a constant $C$, and to condition on $\nu_\phi(D) = 1$. As $C \to -\infty$ this object converges in law to the unit area LQG sphere (with the boundary somehow shrinking---in the metric sense---to a single marked point in the limit) as explained in \cite{zipper}.  Another approach involves starting with an infinite cylinder and two marked points (the cylinder's endpoints), as described in \cite{matingoftrees}.  Yet another approach involves starting with the infinite-volume Polyakov measure (see the next subsection) and ``pinning it down'' at three points, in a manner described in \cite{david2016liouville}. The paper by Aru, Huang and Sun \cite{aruhuangsunsphere} established the equivalence of the  approaches in  \cite{matingoftrees} and  \cite{david2016liouville}. (See also the alternate proof in \cite{ang2021integrability2} and the disk analog in \cite{cercle2021unit}.) A chart, presenting several equivalent definitions and the relationships between them, is included in \cite{finitetree}.

All of these approaches have analogous constructions that produce {\em unrestricted area} LQG spheres (instead of unit area LQG spheres).  In the $C \to -\infty$ construction above, for example, instead of conditioning on $\nu_\phi(D) = 1$ one can (for any $C$ value) multiply the measure by a constant factor to ensure that the measure assigned to $\phi$ distributions with $\nu_\phi(D) \in [1,2]$ is constant, and then take the vague $C \to -\infty$ limit.

\subsection{Polyakov's infinite measure on embedded LQG surfaces} \label{sec::polyakov}
The Polyakov measure is an {\em infinite} measure on the space of unmarked, unrestricted-area LQG spheres {\em embedded} in $\mathbb C$. As mentioned in the introduction, it corresponds to an unrestricted-area LQG sphere embedded in all possible ways (with the embedding chosen from Haar measure on the M\"{o}bius group). The measure is infinite for two reasons: first, we recall from the introduction that the law of the area for an unrestricted-area LQG sphere has the form $A^{-b}dA$ (or $A^{-b}e^{-\mu A}dA$ for some constant $\mu$ in the ``off-critical'' case) which is an infinite measure for the $b$ values that we will encounter (namely $b > 3$, recall Section~\ref{sec::peanospherescaling}). Second, the embedding is chosen from the Haar measure on the  M\"{o}bius group, which itself has infinite volume. Also, just to avoid confusion, let us clarify that it is the Polyakov measure on the space of surfaces --- not the area measure on any individual surface --- that is infinite. In the Polyakov measure, almost all surfaces have finite area (assuming $d \leq 1$ --- the total area is not defined if $d > 1$). Each embedded surface is described by a generalized function $\phi$, so the Polyakov measure can be viewed as an infinite measure on the set of generalized functions.

The Polyakov measure is usually defined in a slightly different way. It is presented as a way to make sense of the expression ``$e^{-S(\phi)}d\phi$'' where $S$ is the so-called {\em Liouville action}, which we will discuss below, see the presentation by David, Kupiainen, Rhodes and Vargas at \cite{david2016liouville} or the lecture notes by Kupianen at \cite{kupiainen2016constructive}. The simplest way to describe it (in the case $\mu=0$) is to say that $\phi$ is an instance of the zero-mean GFF on the sphere plus an independent constant $C$ chosen from the infinite measure $e^{-2Q c} dc$. To obtain the general-$\mu$ measure one simply weights the zero-$\mu$ measure by  $e^{-\mu A}$, where $A$ is the area of the LQG surface.

The fact that the Liouville action produces a measure on embedded LQG surfaces of the form described above (in particular, a measure that is  M\"{o}bius invariant) is counterintuitive at first glance, but it is carefully explained e.g.\ by Ang, Holden and Sun in \cite{ang2021integrability2}, see also the conformal invariance discussion in \cite{kupiainen2016constructive}. There is a certain miracle (related to what we will call ``semi-Gaussian'' measures) that makes everything work out.  Before we discuss the specifics, we will present a couple of simpler semi-Gaussian measures as a warmup.

\subsubsection{Semi-Gaussian measures} \label{sec:semigaussian}
A \textbf{semi-Gaussian} measure is a constant times $e^{-F(v)}dv$ where $F$ is a quadratic that is strictly convex in all directions except one.
For example, $\frac{1}{\sqrt{2\pi}} e^{y-x^2}dxdy$ is semi-Gaussian. It is a product of a normal measure $\frac{1}{\sqrt{2\pi}} e^{-x^2/2}dx$ and an infinite measure $e^{y}dy$. If we restrict this measure to any non-vertical line, we obtain a finite measure, which is a constant times a normal probability measure. If we restrict to any vertical line, we obtain an infinite measure.

In the figures below, the orange lines are level sets of the function $y-x^2$ (and hence also the function $\frac{1}{\sqrt{2\pi}} e^{y-x^2}dxdy$). In each figure, one can imagine ``sampling'' $(X,Y)$ from the infinite measure in two steps.  First one decides which blue line $(X,Y)$ belongs to.  (The blue lines in the three figures are the lines of slope $0$, $-1$ and $1$ respectively. By integrating the density function along the blue lines, we find that the law of the $y$-intercept of the blue line is given by a constant multiple of $e^{y}dy$.) Then {\em given that} one chooses the location {\em on} the blue line.  The {\em conditional law} of the location on the blue line is that of a Gaussian random variable centered at a point on the vertical green line (which is the point on the blue line where $y-x^2$ is largest).

\vspace{-.1in}
{\ssmall
\begin{verbatim}n = 4;{ParametricPlot[{Table[{t,j/n},{j,-2n,2n}],Table[{t,t^2+ j/n},{j,-2n,2n}],{0,t}},{t,-2,2},
PlotRange->{{-1,1},{-1,1}}],ParametricPlot[{Table[{t,-t+j/n},{j,-2n,2n}],Table[{t,t^2+j/n},
{j,-2n,2n}],{-.5,t}},{t,-2,2},PlotRange->{{-1,1},{-1,1}}],ParametricPlot[{Table[{t,t+j/n},
{j,-2n,2n}],Table[{t,t^2+j/n},{j,-2n,2n}],{.5,t}},{t,-2,2},PlotRange->{{-1,1},{-1,1}}]}
\end{verbatim}
}
\vspace{-.2in}
\begin{center}
\includegraphics[scale=.2]{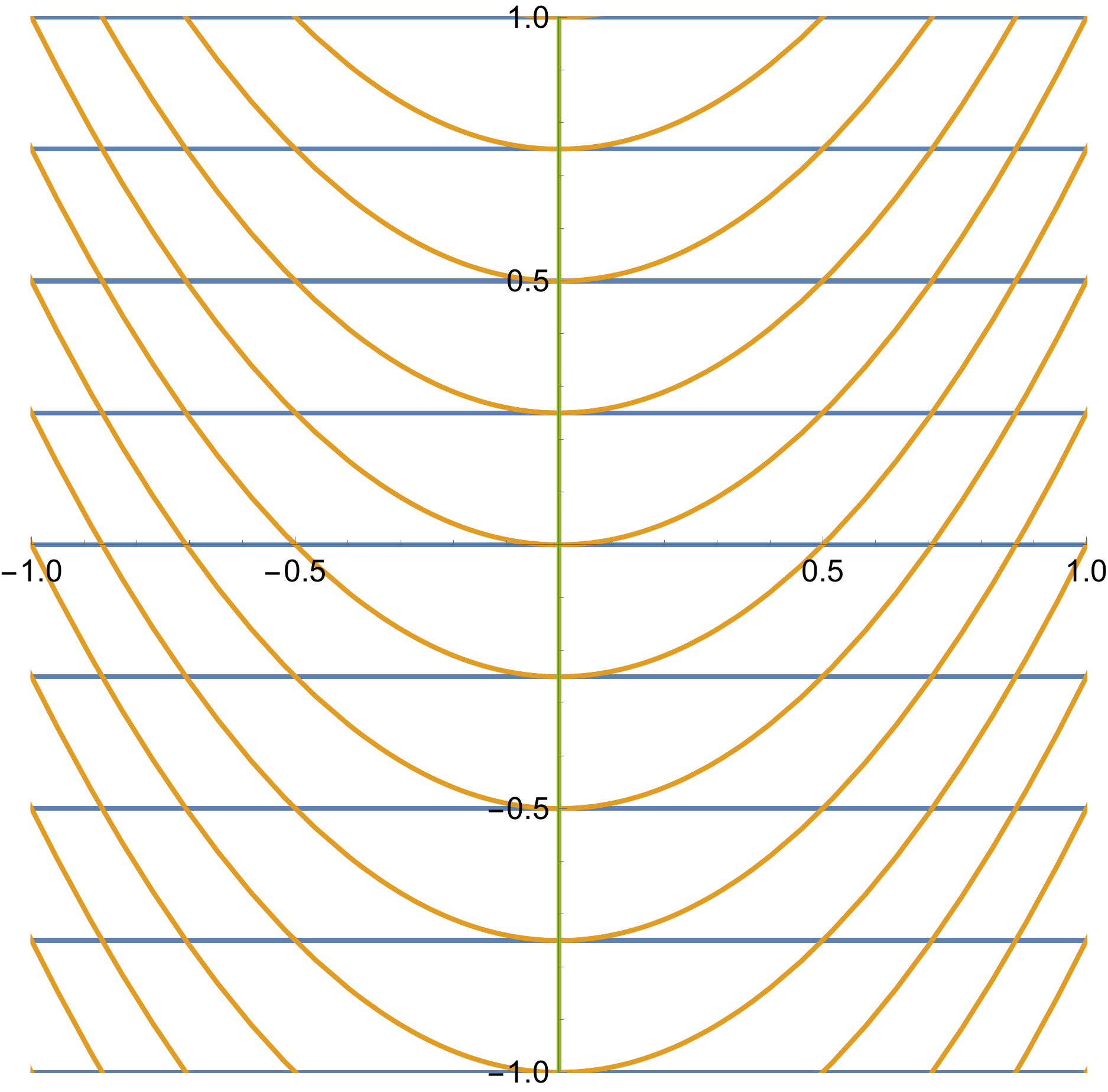} \includegraphics[scale=.2]{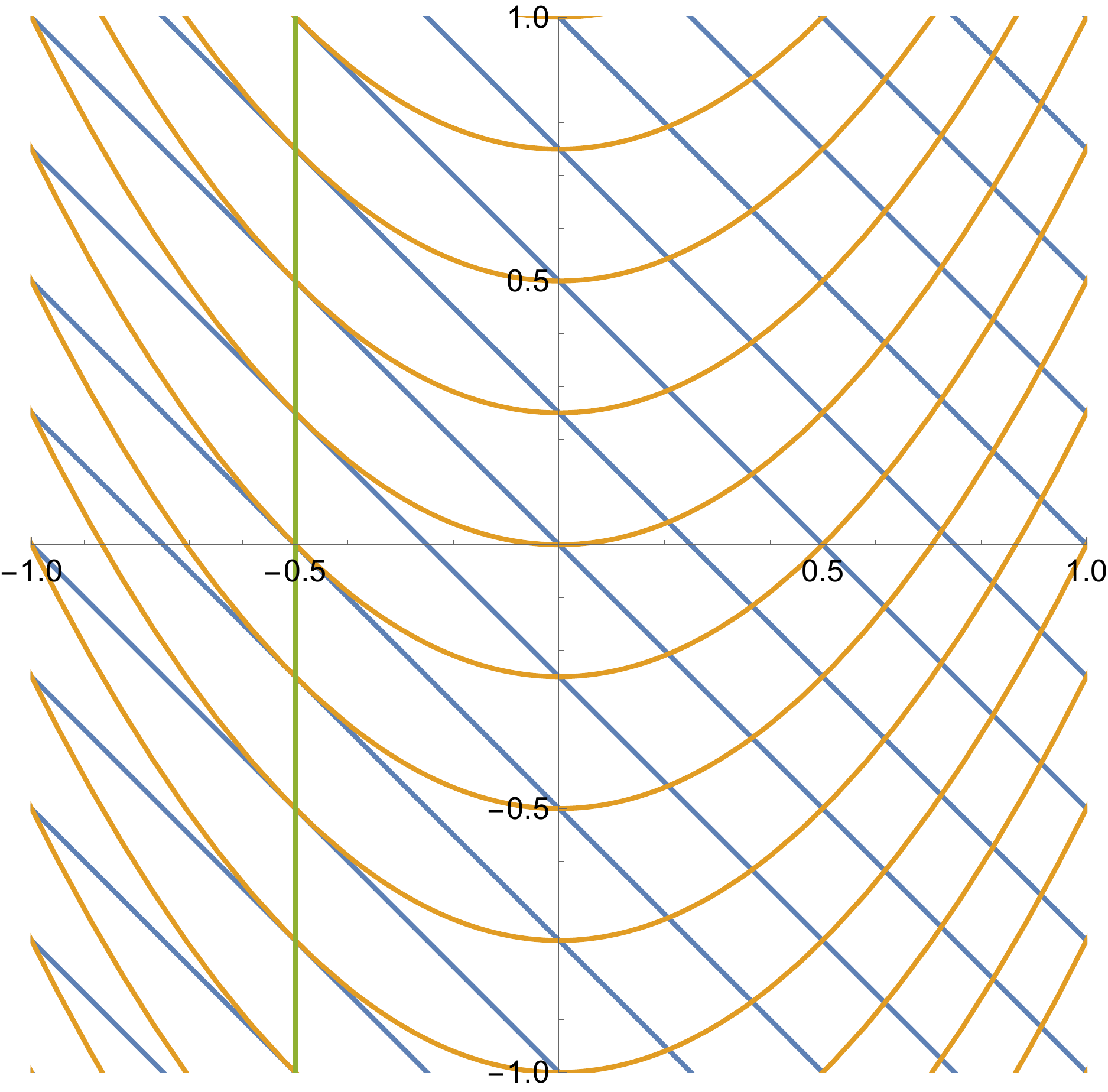} \includegraphics[scale=.2]{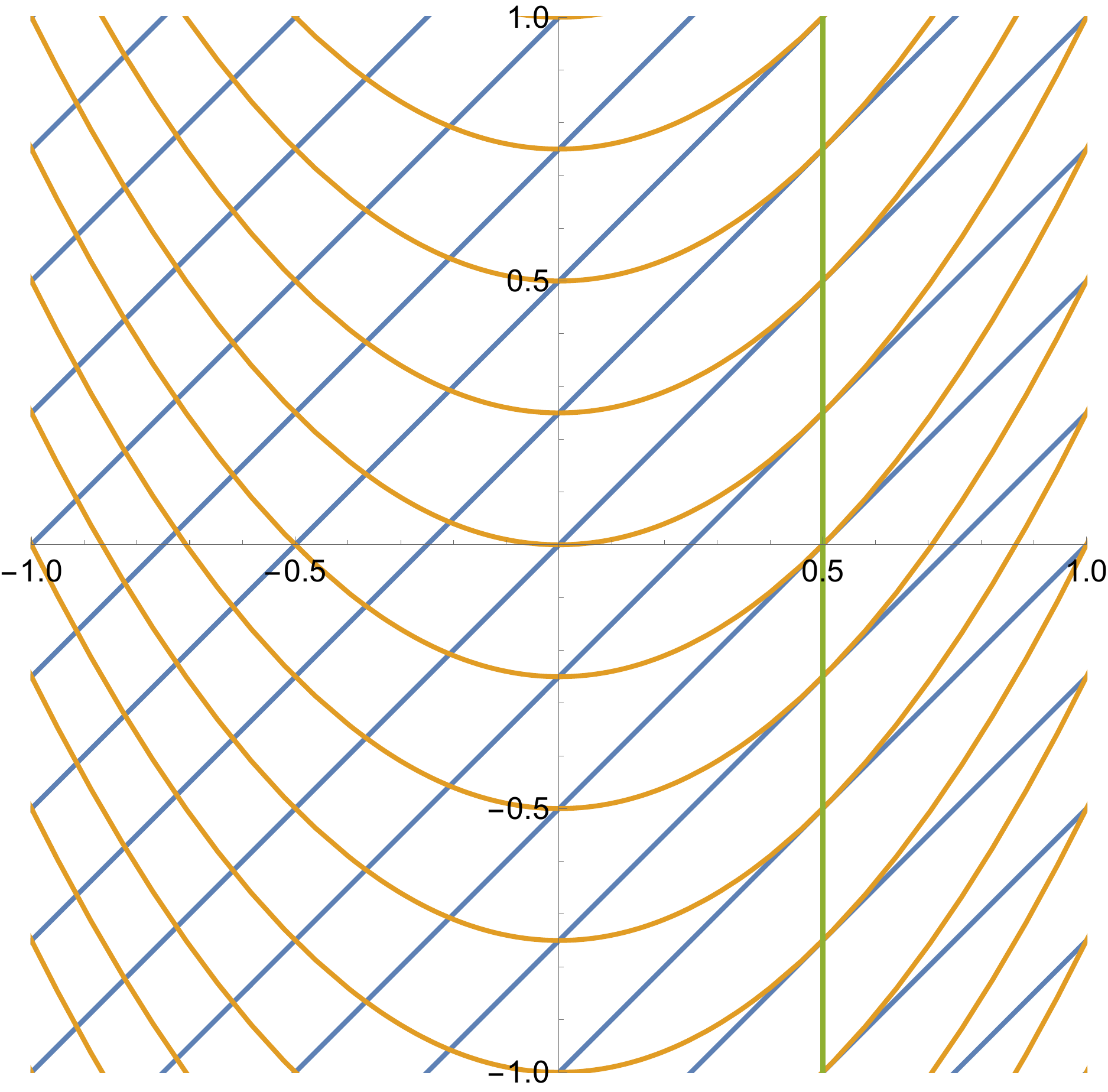}
\end{center}
In the left figure, {\em no matter what} blue line we choose, the conditional expectation of $X$ is $0$.  In the second and third figures, {\em no matter what} blue line we choose, the conditional expectation of $X$ is (respectively) $-1/2$ or $1/2$. Changing the blue-line slope somehow has the effect of ``shifting'' the conditional law of $X$.

This is a bit counter-intuitive. For a closely related example (essentially a rotation of the one above by 45 degrees) suppose $$e^{-(a-b)^2/2t} e^{(a+b)/2}dadb$$ is the (semi-Gaussian) density function for a pair $(A,B)$.  Then we can formally compute $\mathbb E[B|A] = A+t/2$ and $E[A|B] = B+t/2$.  This is because if we restrict to a fixed value of $a$, we obtain a multiple of a Gaussian measure on $b$ values centered at $a+t/2$ (and similarly with $a$ and $b$ reversed). Checking this fact is an elementary complete-the-square calculation: if $b$ is fixed then $-(a^2-2ab+b^2)/2t + (a+b)/2$ is $-\bigl( a^2-2a(b+t/2) + (b+t/2)^2 \bigr)/2t$ plus a term that doesn't depend on $a$.

At a glance, the above seems to suggest that $A$ is $t/2$ units bigger than $B$ on average, and $B$ is $t/2$ units bigger than $A$ on average, which in turn seems contradictory. This is somehow reminiscent of ``envelope switching'' paradoxes, where after one observes the amount of money in either {\em one} of two envelopes, one always expects the {\em other} to contain more, see e.g.\ \cite{broome1995two}. In fact this is just the sort of paradox that one encounters when dealing with infinite measures and/or infinite expectations. We find that the following are equivalent:
\begin{itemize}
\item First ``sample'' $A$ from the infinite measure $e^a da$.  Then choose $B$ as a normal with variance $t$ and mean $A+t/2$.
\item First ``sample'' $B$ from the infinite measure $e^b db$. Then choose $A$ as a normal with variance $t$ and mean $B+t/2$.
\end{itemize}
One can imagine similar constructions in higher dimension.  For example, we could replace the $X$ in the figures above by a vector $(X_1, X_2, \ldots, X_n)$ and replace $y-x^2$ with $y-\sum_{i=1}^n x_i^2$, and replace the blue lines by blue hyperplanes.  One could then argue in a similar way that changing the slope of the blue hyperplanes has the effect of shifting the location of the vertical green line.

By taking a limit of an increasing sequence of finite-dimensional subpaces, we can also make sense of an infinite dimensional semi-Gaussian---precisely the same way we make sense of the GFF as an infinite dimensional Gaussian, using the quadratic function $-\frac12 (\phi, \phi)_{\nabla}$, or the way we define Brownian motion, using the quadratic function $\frac12 \int_{-\infty}^\infty \bigl( \frac{\partial}{\partial t} B(t)\bigr)^2dt$.

Let us give an example.  Define a process (whose law is an infinite measure) by $M(t) = B(t) + |t|/2 + Y$ where $Y$ has density $e^y dy$ and $B(t)$ is an independent standard Brownian motion (defined for all $t \in \mathbb R$, with $B(0) = 0$).  This is an infinite measure on paths that at first glance seems to be ``mostly supported'' on paths that have their minima near $0$.  If we tried to define an ``action'' for $M(t)$ it would have the form $$S(M) = e^{M(0)} + \frac12 \int_{-\infty}^\infty \Bigl( \frac{\partial}{\partial t} \bigl( M(t)- |t|/2 \bigr)^2\Bigr)dt.$$
On the other hand, some thought reveals that, when $t$ is fixed, the (infinite-measure) law of $M(0)$ and $M(t)$ is equivalent to the law of $A$ and $B$ in the example above, and in fact the law is symmetric w.r.t.\ swapping the roles $0$ and $t$.  One can also show that {\em given} $M(0)$ and $M(t)$ the conditional law of the rest of the process is given by a Brownian bridge on $(0,t)$ and Brownian motions with drift on $(t, \infty)$ and (running time backward) on $(-\infty, 0)$, with the given values of $M(0)$ and $M(t)$ as endpoints. One can make a similar argument using $-t$ and $0$ in place of $0$ and $t$ and use this to show that the the law of $M$ is invariant w.r.t.\ to translation by $t$ units (for any $t$).

This is a rather remarkable fact.  The definition of $M(t)$ does not look translation invariant at all---it is clearly centered at $0$.  But somehow translating the definition by $r$ units to the right does two things: it effectively {\em weights} the law of $M$ by $e^{M(r)-M(0)} = e^{B(r) - B(0)}$ (this puts a bias in favor of functions that increase on the interval $[0,r]$) and it deterministically adds the function $f(t) = \frac{|t-r|}{2} - \frac{|t|}{2}$ (which decreases on the interval $[0,r]$ and is flat elsewhere).  And these two changes magically cancel each other out.

We remark that another way to construct this measure is by considering a Brownian bridge measure on processes $\mathcal B(t)$ indexed by $[-T,T]$, with boundary values $\mathcal B(T)=\mathcal B(-T)=T/2$, then multiplying the measure by a constant to ensure that the measure of paths with $\mathcal B(0)\in[-1,1]$, say, is some fixed constant, then taking the vague $T \to \infty$ limit. This construction might make the translation invariance a bit more intuitive.

To extend the story to two dimensions, note that the above construction is still translation invariant if we take the law of $Y$ to be $e^{-2Qy}dy$ and write $M(t) = B(t) - Q|t| + Y$. This just corresponds to taking $\mathcal B(T) = \mathcal B(-T) = -QT$ in the limiting procedure mentioned above. On the cylinder $\mathcal C = \{t+\theta i: t \in \mathbb R, i \in [0,2\pi) \}$ one can then define $\phi(t+\theta i)$ to be $M(t)$ plus the projection of the GFF (on the cylinder) onto the space of functions with mean zero on every fixed-$t$ slice of the cylinder. If restrict to $M(t)$ that have their maximum at $t=0$ (to eliminate the translation symmetry) we obtain the so-called ``$\alpha=0$ quantum sphere'' \cite{matingoftrees}.  If we {\em don't} restrict in this way, we obtain the Polyakov measure (with $
\mu=0$), an infinite measure on surfaces that turns out to be invariant w.r.t.\ {\em all} M\"{o}bius coordinate changes, not only translations.

To see why, note that there are three ``mostly flat'' surfaces often used to parameterize an LQG sphere: the Riemann sphere $\mathbb C \cup \{0\}$, the cylinder $\mathcal C \cup \{ \pm \infty \}$ and the gluing $\mathcal D$ of two unit disks together their boundaries.  There are obvious unit-circle-preserving conformal maps between these three spaces. If $\psi$ is the obvious conformal map from $\mathcal C$ to $\mathcal D$ then $Q \log |\psi'(t+i\theta)| = -Q|t|$. This is the term that comes up in the change of coordinates from $\mathcal D$ to $\mathcal C$. If we instead change coordinates from $\mathcal D$ to $\mathbb C$, then this term is replaced by $-2Q \max \{\log|z|,0 \}$. Working in $\mathbb C$, the proof that the law of the construction is invariant under translations and rotations of $\mathbb C$ (as well as inversions and dilations) is very similar to the one dimensional argument above, and since these operations generate the M\"{o}bius group, full M\"{o}bius invariance follows, see \cite{ang2021integrability2, kupiainen2016constructive}.


\subsubsection{Embedded Polyakov sphere} \label{sec::embeddedpolyakov}
When building physical field theories it is often natural to define a measure on fields $\phi$ by writing formally $e^{-S(\phi)}d\phi$ where $d\phi$ represents some sort of uniform measure on the space of all functions. In some cases this is hard to make mathematically precise because it is not clear what the measure $d\phi$ corresponds to. One might attempt to approximate $d\phi$ by something well defined (like Lebesgue measure on a finite-dimensional space of piecewise-linear functions) but it may be unclear how to normalize the construction to obtain a non-trivial limit, or how to prove such a limit exists. In the special case of the Liouville action, however, these concerns can be overcome. Let us tell the story with a playful dialog.

\vspace{.03in}
\noindent INSTRUCTOR: Consider the measure $e^{-A(\phi) - B(\phi)} d\phi$ where $d\phi$ is the uniform measure on the space of all functions.
\vspace{.03in}

\noindent MATH POLICE: Sorry, there is no such thing as the uniform measure on all functions. Your object is not defined.
\vspace{.03in}

\noindent INSTRUCTOR: But my $A$ is a quadratic function on the space of $\phi$ for which it is finite. If I restrict $A$ to a co-dimension-one subspace then $A$ is the norm for a Hilbert space. So $e^{-A(\phi)}d\phi$ is just a Gaussian Hilbert space cross an infinite measure that looks something like $e^{y}dy$. Check out Janson's book on Gaussian Hilbert spaces \cite{janson1997gaussian}. These things are certainly well defined.
\vspace{.03in}

\noindent MATH POLICE: Okay fine, but you still need to weight by the $e^{-B(\phi)}$ factor.
\vspace{.03in}

\noindent INSTRUCTOR: My friends proved that (with appropriate normalizing) $B(\phi)$ is well-defined and finite for almost all $\phi$ taken from the measure $e^{-A(\phi)}d\phi$.
\vspace{.03in}

\noindent MATH POLICE: So $e^{-B(\phi)}$ is the Radon-Nikodym derivative w.r.t.\ the semi-Gaussian measure? And $B$ is obviously measurable?
\vspace{.03in}

\noindent INSTRUCTOR: That's right.
\vspace{.03in}

\noindent MATH POLICE: Okay, you're free to go.
\vspace{.03in}

We can take the $S=A+B$ in the above dialog to be the {\em Liouville action} defined as follows:
$$S(\phi) = \int \Bigl( \frac{1}{2} \cdot \frac{1}{2\pi}  |\nabla \phi(z)|^2 +  \frac{1}{4\pi} QR(z) \phi(z) +  \mu  e^{\gamma \phi(z)}\Bigr)dz,$$
where $R$ is the curvature associated to the ``reference metric,'' see e.g.\ \cite{david2016liouville}. For example, we may assume that the reference metric is the ordinary sphere so that $R$ is constant (with total integral $8 \pi$ by Gauss-Bonnet---recall that the Ricci curvature is twice the Gaussian curvature) so that the second term is $2Q$ times the mean value of $\phi$. Alternatively, if the reference metric is the glued pair of disks $\mathcal D$ from the previous subsection, then the second term becomes $2Q$ times the mean value of $\phi$ on the unit circle.

The integral of the first two terms is the quadratic part $A(\phi)$ while the integral of the last term is the $B(\phi)$. The fact that the latter is defined (provided $Q > 2$) follows essentially from the ideas of H{\o}egh-Krohn and Kahane---see also the author's later work with Duplantier \cite{lqgkpz} which constructs $e^{\gamma \phi(z)}dz$ as a measure-valued function of the field $\phi$ and explains the LQG context. In light of the above dialog, the Polyakov measure $e^{-S(\phi)}d\phi$ is in fact rigorously defined; this is explained in more detail in \cite{david2016liouville} (see also the higher genus version in \cite{david2016liouvilletori}).

Stories like the above, where $S$ has a quadratic term (relatively simple to handle) and a non-quadratic term (requiring more thought) are relatively common in quantum field theory, and it is certainly not {\em always} the case that the non-quadratic part simply modifies the Gaussian part in an absolutely continuous way.

We stress that the $\mu=0$ measure is already rather interesting. This is the measure shown to be M\"{o}bius invariant (following the LQG coordinate change rules), despite appearing at first glance to be ``centered'' at a specific location within the M\"{o}bius group.  The miracle is that when one applies a M\"{o}bius transformation $\psi$, the $Q \log |\psi'|$ factor one has to add exactly compensates for the effect of the recentering (swapping the curvature measure $R(z)dz$ for its image under $\psi$), see \cite{ang2021integrability2}. This implies that the Polyakov measure factors as the product of a measure on (unembedded) LQG surfaces and a measure on embeddings (given by the infinite-volume Haar measure on the M\"{o}bius group) see \cite{ang2021integrability2}. The same is true if one considers a non-zero $\mu$. However, the Polyakov measure with a non-zero $\mu$ only makes sense as described above if $d \leq 1$ (since otherwise the natural volume measure is infinite), while the construction of the zero-$\mu$ measure makes sense for any $Q \geq 0$.   Changing the choice of reference metric is analogous to changing the slope of the blue lines/hyperplanes in the previous section and (up to a constant multiplicative factor) it leaves the measure construction unchanged \cite{david2016liouville}.

On a flat reference metric, the quadratic action is simply the Dirichlet energy of $\phi$ while the $\mu e^{\gamma \phi(z)}$ term corresponds to the $K e^{\rho}$ term associated to the Liouville equation with a non-zero $K$. Polyakov was not the first to work with variants of the free field action, and in \cite{polyakov2005confinement} and \cite{polyakov2008quarks} he attributes a closely related ``quadratic action'' to Douglas, who used it in his work on Plateau's problem for minimal surfaces. (Douglas was awarded one of the first two Fields medals for this work in 1936 \cite{douglas1931solution}.) See also \cite{brink1976locally, deser1976complete} which introduce the action now called the {\em Polyakov action}, which Polyakov used later in \cite{polyakov1981quantum}.

\subsection{More history}
The Brownian map and peanosphere constructions were directly motivated by discrete objects---namely random planar map models.  The construction of LQG surfaces, on the other hand, was motivated within physics as a ``quantization'' based on the Dirichlet energy or on the Liouville equation. Why would one expect these objects to be equivalent?

The discrete planar map models were well studied in physics due to their relationship with random matrices and the random particle systems corresponding to random matrix eigenvalues, see the seminal papers from the 1970's by t'Hooft and by Br\'ezin, Itzykson, Parisi and Zuber \cite{thooft1974planar, brezin1978planar}, along with more recent overviews by Eynard, Guionnet and Maurel-Segala \cite{eynard2016counting, guionnet2005combinatorial}. However, it took some time for people to be persuaded that the Liouville theory corresponded to the scaling limit of these models. For example, Polyakov wrote in a memoir \cite{polyakov2008quarks} (see also his previous memoir \cite{polyakov2005confinement}) that he did not become convinced of the connection between the discrete models and Liouville quantum gravity until the late 1980s after jointly deriving, with Knizhnik and Zamolodchikov, the so-called {\em KPZ formula} for certain Liouville quantum gravity scaling dimensions \cite{knizhnik1988fractal} and comparing them with known combinatorial results for the discrete models (for a rigorous approach see also \cite{lqgkpz, benjamini2009kpz, rhodes2011kpz, garban2013quantum, duplantier2014critical}).

At this point the relationship between the planar map models and the continuum objects is on more solid mathematical ground. We have convergence results of various kinds and formal equivalence proofs between the different continuum objects, see Section~\ref {sec::relationships}. But these rely on a lot of machinery that was not yet available in the 1980s. The quantum gravity zipper perhaps gives the cleanest way to rigorously relate the peanosphere construction and the Liouville quantum gravity construction, as described in \cite{zipper, matingoftrees}, see also \cite{astala2011random}.

\subsection{Computing the scaling exponent} \label{sec::peanospherescaling}
On the other hand, just as in the peanosphere case, one can compute the scaling exponent quite easily. This is already enough to show that {\em if} for some $\gamma$ the unrestricted-area LQG sphere is the scaling limit of the unrestricted-area discrete models, then we must have $\gamma = \sqrt{8/3}$ in the undecorated case and $\gamma = \sqrt{2}$ in the tree-decorated case.

Let us sketch that calculation here. Say we fix the Gaussian free field up to the mean value $X$ on the unit circle. Then it is not hard to see (recall Section~\ref{sec::embeddedpolyakov}) that this mean value has a conditional law given by the infinite measure $e^{-2Q x}dx$. The total area of $e^{\gamma \phi(z)}dz$ will then be some constant (depending on the given value of $\phi$ modulo constant) times $e^{\gamma X}$. To understand the conditional law of this quantity, we then have to consider the image of the measure $e^{-2Q x}dx$ under the map $\psi(x) =  e^{\gamma x}$.  By a standard change of variables the new measure is (up to a constant factor) $$\frac{1}{(\psi^{-1})'(A)} e^{-2Q \psi^{-1}(A)} dA = \frac{1}{A} e^{-2Q \log (A)/\gamma} = A^{-1 - 2Q/\gamma} dA.$$
Here $1+2Q/\gamma = (4/\gamma^2 + 2) $ which comes to $7/2$ if $\gamma = \sqrt{8/3}$ and $4$ if $\gamma = \sqrt{2}$.

Since this calculation works for any fixed choice of $\phi$ modulo additive constant in the Polyakov measure, it is not hard to deduce that the same scaling rate must hold for the unrestricted-area LQG sphere overall.


\subsection{Random surfaces embedded in $d$-dimensional space}

Suppose we accept, based on the previous discussion, that the $\gamma = \sqrt{8/3}$ theory (which corresponds to $25-6Q^2 = 0$) describes the scaling limit of the undecorated $d=0$ model. We can call this {\em pure Liouville quantum gravity} model. We would then like to argue that if we {\em weight} the law of the pure model by the $d$th power of the GFF partition function (or the corresponding and equivalent ``loop soup'' partition function) we obtain a new LQG model with a $Q$ parameter that satisfies $d=25-6Q^2$.

This has been explained heuristically in various ways over the years. The author with Ang, Park and Pfeffer gave a rigorous version of this statement that applied to a certain way of regularizing the random surfaces \cite{ang2020brownian}. We won't give details here but we mention below a few of the ingredients used to make the connection between loop soup weightings and changes to $d=25-6Q^2$.

On a compact surface with boundary, the heat kernel trace can be written $Z = Z(t) = \textrm{sp}\, e^{t\Delta} = \sum e^{t \lambda_n}$ where $\lambda_n$ are the eigenvalues of the Laplace-Beltrami operator $\Delta$.  By standard Tauberian theory, the asymptotics of $Z$ (as $t \to 0$) are closely  related to the asymptotics of $\lambda_n$ (as $n \to \infty$). Weyl addressed the latter for bounded planar domains $D$ in 1911 \cite{weyl1911asymptotische} (see discussion in \cite{mckean-singer}) by showing $-\lambda_n \sim \frac{2\pi n}{\textrm{area}(D)}$ as $n \to \infty$ which is equivalent to \begin{equation*} Z \sim \frac{\textrm{area}(D)}{4 \pi t}\end{equation*} as $t \to 0$. In 1966 Kac gave higher order correction terms for $Z$ on domains with piecewise linear boundaries (accounting for boundary length and corners) in his famously titled ``Can you hear the shape of a drum?''\ which asks what features of the geometry of $D$ can be deduced from the $\lambda_n$ or equivalently from $Z$ \cite{kac1966can}.  (Short answer: some but not all.)

McKean and Singer (among others) extended these asymptotics from planar domains to smooth manifolds with non-zero curvature \cite{mckean-singer} where the constant order correction term is a certain curvature integral. For two dimensional surfaces with boundary, the integral $\int_\delta^\infty Z(t)/t dt$ turns out to describe the {\em Brownian loop soup measure} of the set of loops longer than $\delta$ (as developed and explored by Lawler, Werner, Dub\'edat and others -- see \cite{ang2020brownian} for further explanation and references).  When the metric takes the form $e^{\rho}$ times a flat metric, the small $\delta$ asymptotics have a constant order correction term that corresponds to the Dirichlet energy of $\rho$. This is the so-called Polyakov-Alvarez formula, also known as the Polyakov-Ray-Singer or Weyl anomaly formula, and it has played a major role in conformal field theory. See e.g.\ the early discussion of Laplacian determinants and this formula by Osgood, Phillips and Sarnak \cite{sarnak1990determinants, osgood1988extremals}.

\section{Conformal field theory and multipoint correlations} 
\label{sec::cft}

Conformal field theory is a huge subject.  For a broader overview of conformal field theory, beyond just Liouville theory, the reader might begin with the well-known (and very long) textbook by Di Francesco,  Mathieu, and S\'en\'echal \cite{francesco2012conformal}. Alternatively, the overview at \cite{ribault2014conformal} begins with a list of several conformal field theory textbooks.

Liouville theory is one of many conformal field theories, but it is by itself a large and highly studied subject. Here one might begin with the 1990 survey by Seiberg \cite{seiberg1990notes}, the 1993 lecture notes by Ginsparg and Moore \cite{ginsparg1993lectures}, the 1995 survey by Di Francesco, Ginsparg and Zinn-Justin \cite{di19952d}, the 2004 retrospective by Nakayama \cite{nakayama2004liouville} or the 2014 textbook by Ribault \cite{ribault2014conformal} (which among other things explains why Liouville conformal field theory is uniquely characterized by certain axioms).

What is a conformal field theory? An internet search for ``A conformal field theory is'' reveals several definitions, the first two from Wikipedia:
\begin{itemize}[leftmargin=.16in]
\item a quantum field theory that is invariant under conformal transformations.
\item a set of correlation functions that obey a number of axioms.
\item a functor \cite{moore1989classical} between categories satisfying certain ``sewing axioms.''
\item a Virasoro module \cite{schottenloher2008mathematical}

$$V = \oplus_{i=B_1} W(c_i,\phi_i) \otimes W(\bar c_i, \bar \phi_i)$$
with unitary highest weight modules $W(c_i, \phi_i)$, $W(\bar c_i, \bar \phi_i)$ subject to [certain] axioms.
\end{itemize}
The first definition is standard but its precise meaning depends on how one defines a quantum field theory. The third and fourth definitions represent formalization efforts that would go somewhat beyond the scope of this note. So let us focus on the second definition, sometimes called the {\em conformal bootstrap} approach to CFT.  Here a CFT is no more or less than a set of {\em correlation functions}, and the interpretation of these functions is that they represent (in some sense) expected products of random generalized functions called {\em fields} evaluated at different points.

The famous 1984 paper by Belavin, Polyakov and Zamolodchikov (BPZ) \cite{belavin1984infinite} argued that conformal invariance symmetries should imply certain properties for these correlation functions, and that this should be sufficient to allow one to explicitly compute the correlation functions for {\em some} conformal field theories (the so called {\em minimal models}) including a theory that one would expect to describe the scaling limit of the Ising model.  
About a decade later Dorn, Otto, Zamolodchikov and Zamolodchikov were able to compute certain three-point correlations for the Liouville theory \cite{dorn1994two, zamolodchikov1996conformal}. Other correlation functions could then be deduced from these using the BPZ theory \cite{belavin1984infinite} and further input proposed by
Teschner \cite{teschner1995liouville, teschner2001liouville, teschner2004lecture}.
These arguments involved mathematically non-rigorous steps, such as assuming without proof that formulas defined in one setting could be analytically continued and applied in other settings.

Building on \cite{david2016liouville}, Guillarmou, Kupianen, Rhodes, and Vargas have produced a series of papers that define and derive the correlation functions for Liouville theory mathematically, building on earlier derivations from the physics literature that we mentioned above. The expressions describing the correlation functions are complicated (integrals, special functions, recursive definitions, etc.) but nonetheless explicit. The first work in this series is a proof of the DOZZ formula \cite{kupiainen2020integrability}. The second paper derives an analog of Plancharel's theorem (which states that the Fourier transform preserves $L^2$ norm) along with a certain ``spectrum'' relevant to this context \cite{guillarmou2020conformal} . The final paper completes the bootstrap program with an extension to $n$-point functions and higher genus surfaces \cite{guillarmou2021segals}. We recommend that the reader take a look at the introduction to \cite{guillarmou2020conformal}, which summarizes this viewpoint and situates it within the larger enterprise of quantum field theory. Here we will give a much shorter overview of this viewpoint, aiming only to give a simple account of the relationship to the other perspectives in this paper.

The physics literature on conformal field theory can be challenging for mathematicians to follow. It comes with a large and very specialized jargon, and it does not always proceed in the order a mathematician would expect (where one first produces the measure space and the $\sigma$-algebra, then constructs the measure, then begins doing calculations). Partly this is because (when quantum wave functions are involved) not everything in quantum field theory {\em can} be described in a simple probabilistic way---sometimes ``observables'' are non-commuting operators that can only be defined in indirect ways. Fortunately, Liouville theory does have a simple probabilistic interpretation.

Let us make one more comment on nomenclature. In quantum field theory a (generalized) function on $\mathbb R^{n+1}$ can be interpreted as a {\em path} on the space of functions on $\mathbb R^n$, and the term {\em integral} is often a shorthand for {\em a measure w.r.t.\ which one integrates}. In this context, the object we call the Polyakov measure (on the space of surfaces) in Section \ref{sec::polyakov} is also called the {\em Polyakov path integral} (with $n$ taken to be $1$). This language evokes the Feynman path integral, an integral over particle trajectories that appears in quantum mechanics. Similar interpretations arise in string theory, where an integral over the space of string trajectories is seen as an integral over the space of embedded surfaces (a.k.a.\ ``worldsheets''). A very informal overview of this viewpoint (featuring Vargas) appears in an online video produced by Quantum Magazine (starting at time 7:20) \url{https://youtu.be/9uASADiYe_8}.

\subsection{Gaussian case} \label{subsec::gaussiancase}

It is generally instructive to do something easy before doing something hard. So let us start in the simple setting where $\phi$ is a whole plane Gaussian field with additive constant chosen so that the mean value on the unit circle is $0$. (This is also the starting point in \cite{guillarmou2020conformal} for example.)  In this case the Green's function is given by
$$G(x,y) =\ln \frac{1}{|x-y|} + \ln|x|_+ + \ln |y|_+$$
where $|x|_+ = \max\{1,|x| \}$ and the construction of ``random fields'' of the form $V_\alpha(x) := e^{\alpha \phi (x)}$ is relatively straightforward. In some sense it can already be seen in the work of quantum field theorists of the 1970's, beginning with the work of H{\o}egh-Krohn. As we have already discussed, one way to make sense of this is by writing
 $V^\epsilon_\alpha(x) := \epsilon^{\alpha^2/2} e^{\alpha \phi_\epsilon(z)}$ and $V_\alpha(x) := \lim_{\epsilon \to 0} V^\epsilon_\alpha(x)$, where $\phi_
 \epsilon(z)$ is the mean value of $\phi$ on $\partial B_\epsilon(z)$ and the convergence holds locally a.s\ in the space of random generalized functions (or in the space of random measures). Then for any sufficiently small $\epsilon$ (so that the balls $B_\epsilon(x_i)$ don't overlap) we have
$$\Big \langle \prod V^\epsilon_{\alpha_i}(x_i) \Big\rangle = e^{\sum \alpha_i \alpha_j \tilde G(x_i, x_j)}$$ where $$\tilde G(x,y)= \begin{cases} G(x,y) & x\not = y \\ \log |x|_+ + \log |y|_+ & x = y \end{cases}.$$ 
This comes up out to be $$\prod_{i \not = j} |x_i-x_j|^{\alpha_i \alpha_j} \prod_{i,j} \Bigl( \max \{ |x_i|, 1 \} \max \{|x_j|,1 \} \Bigr)^{\alpha_i \alpha_j}$$
and in the case that all the $x_i$ lie inside the unit disc, the expression is simply
$$\prod_{i \not = j} |x_i - x_j|^{\alpha_i \alpha_j}.$$
By starting with $V_{\alpha_i}^\epsilon$ and taking $\epsilon \to 0$ we can give meaning to the multipoint correlation function (a.k.a.\ {\em Schwinger function}) written
$$\Big \langle \prod V_{\alpha_i}(x_i) \Big\rangle = \prod_{i \not = j}  |x_i - x_j|^{\alpha_i \alpha_j}.$$
This can be interpreted as the density function for a two-dimensional ``Coulomb gas'' of particles with charges $\alpha_i$.  The overall integral over all $x = (x_1,\ldots, x_n)$ may be infinite (depending on the $\alpha_i$ values) but it is finite if one restricts to a subset of $(x_1, \ldots, x_n)$ values such that the $|x_i|$ are bounded above and the $|x_i-x_j|$ are bounded below. Just to clarify: for now we are making this Coulomb gas calculation only for the centered GFF $\phi$, not for the full Polyakov measure from Section~\ref{sec::polyakov}.

As a formal collection of correlation functions (obtained in the $\epsilon \to 0$ limit) this expression makes perfect sense for any $\alpha_i \in \mathbb C$.  Differentiation commutes with expectation, so one can use the same formula to compute expectations of products involving derivatives like $\frac{\partial}{\partial x_i} V_{\alpha_i}(x_i)$. Expectations of fields involving $\phi$ itself (or polynomials in $\phi$) can also be defined mathematically.

\subsection{Incorporating the Liouville term or area conditioning}

The $\phi$ used in Section~\ref{subsec::gaussiancase} (plus a deterministic function) can be obtained by restricting the Polyakov measure to the set of $\phi$ whose mean value on the unit circle is zero. But what if we instead restrict the Polyakov measure to the set of surfaces of total quantum area $1$ (or weight by $e^{-\mu A}$ where $A$ is the surface area)?  The answer is that after such a conditioning or such a weighting, $\phi$ is no longer Gaussian, and the $n$-point correlation computation transforms from {\em easy} to {\em doable but only barely}.  On the other hand, in order to understand some fundamental things like the law of conformal modulus of four points (or $n$ points) sampled independently from the measure on an LQG sphere, one has to address the harder question, and this is precisely what is done in \cite{kupiainen2020integrability, guillarmou2020conformal, guillarmou2021segals}. Although these results are rather recent, they have already inspired a tremendous amount of activity, establishing {\em exact solvability} for many problems that could previously only be addressed more qualitatively.

Papers by Ang and Sun, including some with co-authors Holden and Remy, give applications of LCFT results to various domains (SLE, CLE, and a variance formula for the peanosphere) \cite{ang2020conformal, ang2021integrability, ang2021fzz, ang2021integrability2}. These remarkable papers combine mating-of-trees and conformal-welding techniques with LCFT techniques, leading to rigorous proofs of physics results (such as the FZZ formula and the imaginary DOZZ formula) as well as entirely new results. They have to extend the welding/mating theory to finite-volume surfaces/trees, which is harder than the infinite-volume work, since the finite-volume surfaces lack scale invariance.  Recent integrability achievements for the disk include works by Remy and Zhu \cite{remy2020integrability, remy2020distribution} and a proof of the Fyodorov-Bouchaud formula by Remy \cite{remy2020fyodorov}. See also the work of Ghosal, Remy, Sun and Sun on the torus \cite{ ghosal2020probabilistic}.

\section{Relationships} \label{sec::relationships}
Although all four viewpoints are equivalent in some sense, the relationships are somewhat involved.  The following diagram summarizes a few of the associated keywords, which we will discuss briefly below.

\begin{center}
\includegraphics[width=\textwidth]{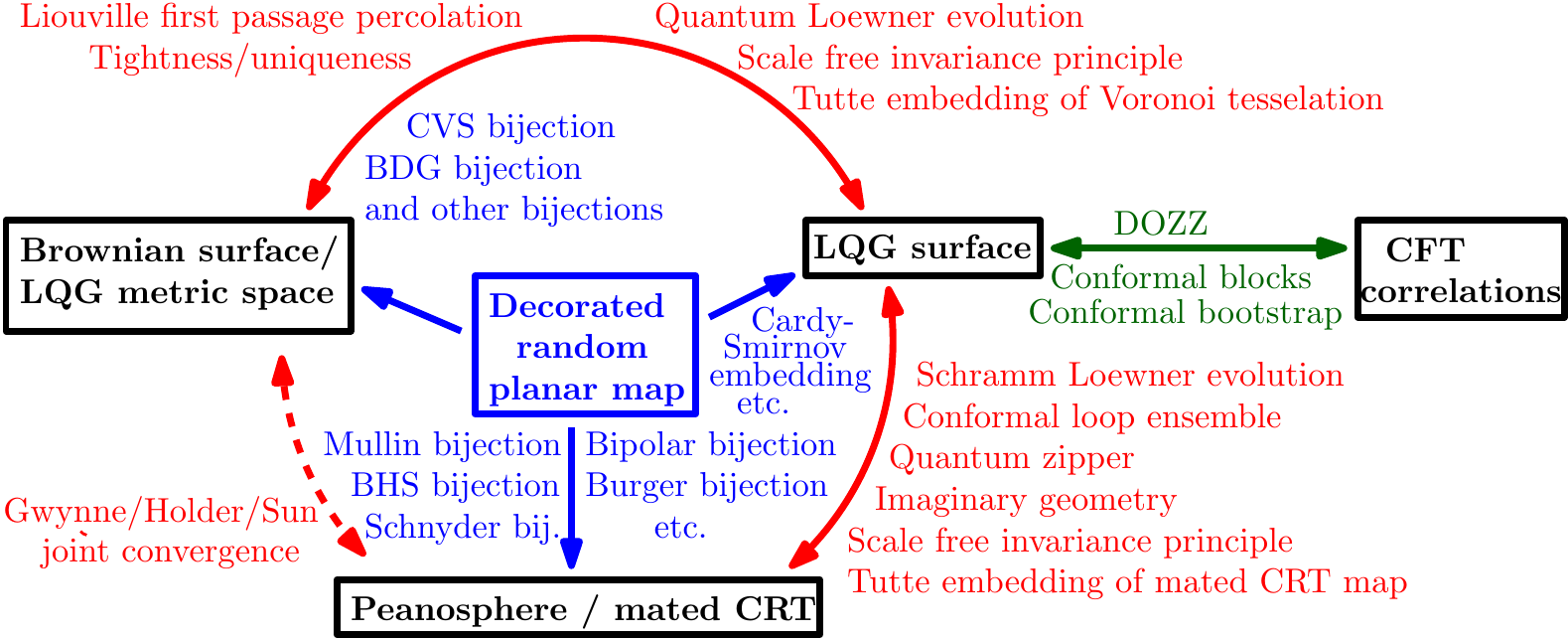}
\end{center}

Miller and the current author recently showed that the LQG sphere and the Brownian sphere are {\em a posteriori} equivalent \cite{lqgtbm1, lqgtbm2, lqgtbm3}. That is, we showed that there is a canonical way to endow each object with the other object's structure, and that once this is done the two objects agree in law. This is not all obvious. The papers are long and difficult and build on hundreds of pages of prior work. The author's other joint works on the LQG construction include \cite{lqgkpz, lqggff, fieldmeasurecorrespondencecommutes, criticalchaos2, duplantier2014critical}.  This establishes the upper arrow in the special case of the Brownian map.  The construction of the metric space structure for general LQG surfaces (with $d \not = 0$) was discussed earlier.

In \cite{zipper} the author showed that when two infinite half-plane-homeomorphic $\gamma$-LQG surfaces are ``conformally welded'' to one another along their boundaries --- and the new surface is conformally mapped to the half-plane --- the interface between these curves becomes an SLE$_\kappa$ curve with $\kappa = \gamma^2$. The sphere version is established in a follow-up work \cite{finitetree}, and a disk version appears in \cite{ang2021liouville}. Establishing the correlation formula for $\kappa > 8$ was achieved in \cite{gwynne2017brownian}.

With Duplantier and Miller, the author showed the equivalence of the peanosphere approach and the SLE-decorated-LQG approach in \cite{matingoftrees}, which draws from the imaginary geometry results in \cite{IG1, IG2, IG3, IG4} and the quantum zipper construction in \cite{zipper}. A remarkable series of papers by Holden and Sun has shown that if one embeds the uniformly random triangulation in the plane in a natural way (inspired by the conformal coordinates Smirnov developed for his proof of Cardy's formula) then the counting measure on the vertices converges (in an appropriate scaling limit) to the Liouville quantum gravity measure, see the overview at \cite{holden2022natural, holden2019convergence}.

Random planar maps converge to continuum objects in the limit, but there are also natural ways to generate random planar maps from the continuum constructions by ``coarse graining'' in some sense. The Poisson-Voronoi tesselation (constructed from the Brownian map) and mated-CRT map (constructed from the peanosphere) have been shown to converge to LQG when they are embedded according to the Tutte embedding.  This was done in a series of four papers with Gwynne and Miller \cite{harmonicmatedcrt, tuttebrownian,  gwynne2021invariance, tuttelqg}. Effectively, this gives a way of putting a conformal structure on the Brownian map or the matings of trees that is more concrete than the one guaranteed by \cite{lqgtbm1, lqgtbm2, lqgtbm3}.  The heart of all of this is \cite{gwynne2021invariance} which gives an invariance principle (i.e., Brownian motion convergence) for random walks in random environments that are ``scale free'' in the sense that there is no universally typical length scale. (All of the most natural discretizations of Liouville quantum gravity measures are scale free in this sense.)

Gwynne, Holden and Sun established the {\em joint convergence} of random triangulations in the metric and peanosphere sense \cite{gwynne2021joint}. (Convergence in one topology coupled with convergence in another topology does not imply joint convergence in the product of the two topologies, but it was established in this particular case.) This result builds on earlier convergence work by Gwynne and Miller for percolation-interface-decorated surfaces \cite{gwynne2021convergence, gwynne2021characterizations}.

A series by the author plus Miller and Werner explores conformal loop ensembles on random surfaces, the surfaces obtained by cutting along the boundaries of these loops, and so forth \cite{miller2020simple, miller2021non}. These papers extend the conformal welding stories described earlier and show that much of the intuition one derives from discrete planar maps is correct in the continuum as well (e.g., the surfaces inside and outside of a CLE loop are independent given the interface length) and this leads to a number of interesting computations.  These papers build on another recent paper by the same authors \cite{clepercolations} which concerns continuum versions of the classical Edwards-Sokal couplings involving FK-models and Ising/Potts models and their variants, see also \cite{notbyrange}.

\section{Gauge theory}
Random surfaces are related to many areas of math and physics, including random matrix theory, two-dimensional statistical physics, string theory, and so on.  But we should also note that Polyakov's influential 1981 paper began by mentioning an interest in gauge theory:
\begin{quote}“In my opinion at the present time we have to develop an
art of handling sums over random surfaces. These sums replace the old-fashioned
(and extremely useful) sums over random paths. The replacement is necessary,
because today gauge invariance plays the central role in physics. Elementary
excitations in gauge theories are formed by the flux lines (closed in the absence
of charges) and the time development of these lines forms the world surfaces. All
transition amplitude[s] are given by the sums over all possible surfaces with fixed
boundary.” (A.M. Polyakov, Moscow, 1981.) \cite{polyakov1981quantum}
\end{quote}
Over 40 years later, many fundamental gauge theory problems remain unsolved, including the famous Clay Millenial Prize Problem, and it is unclear how much random surface theory will help---see e.g.\ the ``skeptic vs.\ enthusiast'' dialog in Section 2 of \cite{cordes1995lectures}, written in 1995. Nonetheless, there has been some progress on so-called {\em gauge string duality}. For more on these efforts, the reader may consult the literature on the AGT conjecture and the AdS/CFT correspondence, or see the recent works of Chatterjee and Jafarov on lattice string trajectories and Yang-Mills theory \cite{chatterjee2019rigorous, chatterjee2016leading, chatterjee20161, jafarov2016wilson}. Good entry points to the subject for probabilists include Chatterjee's recent survey \cite{chatterjee2016yang} and Thierry L\'evy's recent books about continuum Yang-Mills theory in two dimensions \cite{levy2003yang, levy2017master, levy2020two}. The 1986 article by Bridges, Giffen, Durhuus and Fr{\"o}hlich may also be read as a first step toward realizing Polyakov's vision \cite{brydges1986surface}.

\begin{ack}
Random surfaces and curves have featured prominently in most of my own work over the past decade, and I thank the many amazing co-authors and students who have collaborated with me on this subject, including Tom Alberts, Morris Ang, Nathana\"el Berestycki, Manan Bhatia, Bertrand Duplantier, Ewain Gwynne, Nina Holden, Richard Kenyon, Sungwook Kim, Greg Lawler, Asad Lodhia, Oren Louidor, Jason Miller, Andrei Okounkov, Minjae Park, Yuval Peres, Joshua Pfeffer, R\'emi Rhodes, Oded Schramm, Nike Sun, Xin Sun, Vincent Vargas, Sam Watson, Menglu Wang, Wendelin Werner, David Wilson, Catherine Wolfram, Hao Wu and Pu Yu. We also thank Tom Alberts, Morris Ang, Nathanael Berestycki, Olivier Bernardi, Sky Cao, Nicolas Curien, Bertrand Duplantier, Ewain Gwynne, Jean-Fran{\c{c}}ois Le Gall, Gr\'egory Miermont, Jason Miller, Minjae Park, Guillaume Remy, R\'emi Rhodes, Steffen Rohde, Stanislav Smirnov, Yilin Wang and Wendelin Werner for their help in reading and improving the manuscript. And finally we thank Thomas Budzinski for allowing us to use his wonderful illustration.
\end{ack}

\begin{funding}
This work was partially supported by NSF Award: DMS 1712862.

\end{funding}

\bibliographystyle{emss}
\bibliography{ICM-randomsurface-capitals}

\end{document}